\documentclass[12pt]{amsart}
\usepackage{amsmath, amsthm, amssymb,cite,enumerate,slashed}
\usepackage{fullpage}
\usepackage{color}
\usepackage{hyperref}
\usepackage[english]{babel}
\usepackage{framed}

\allowdisplaybreaks

\newtheorem{thm}{Theorem}
\newtheorem{lem}{Lemma}[section]

\newtheorem{prop}[lem]{Proposition}
\newtheorem{cor}[lem]{Corollary}
\newtheorem{rmkn}{Remark}[thm]

\newtheorem{rmk}{Remark}[lem]

\newtheorem{que}[lem]{Question}

\newtheorem{defn}[lem]{Definition}
\numberwithin{equation}{section}

\newcommand{\R}{\mathbb{R}}
\newcommand{\D}{\mathcal{D}}

\newcommand{\Oc}{\mathcal{O}}

\newcommand{\eps}{\varepsilon}
\newcommand{\norm}[1]{\left\lVert#1\right\rVert}

\newcommand{\mcl}[1]{\mathcal{ #1 }}
\newcommand{\fm}[1]{\begin{align*} #1 \end{align*}}
\newcommand{\eq}[1]{\begin{equation}\begin{aligned} #1 \end{aligned}\end{equation}}
\newcommand{\lra}[1]{\langle #1 \rangle}
\newcommand{\abs}[1]{\left| #1 \right|}
\newcommand{\kh}[1]{\left( #1 \right)}

\newcommand{\wt}[1]{\widetilde{ #1 }}
\newcommand{\wh}[1]{\widehat{ #1 }}

\DeclareMathOperator{\supp}{supp}

\allowdisplaybreaks[1]

\begin{document}

\title{Timelike asymptotics for global solutions to a scalar quasilinear wave equation satisfying the weak null condition}

\author{Dongxiao Yu}
\address{Department of Mathematics, Vanderbilt University}
\thanks{}
\email{dongxiao.yu@vanderbilt.edu}


\begin{abstract}

We study the timelike asymptotics for global solutions to a scalar quasilinear wave equation satisfying the weak null condition. Given a global solution $u$ to the scalar wave equation with sufficiently small $C_c^\infty$ initial data, we derive an asymptotic formula for this global solution inside the light cone (i.e.\ for $|x|<t$). It involves the scattering data obtained in the author's asymptotic completeness result in \cite{MR4772266}. Using this asymptotic formula, we prove that $u$ must vanish under some decaying assumptions on $u$ or its scattering data, provided that the wave equation violates the null condition.

\end{abstract}

\maketitle
\tableofcontents\addtocontents{toc}{\protect\setcounter{tocdepth}{1}}

\section{Introduction}

This paper is devoted to the study of the long time dynamics for a scalar quasilinear wave equation in $\R^{1+3}$, of the form
\eq{\label{qwe} g^{\alpha\beta}(u)\partial_\alpha \partial_\beta u=0}
with small, smooth, and localized initial data
\eq{\label{init} (u,u_t)|_{t=0}=(\eps u_0,\eps u_1),\qquad \text{for fixed }u_0,u_1\in C_c^\infty(\R^3),\text{ and for } \eps\ll1.}
Here we use the Einstein summation convention, with the sum taken over $\alpha,\beta=0,1,2,3$ with  $\partial_0=\partial_t$, $\partial_i=\partial_{x_i}$, $i=1,2,3$. We assume that  the $g^{\alpha\beta}(u)$ are smooth functions of $u$, such that $g^{\alpha\beta}=g^{\beta\alpha}$ and $(g^{\alpha\beta}(0))=(m^{\alpha\beta})=\text{diag}(-1,1,1,1)$. Since we expect $|u|\ll1$, without loss of generality we also assume that $g^{00}\equiv -1$.

The equation \eqref{qwe} satisfies the weak null condition introduced by Lindblad and Rodnianski \cite{MR1994592}. Moreover, Lindblad \cite{MR2382144} proved that the Cauchy problem \eqref{qwe} and \eqref{init} has a unique global solution if the $\eps$ in \eqref{init} is sufficiently small. We refer to \cite{MR1177476,MR2003417} for earlier works on the global existence for \eqref{qwe}, and to \cite{keir2018weak} for an alternative proof of the global existence using the $r^p$-weighted energy method introduced in \cite{MR2730803}.

Recently, the author \cite{MR4232783,MR4315017,MR4772266} studied the asymptotic behaviors of global solutions to \eqref{qwe} near the light cone.   We first identified a new notion of asymptotic profile and an associated notion of scattering data for \eqref{qwe} by deriving a new reduced system called the geometric reduced system. Then, we proved the existence of the modified wave operators and the asymptotic completeness for~\eqref{qwe}. In summary, the modified scattering for \eqref{qwe} was established.

In this paper, we seek to study the asymptotic behaviors of global solutions to \eqref{qwe} in the interior of the light cone (i.e.\ for $|x|<t$). Our main result (Theorem \ref{mthm}) is an asymptotic formula for a global solution to \eqref{qwe}.  If $u=u(t,x)$ is a global solution to the Cauchy problem \eqref{qwe} and \eqref{init}, we have
\eq{\label{intro:asyfor}u(t,x)&=\frac{\eps}{2\pi}\int_{\mathbb{S}^2}\wh{A}(x\cdot\theta-t,\theta)\ dS_\theta+O(\lra{t-|x|}^{-2}t^{C\eps})}
whenever $t\geq e^{\delta/\eps}$ and $|x|<t$\footnote{Though the estimate \eqref{intro:asyfor} holds for all $|x|<t$, it is useful only when $|x|<t-t^{1/2+}$. See Remark \ref{mthm_rmk1.2}.}. Here $\delta\in(0,1)$ is a fixed parameter, and $\wh{A}=\wh{A}(q,\omega)$\footnote{\label{footnoteqomega}In \cite{MR4772266}, the variable $q$ corresponds to an optical function and $\omega$ corresponds to the usual angular coordinate. Note that we take $q\approx |x|-t$, so for a fixed $t$, we have $q\to-\infty$ as $|x|\to0$.} is the scattering data obtained from the asymptotic completeness result in \cite{MR4772266}.  We refer to Section \ref{sec1.2.3ac} for the definition of $\wh{A}$. Here one can view $\wh{A}$ as a nonlinear version of the Friedlander radiation field for $\partial u$.

As a corollary of \eqref{intro:asyfor}, if the equation \eqref{qwe} violates the null condition\footnote{Since $g^{00}\equiv -1$, we can equivalently assume that there exist $\alpha,\beta\in\{0,1,2,3\}$ such that $\frac{d}{du}g^{\alpha\beta}(u)|_{u=0}\neq 0$. }, 
then the scattering data $\wh{A}(q,\omega)$ of a nonzero global solution $u$ to \eqref{qwe} and \eqref{init} cannot decay arbitrarily fast in all directions (i.e.\ for all $\omega\in\mathbb{S}^2$) as $q\to-\infty$ (see Footnote \ref{footnoteqomega}). In this paper, we will prove that, if the null condition is violated, then the only global solution $u$ to \eqref{qwe} and \eqref{init} satisfying
\eq{\label{intro:adecay}|\wh{A}(q,\omega)|\lesssim \lra{q}^{-1-},\qquad \forall(q,\omega)\in\R\times\mathbb{S}^2}
is the zero solution as long as $\eps\ll1$. 
One can also obtain similar results for the pointwise decay rates of $u$ and $\partial u$. For example, assuming that
\fm{|(u_t-u_r)(t,(t-t^{1/2})\omega)|\lesssim \eps t^{-3/2-},\qquad \forall t\geq 1,\ \omega\in\mathbb{S}^2,}
we have \eqref{intro:adecay}. As a result, if the null condition is violated, under the pointwise bound for $\partial u$ above, we have $u\equiv 0$ as long as $\eps\ll1$.
We refer to Theorems \ref{sdltthm} and \ref{ltthm_weak} for the precise statements of these results.

\subsection{Background}\label{sec:lifespan:intro}
Let us consider  a generalization of the equation \eqref{qwe} in~$\R^{1+3}$
\begin{equation}\label{nlw}
\square u^I=F^I(u,\partial u,\partial^2u),\qquad I=1,2,\dots,M.
\end{equation}
The nonlinear terms are assumed to be smooth with the Taylor expansions
\begin{equation}\label{nonlinearity}F^I(u,\partial u,\partial^2u)=\sum a_{\alpha\beta,JK}^I\partial^\alpha u^J\partial^\beta u^K+O(|u|^3+|\partial u|^3+|\partial^2 u|^3).\end{equation}
The sum  is taken over all $1\leq J,K\leq M$ and all multiindices $\alpha,\beta$ with $|\alpha|\leq |\beta|\leq 2$, $|\beta|\geq1$ and $|\alpha|+|\beta|\leq 3$. Besides, the coefficients $a_{\alpha\beta,JK}^I$'s are all universal constants.

There have been several results on the lifespans of solutions to the Cauchy problem \eqref{nlw} with $C_c^\infty$ initial data of size $\eps\ll1$. For example, John \cite{MR808321,MR600571} proved that \eqref{nlw} does not necessarily have a global solution. There he proved finite time blowup results for $\Box u=u_t^2$ and  $\Box u=u_tu_{tt}$.
We also refer to \cite{MR1867312,MR3561670,MR3474069,MR2284927,MR4047642,MR3595936,MR3694013} for the blowup mechanisms for several quasilinear wave equations. 
For arbitrary nonlinearities, we expect to have almost global existence: the solution exists for all $t\in[0,e^{c/\eps}]$ where  $c>0$ is a small constant. We refer to \cite{MR1047332,MR0745325,MR0784477,MR0897781,MR1466700,MR2015331,MR1945285} for several cases where this expectation holds, but we remark that it has not been proved for the most general system \eqref{nlw}; see the discussions in \cite{MR4565922,MR4046191}.

To obtain global existence, one needs extra assumptions on the system \eqref{nlw}. It was proved by Klainerman \cite{MR837683,MR804771} and Christodoulou \cite{MR820070} that the \emph{null condition} is sufficient for global existence. We say that \eqref{nlw} satisfies the null condition if for each $1\leq I,J,K\leq M$ and for each $0\leq m\leq n\leq 2$ with $n\geq 1$ and $m+n\leq 3$, we have
\eq{\label{nullcond}A_{mn,JK}^I(\omega):=\sum_{|\alpha|=m,|\beta|=n}a_{\alpha\beta,JK}^I\widehat{\omega}^{\alpha}\widehat{\omega}^{\beta}=0,\qquad\text{whenever }\widehat{\omega}=(-1,\omega)\in\R\times\mathbb{S}^2.}
We also refer to \cite{MR2455195,MR2110542,MR1857981} and the references therein for extensions of this global existence result to systems of nonlinear wave equations with multiple wave speeds.

The null condition is not necessary for global existence. In \cite{MR2382144,MR2134337}, two examples that violate the null condition in general but admit global existence were presented. One example is the scalar equation \eqref{qwe} and the other is the Einstein vacuum equations in wave coordinates. These examples instead satisfy the \emph{weak null condition} introduced in \cite{MR1994592}. Note that the null condition implies the weak null condition and that both $\Box u=u_t^2$ and $\Box u=u_tu_{tt}$ violate the weak null condition. There is a conjecture that the weak null condition is sufficient for global existence\footnote{To be more precise, we have different versions of this conjecture. Some of them have been disproved, while others remain open. We refer to \cite[Section 1.1]{yu2022nontrivial} for a discussion on the current status of this conjecture.}, and we refer to \cite{keir2018weak,keir2019global,kadar2023small} for recent progress on it.

To define the weak null condition, we introduce a type of asymptotic equations derived by H\"ormander  \cite{MR0897781,MR1120284,MR1466700}. Suppose that $u$ is a global solution to \eqref{nlw} and we make the ansatz\begin{equation}\label{ansatz}
u^I(t,x)\approx \eps r^{-1}U^I(s,q,\omega),\hspace{1cm}r=|x|,\ \omega=x/r,\ s=\eps\ln t,\ q=r-t,\ 1\leq I\leq M.
\end{equation}
Assuming that $t=r\to\infty$, we substitute this ansatz into \eqref{nlw} and compare the coefficients of terms of order $\eps^2 t^{-2}$.  We thus obtain the following asymptotic PDE's (called \emph{H\"ormander's asymptotic equations} for \eqref{nlw}) for $U=(U^{I}(s,q,\omega))$:
\begin{equation}\label{asypde11}2\partial_s\partial_q U^I=\sum A_{mn,JK}^I(\omega)\partial_q^mU^J\partial_q^nU^K.\end{equation}
Here $A_{mn,JK}^I$ is defined by  \eqref{nullcond} and the sum in \eqref{asypde11} is taken over $0\leq m\leq n\leq 2$ with $n\geq 1$ and $m+n\leq 3$.
We say that  \eqref{nlw} satisfies the weak null condition if the corresponding \eqref{asypde11} has a global solution for all $s\geq 0$ and if the solution and all its derivatives grow at most exponentially in $s$, provided that the initial data decay sufficiently fast in $q$.  
For example, for the scalar equation \eqref{qwe}, the corresponding asymptotic equation \eqref{asypde11} becomes
\begin{equation}\label{asypde1}2\partial_s\partial_q U=G(\omega)U\partial_q^2U,\end{equation}
where \eq{\label{gomedefn}G(\omega):=g_0^{\alpha\beta}\widehat{\omega}_\alpha\widehat{\omega}_\beta,\qquad g_0^{\alpha\beta}=\frac{d}{du}g^{\alpha\beta}(u)|_{u=0},\ \widehat{\omega}=(-1,\omega)\in\R\times\mathbb{S}^2.}
One can prove the global existence and growth control for \eqref{asypde1}, so \eqref{qwe} satisfies the weak null condition.

\subsection{The modified scattering for \eqref{qwe}}

Following the global existence result in \cite{MR2382144}, we are interested in the long time dynamics for the equation \eqref{qwe}. Recently, the modified scattering for  \eqref{qwe} was established. In this subsection, we give an overview of the results in \cite{MR4772266,MR4232783,MR4315017} and discuss some related works on the long time dynamics for general quasilinear wave equations in $\R^{1+3}$.

\subsubsection{The geometric reduced system}
Instead of H\"ormander's asymptotic equation \eqref{asypde1}, we used a new system of asymptotic equations for \eqref{qwe} to study its long time dynamics in \cite{MR4772266,MR4232783,MR4315017}. It is called the \emph{geometric reduced system}, and it is expected to describe the asymptotic behaviors of global solutions to \eqref{qwe} more accurately than \eqref{asypde1} does.

To derive the geometric reduced system, we return to the ansatz \eqref{ansatz} (with $M=1$) used in the derivation of \eqref{asypde1}. We still make the ansatz
\eq{\label{ansatz:new}u(t,x)\approx \eps r^{-1}U(s,q,\omega),\qquad r=|x|,\ \omega=x/r,\ s=\eps\ln t,}
but now we let $q(t,x)$ be an optical function that is close to $r-t$ to some extent. By an optical function, we mean that $q$ solves the eikonal equation \eq{\label{intro:eik}g^{\alpha\beta}(u)\partial_\alpha q\partial_\beta q=0.}
Eikonal equations have been widely used in the study of nonlinear wave equations and the Einstein equations; see, e.g.,\ \cite{MR2003417,MR2382144,MR1316662,MR3638312,MR2178963,MR3402797,keir2018weak}. 
Our new ansatz is related to the geometry of the null cone with respect to the Lorentzian metric $(g_{\alpha\beta})$ that is the inverse of the $4\times 4$ matrix $(g^{\alpha\beta}(u))$. Since the geometric information from \eqref{qwe} is considered in our derivation, we expect to obtain a more accurate type of asymptotic equations for \eqref{qwe}.

When substituting the ansatz above into \eqref{qwe}, we obtain an equation involving $\partial q$. However, the optical function has no explicit formula because the eikonal equation is fully nonlinear. Alternatively, we define an auxiliary function $\mu=q_t-q_r$ and express $\partial q$ approximately in terms of $\mu$ and $U$ using the eikonal equation. By sending $t=r\to\infty$ and considering terms of order $\eps^2t^{-2}$ in \eqref{qwe}, we obtain the following reduced system for $(\mu,U)(s,q,\omega)$:
\eq{\label{asyeqn}\left\{\begin{array}{l}
\displaystyle \partial_s (\mu U_q)=0,\\
\displaystyle \partial_s\mu=\frac{1}{4}G(\omega)\mu^2U_q.\end{array}\right.}
Here $G(\omega)$ is defined by \eqref{gomedefn}. The system \eqref{asyeqn} is the geometric reduced system for \eqref{qwe}. By setting $(\mu,U_q)|_{s=0}(q,\omega)=(A_1,A_2)(q,\omega)$, we obtain an explicit solution to \eqref{asyeqn}:
\eq{\label{asyeqn:sol}\left\{\begin{array}{l}\displaystyle \mu(s,q,\omega)=A_1(q,\omega)\cdot \exp(\frac{1}{4}G(\omega)(A_1\cdot A_2)(q,\omega)s),\\[1em] \displaystyle  U_q(s,q,\omega)=A_2(q,\omega)\cdot \exp(-\frac{1}{4}G(\omega)(A_1\cdot A_2)(q,\omega)s).\end{array}\right.} In fact, the first equation in \eqref{asyeqn} yields $\mu U_q=A_1A_2$ for all $s\geq 0$, so the second equation in \eqref{asyeqn} is reduced to a linear ODE. From this perspective, the reduced system \eqref{asyeqn} is of a simpler form than \eqref{asypde1} which is a nonlinear PDE.

We  finally remark that the derivation above can be extended. Given an arbitrary system of quasilinear wave equations in $\R^{1+3}$, we can derive the corresponding geometric reduced system. See \cite[Chapter 2]{MR4315017} and \cite[(1.11)]{yu2022nontrivial}. 

\subsubsection{Existence of modified wave operators}\label{sec1.2.2emwo}

Making use of  \eqref{asyeqn}, we proved the existence of the modified wave operators for \eqref{qwe} in \cite{MR4232783}. Given an asymptotic profile, we seek to find a global solution to \eqref{qwe} such that this global solution matches the given asymptotic profile at infinite time. Such a problem is sometimes referred to as a backward scattering problem. To the author's knowledge, \cite{MR4232783} seems to be the only result on the modified wave operators for \eqref{qwe}. We refer to 
Lindblad-Schlue \cite{MR4591580,lindblad2023scattering} for backward scattering for semilinear wave equations satisfying the null condition or the weak null condition and some related models. Also see \cite{dafermos2013scattering,MR4660999,MR4299134,MR4746634,kadar2024scattering,MR3952698,MR4577274} and the references therein for similar results for the Einstein equations and other wave-type models.

Let us specify how the asymptotic profile is chosen in \cite{MR4232783}. Fix an arbitrary $A(q,\omega)\in C_c^\infty(\R\times\mathbb{S}^2)$ and suppose that $\supp A\subset[-R,R]\times\mathbb{S}^2$ for some constant $R>0$. Define $(\mu,U_q)(s,q,\omega)$ by \eqref{asyeqn:sol} with $(A_1,A_2)=(-2,A)$\footnote{\label{footnotegaugefreedom}We can take $\mu|_{s=0}\equiv-2$ because of the gauge freedom. In fact, one can prescribe $\mu|_{s=0}$ freely using the fact that $\overline{q}_t-\overline{q}_r=H_q\cdot (q_t-q_r)$ if $\overline{q}=H(q,\omega)$ for some function $H$.}. That is, we set 
\eq{\label{asyeqn:sol:2A}\left\{\begin{array}{l}\displaystyle \mu(s,q,\omega)=-2\exp(-\frac{1}{2}G(\omega)A(q,\omega)s),\\[1em] \displaystyle  U_q(s,q,\omega)=A(q,\omega)\cdot \exp(\frac{1}{2}G(\omega)A(q,\omega)s).\end{array}\right.}
We call this $A$ the \emph{scattering data} for the modified wave operator problem. To uniquely determine $U$, we add a boundary condition $\lim_{q\to-\infty}U(s,q,\omega)=0$\footnote{Since $A|_{q\leq -R}\equiv 0$, we have $U|_{q\leq -R}\equiv 0$. In a backward scattering problem, this property and the finite speed of propagation imply that the corresponding global solution to \eqref{qwe} vanishes for $r-t\leq -R$. This is a key property used in the proof in \cite{MR4232783}. }.
As a result, we obtain a solution $(\mu,U)(s,q,\omega)$ to the geometric reduced system \eqref{asyeqn} for all $s\geq 0$.

We now fix $\delta\in(0,1)$ which corresponds to the initial time when our reduced system starts to play a role\footnote{We will set $s=\eps\ln t-\delta$, so $s=0$ if and only if $t=e^{\delta/\eps}$.} and $\eps\ll1$ which is the size of the global solution we will construct.  We recover an approximate optical function $q(t,x)$ by solving the transport equation introduced in the derivation of \eqref{asyeqn} \fm{(q_t-q_r)(t,x)=\mu(\eps\ln t-\delta,q(t,x),x/|x|).} 
We choose the boundary condition so that $q\approx r-t$ to some extent.
Motivated by the new ansatz \eqref{ansatz:new}, we define
\eq{\label{asymprofile}\eps |x|^{-1}U(\eps\ln t-\delta,q(t,x),x/|x|)}
as the \emph{asymptotic profile} in this problem. One can show that \eqref{asymprofile} is an approximate solution to \eqref{qwe} whenever $|r-t|<t/2$ and $t>e^{\delta/\eps}$.  In \cite{MR4232783}, we proved that there exists a global solution $u$ to \eqref{qwe} for all $t\geq 0$ such that $u$ matches the asymptotic profile \eqref{asymprofile} at infinite time. See \cite[Theorem 1]{MR4232783}. For example, we have
\eq{\label{mwo:est}|u-\eps r^{-1}U|\lesssim \eps t^{-3/2+C\eps}\lra{r-t}^{1/2}}
for all $t\gtrsim1$ and $|x|\leq 5t/4$. Note that $u|_{|x|\leq t- R}\equiv 0$ and $U|_{q\leq -R}\equiv 0$, so the estimate \eqref{mwo:est} is trivial for $|x|\leq t-R$. Recall that the constant $R$ comes from the support of $A$. By comparing \eqref{mwo:est} with the pointwise decay $u=O(\eps t^{-1+C\eps})$, we conclude that the asymptotic profile \eqref{asymprofile} approximates the global solution $u$ well whenever $|r-t|\lesssim t^{1-}$.

In a recent work \cite{yu2022nontrivial}, the main theorem in \cite{MR4232783} was extended in two different ways. First, the assumption  $A\in C_c^\infty$ is relaxed. To prove the existence of the modified wave operators for \eqref{qwe}, we only need to assume that $A\in C^\infty$ and that \eq{\label{weakassumA}|\partial_q^a\partial_\omega^cA(q,\omega)|&\lesssim_{a,c} \lra{q}^{-\gamma_--a}\cdot 1_{q<0}+\lra{q}^{-\gamma_+-a}\cdot 1_{q\geq 0},\qquad \forall (q,\omega)\in\R\times\mathbb{S}^2,\ a,c\geq 0}where $\gamma_->2$ and $\gamma_+>1$ are two fixed parameters. See \cite[Corollary 8.1]{yu2022nontrivial}.
Moreover, the method in \cite{MR4232783} can be applied to a larger class of quasilinear wave equations in three space dimensions. In particular, one  obtains nontrivial global solutions to $\Box u=u_t^2$ and $\Box u=u_tu_{tt}$ for all $t\geq 0$, despite the finite time blowup results by John. See \cite[Corollaries 8.3 and 8.9]{yu2022nontrivial}. 

\subsubsection{Asymptotic completeness}\label{sec1.2.3ac}
In \cite{MR4772266}, we proved the asymptotic completeness for \eqref{qwe}. Given a global solution to \eqref{qwe} with initial data \eqref{init}, we seek to find the corresponding asymptotic profile. Before \cite{MR4772266}, Deng and Pusateri \cite{MR4078713} proved a partial scattering result for \eqref{qwe} using H\"ormander's asymptotic equation \eqref{asypde1}.  They applied the spacetime resonance method, and we refer to \cite{MR3084697,MR3202837} for some earlier applications of this method. We also refer to \cite{MR3638312,MR4746634,MR4660999,MR3936130,MR2928109,MR2400261,MR3952698} and the references therein for similar results for the Einstein equations and other wave-type models.

We now briefly explain the proof. More details are provided in Section \ref{secreview} below. Suppose that $u$ is the global solution to \eqref{qwe} with data \eqref{init} of size $\eps\ll1$. Let $R>0$ be a constant such that $\supp(u_0,u_1)\subset\{|x|\leq R\}$. We first construct a global optical function $q(t,x)$ by solving the eikonal equation \eqref{intro:eik} for all $|x|>t/2$ and $t>e^{\delta/\eps}$ where $\delta\in(0,1)$ is a small fixed parameter\footnote{Here we use $\{|x|>t/2,t>e^{\delta/\eps}\}$ for simplicity. In the actual proof, we use a region defined by \eqref{secreview:defnOmega}. }. The boundary condition is chosen so that $q\approx r-t$. In the proof, we apply both the method of characteristics and the method from \cite{MR1316662,MR2178963}. The key step is to carefully estimate the second fundamental form $(\chi_{ab})_{a,b=1,2}$ with the help of the Raychaudhuri equation.

Next, motivated by the ansatz \eqref{ansatz:new}, we define $(\mu,U):=(q_t-q_r,\eps^{-1}ru)$ . Both $\mu$ and $U$ are originally functions of $(t,x)$ for $|x|>t/2$ and $t>e^{\delta/\eps}$, and can be viewed as functions of $(s,q,\omega)$ via the inverse of the map
\fm{(t,x)\mapsto(\eps\ln t-\delta,q(t,x),x/|x|)=(s,q,\omega).}
This $(\mu,U)(s,q,\omega)$ is an approximate solution to the geometric reduced system \eqref{asyeqn}, and there is an exact solution $(\wt{\mu},\wt{U})(s,q,\omega)$\footnote{There is a different boundary condition $\lim_{q\to\infty}\wt{U}(s,q,\omega)=0$. In a forward problem, we have $u\equiv 0$ for $r-t\geq R$ by the finite speed of propagation.} to \eqref{asyeqn} matching $(\mu,U)(s,q,\omega)$ at infinite time; see \eqref{prevac:limita} and \eqref{prevac:defnwtmuU} below. Following the discussion in Section \ref{sec1.2.2emwo}, we hope that $(\wt{\mu},\wt{U})$ is also of the form \eqref{asyeqn:sol:2A}, so we can define the scattering data as the initial data of $\wt{U}_q$ at $s=0$. However, this is impossible since we may have $\wt{\mu}|_{s=0}\not\equiv -2$. To handle this issue, we use the gauge freedom (see Footnote \ref{footnotegaugefreedom}) to transform $(\wt{\mu},\wt{U})$ to an equivalent exact solution $(\wh{\mu},\wh{U})$ to \eqref{asyeqn} with $\wh{\mu}|_{s=0}\equiv-2$; see \eqref{prevac:defF} and \eqref{prevac:defhata} below. This new solution $(\wh{\mu},\wh{U})$ is of the form \eqref{asyeqn:sol:2A} with $A$ replaced by $\wh{A}=\wh{U}_q|_{s=0}$. This function $\wh{A}$ is called the \emph{scattering data} in the asymptotic completeness problem\footnote{This definition of scattering data is different from the one in \cite{MR4772266}. With this new definition, the gauge independence result can be stated in a much cleaner way.} and will play the same role as the scattering data $A$ in the modified wave operator problem. In general, one cannot show that $\wh{A}$ is $C^\infty$. Instead, for each integer $N\geq 0$, there exists an $\eps_N>0$ \emph{depending on $N$}, such that for all $\eps\in(0,\eps_N)$\footnote{Later we will simply write $\eps\ll_N1$. See Section \ref{prel:depeps}.},  we have $\wh{A}\in C^N$ and
\eq{\label{estwhA}|\partial_q^a\partial_\omega^c\wh{A}(q,\omega)|&\lesssim_{a,c} \lra{q}^{-1-a+C_{a,c}\eps},\qquad \forall (q,\omega)\in\R\times\mathbb{S}^2,\ a,c\geq 0,\ a+c\leq N.}
The constants $C_{a,c}$ and the implicit constants in $\lesssim_{a,c}$ are all uniform in $\eps$.
We also show a gauge independence result: for each fixed parameter $\delta$, the scattering data $\wh{A}$ is independent of the choice of the optical function; see Corollary~\ref{secreview:gauge:prop:cor}. 

Finally, we construct an approximate solution $\wt{u}$ to \eqref{qwe} from $(\wt{\mu},\wt{U})$\footnote{Let $\wh{u}$ be the approximate solution we obtain if we start from $(\wh{\mu},\wh{U})$ and follow the construction below. By \cite[Lemma 7.2]{MR4772266} and \eqref{wtuwhuequal}, we have $\wt{u}=\wh{u}$. We will thus not distinguish $\wt{u}$ from $\wh{u}$ in this paper.} by following the construction in \cite[Section 4]{MR4232783}. That is, we construct an approximate optical function $\wt{q}$ by solving $\wt{q}_t-\wt{q}_r=\wt{\mu}$, and  define $\wt{u}$ by \eqref{asymprofile} with $(U,q)$ replaced by $(\wt{U},\wt{q})$. Then, $\wt{u}$ is an approximate solution to \eqref{qwe} for $|x|>t/2$ and $t>e^{\delta/\eps}$. In addition, we prove that $\wt{u}$ offers a good approximation of $u$ at least near the light cone. See \cite[Theorem 1]{MR4772266} for a precise statement. For example, we have
\eq{\label{ac:est}|u-\wt{u}|\lesssim_\gamma \eps t^{-2+C\eps}\lra{r-t}} whenever $t>e^{\delta/\eps}$ and $|x|\geq t-t^\gamma$ where $\gamma\in(0,1)$ is a fixed parameter. Note that $u\equiv \wt{u}\equiv 0$ whenever $|x|\geq t+R$ by the finite speed of propagation, so the estimate \eqref{ac:est} is trivial for $|x|\geq t+R$. Recall that the constant $R$ comes from the support of the initial data \eqref{init}.

\subsubsection{More about the long time dynamics}
From \cite{MR4772266,MR4315017,MR4232783,yu2022nontrivial}, one may contend that the asymptotic behaviors of global solutions to \eqref{qwe} have been fully understood. This is however not true. For example, the following question is not answered in these papers.
\que{\label{que1}\rm In both the modified wave operator problem and the asymptotic completeness problem, can we give a full description of the asymptotic behaviors of the global solutions to \eqref{qwe}?}\rm

\bigskip

In the modified wave operator problem, starting from a scattering data $A$, we construct an asymptotic profile \eqref{asymprofile} and find the matching global solution $u$ to \eqref{qwe}. The estimate \eqref{mwo:est} suggests that $u\approx \eps r^{-1}U$ for $t\gg 1$ and $r\leq t+t^{1-}$. However, we do not have a good way to describe the asymptotic behaviors of $u$ outside the light cone. From \cite{yu2022nontrivial,MR4232783}, we only know that $|u|\lesssim \eps t^{-1/2+C\eps}r^{-1/2}$ whenever $t\gg1$ and  $r\geq 5t/4$. 

A similar issue arises in the forward problem. Let $u$ be the global solution to the Cauchy problem \eqref{qwe} and \eqref{init}. In the asymptotic completeness problem, we find the scattering data $\wh{A}$ and define an approximate solution $\wt{u}$ to \eqref{qwe}. The estimate \eqref{ac:est} suggests that $u\approx \wt{u}$ for $t\gg_\eps 1$ and $r\geq t-t^{1-}$. The asymptotic behaviors of $u$ inside the light cone (e.g.\ for $r<t/2$), however, are not discussed in \cite{MR4772266,MR4078713}. We only know from \cite{MR2382144} that $u=O(\eps t^{-1+C\eps})$.

Question \ref{que1} arises because of the limitations of H\"ormander's asymptotic equations and geometric reduced systems. Both these asymptotic equations are derived under the assumption that $r=t\to\infty$, so they are expected to approximate the original wave equations well only at null infinity. It is thus necessary to find a different way to describe the asymptotic behaviors of global solutions away from the light cone.

We now discuss some previous results related to Question \ref{que1}. For simplicity, our discussions will be restricted to the wave equations (e.g.\ \eqref{nlw}) in $\R^{1+3}$. For the backward problem, Lindblad and Schlue \cite{lindblad2023scattering}  constructed a global solution to $\Box u=0$ with prescribed radiation fields. Unlike the results in \cite{MR4591580}, in \cite[Theorem 1.5]{lindblad2023scattering} they gave a complete characterization of the asymptotics of the solution towards timelike and spacelike infinity.

In the forward problem, we first recall the strong Huygens principle: a solution to $\Box u=0$ with compactly supported data must vanish for $|r-t|\gtrsim 1$. This principle is however unstable \cite{MR0946226,mathisson1939probleme}. A related result is Price's law which was first conjectured in \cite{PhysRevD.5.2419}. It states that we have a $t^{-3}$ local uniform decay rate for linear waves on the Schwarzschild background, and we refer to \cite{MR2199010,MR2735767,MR2864787,MR3038715,MR2921169} for its proof. For a system \eqref{nlw} satisfying the null condition along with $C_c^\infty$ data of size $\eps\ll1$, Christodoulou \cite{MR820070} obtained a pointwise decay $u=O(\eps\lra{r-t}^{-1}\lra{r+t}^{-1})$. This bound is not optimal for a scalar semilinear or quasilinear wave equation satisfying the null condition with constant-coefficient null forms. In those cases, we have $u=O(\eps\lra{r-t}^{-2}\lra{r+t}^{-1})$; see \cite{dong2022generically,luk2024late,luk2024latenew}. On the other hand, if we pose the initial data on a hyperboloid instead of a time slice, and if we do not assume that the data is compactly supported, then the bound $u=O(\eps\lra{r-t}^{-1}\lra{r+t}^{-1})$ is optimal; see \cite{dong2022generically}. A similar result on the sharp pointwise decay for a semilinear wave equation with the null condition in nonstationary and asymptotically flat spacetimes can be found in \cite{looi2022global}. We also remark that a closely related topic is the \emph{late time tail} for solutions to wave equations on asymptotically flat spacetimes. In this direction, we refer to, e.g.,  \cite{luk2024late,MR3725885,MR4135420,MR4554671,MR4620915,looi2024asymptotic}, and the references therein (in particular, see Luk-Oh \cite[Section 1.3]{luk2024late}).

Our final remark is that all the references listed in the previous two paragraphs do not apply to \eqref{qwe}. This is due to its weak null structure.

\subsection{The main theorems}

In this paper, we will partially answer Question \ref{que1} by providing a description of the asymptotic behaviors of a global solution $u$ to the forward Cauchy problem \eqref{qwe} and \eqref{init} inside the light cone. We will present an asymptotic formula for $Z^Iu$ for each multiindex $I$. Here  $Z$ denotes one of the commuting vector fields: translations $\partial_\alpha$, scaling $S=t\partial_t+r\partial_r$, rotations $\Omega_{ij}=x_i\partial_j-x_j\partial_i$, and Lorentz boosts $\Omega_{0i}=x_i\partial_t+t\partial_i$. Besides,  $Z^I$ denotes a product of $|I|$ commuting vector fields for each multiindex $I$.

Fix $(u_0,u_1)\in C_c^\infty(\R^3)$. From \cite{MR2382144}, we know that the Cauchy problem \eqref{qwe} and \eqref{init} admits a unique global $C^\infty$ solution $u$ for all $t\geq 0$ as long as $\eps\ll1$. In Section \ref{sec1.2.3ac}, for a fixed $\delta\in(0,1)$, we explained how to construct the scattering data $\wh{A}$ for the global solution~$u$. For each multiindex $I$, we define $A_I=A_I(q,\omega)$ inductively by $A_0=-2\wh{A}$ and
\eq{\label{s1:AIdefn}
A_I(q,\omega)=\left\{
\begin{array}{ll}
    q\partial_qA_{I'}, & Z^I=SZ^{I'};  \\
    (\omega_i\partial_{\omega_j}-\omega_j\partial_{\omega_i})A_{I'}, & Z^I=\Omega_{ij}Z^{I'},\ 1\leq i<j\leq 3;\\
     -q\omega_i\partial_qA_{I'}+\partial_{\omega_i}A_{I'}-2\omega_iA_{I'}, & Z^I=\Omega_{0i}Z^{I'},\ 1\leq i\leq 3;\\
   0, & Z^I=\partial Z^{I'};
\end{array}\right.
}
This definition comes from Proposition \ref{prop:trzwtu}. We also recall from Section \ref{sec1.2.3ac} that for each integer $N\geq 0$, we have $\wh{A}\in C^N$ as long as $\eps\ll_N1$ and the estimate \eqref{estwhA}. Thus, we can define $A_I$ as long as $\eps\ll_{|I|}1$. Moreover, the estimate \eqref{estwhA} holds with $\wh{A}$ replaced by $A_I$.

We now explain the meanings of $A_I$ and \eqref{s1:AIdefn}. By Proposition \ref{prop:trzwtu}, we have \fm{(\partial_t-\partial_r)Z^I\wt{u}= \eps r^{-1}A_I(r-t,\omega)+\text{an error term},\qquad \forall t\gg_{\delta,\eps}1,\ t/2<r<t+R.} where $\wt{u}$ is the approximate solution constructed from the scattering data $\wh{A}$ in Section~\ref{sec1.2.3ac}. Moreover, given two multiindices $I$ such that $Z^I=ZZ^{I'}$, we seek to prove
\fm{Z(\eps r^{-1}A_{I'}(r-t,\omega))=\eps r^{-1}A_{I}(r-t,\omega)+\text{an error term}.}
Apply the chain rule to compute the left hand side. This gives us \eqref{s1:AIdefn}.

We are ready to state the first main theorem.

\thm{\label{mthm} Fix $(u_0,u_1)\in C_c^\infty(\R^3)$. Let $u$ be the global $C^\infty$ solution to \eqref{qwe} with initial data \eqref{init} for $\eps\ll1$. Fix a small constant $\delta\in(0,1)$, and let $\wh{A}$ be the scattering data corresponding to $u$. Then, for each multiindex $I$, as long as $0<\eps\ll_{\delta,|I|}1$, we have
\eq{\label{mthmcl}Z^Iu(t,x)&=-\frac{\eps}{4\pi}\int_{\mathbb{S}^2}A_I(x\cdot \theta-t,\theta)\ dS_\theta+O_{|I|}(\lra{t-|x|}^{-2}t^{C_{|I|}\eps}),\qquad \forall t\geq e^{\delta/\eps}, |x|<t.}
The function $A_I$ is defined inductively by $A_0=-2\wh{A}$ and \eqref{s1:AIdefn}.}
\rm

\rmkn{\rm The scattering data $\wh{A}$ describes the asymptotic behaviors of a global solution $u$ to \eqref{qwe} and \eqref{init} near the light cone. 
Thus, the formula \eqref{mthmcl} connects the asymptotic behaviors of $u$ in the interior of the light cone with the asymptotics of $u$ at null infinity. The reason for this will be clear later when we discuss the proofs of the main theorems.}

\rmkn{\rm \label{mthm_rmk1.2} The estimate \eqref{mthmcl} holds whenever $t\geq e^{\delta/\eps}$ and $|x|<t$. However, we recall from \cite{MR2382144} that $Z^Iu=O(\eps t^{-1+C\eps})$ for each multiindex $I$ as long as $\eps\ll_{|I|}1$. Thus, it provides a good approximation for $Z^Iu$ only if the remainder term $O(\lra{t-|x|}^{-2}t^{C\eps})$ is far less than $\eps t^{-1+C\eps}$. For example, in the region where $|x|<t-t^{\gamma}$ for a fixed constant $\gamma>1/2$, we have $\lra{t-|x|}^{-2}t^{C\eps}\leq t^{-2\gamma+C\eps}\ll\eps t^{-1+C\eps}$. 
}

\rmkn{\label{rmk:realasymp}\rm
It is natural to ask whether \eqref{mthmcl} is a \emph{real} asymptotic formula for $Z^Iu$. For simplicity, we only discuss the case when $|I|=0$ here. In Theorem \ref{mthm}, we do not exclude the possibility that \eq{\label{rmk:realasymp:est}|\int_{\mathbb{S}^2}\wh{A}(x\cdot\theta-t,\theta)\ dS_\theta|\lesssim \lra{|x|-t}^{-2}t^{C\eps},\qquad \forall t\geq e^{\delta/\eps},\ |x|<t.}
In fact, if the equation \eqref{qwe} satisfies the null condition, then the estimate \eqref{rmk:realasymp:est} always holds; see Remark \ref{sdltthm_rmk} below. In this case, the estimate \eqref{mthmcl} reduces to a pointwise bound $|u|\lesssim \lra{r-t}^{-2}t^{C\eps}$ and is therefore not accurate enough.

One can ask whether this is also the case when the null condition is violated. The answer is no. Later we will prove that, if $u$ is a nonzero global solution to \eqref{qwe} (violating the null condition) and \eqref{init} with $\eps\ll1$, then \eqref{rmk:realasymp:est} fails for a sequence of points $\{(t_n,x_n)\}$ with $t_n\uparrow\infty$ and $|x_n|<t_n$. This is a corollary of Theorem \ref{ltthm_weak} and we refer to Remark \ref{rmk:ltthm_weak:realasymp} below.
}

\rmkn{\label{mthm_rmk1.3}\rm By \eqref{s1:AIdefn}, we have $A_I\equiv 0$ whenever $Z^I$ contains a translation $\partial$. This does not indicate that we cannot approximate terms like $\partial Z^Iu$. Given a multiindex $I$, we apply \cite[(3.1)]{MR2382144}:
\fm{\partial_t&=\frac{tS-\sum_{i=1}^3x_i\Omega_{0i}}{(t-r)(t+r)},\qquad \partial_j=\frac{\sum_{i=1}^3 x_i\Omega_{ij}+t\omega_j\Omega_{0j}-x_iS}{(t-r)(t+r)}}
to replace each translation in $Z^I$ with a finite sum of $S,\Omega_{ij},\Omega_{0i}$ multiplied by an explicit function of $(t,x)$. By Leibniz's rule, one can write \eq{\label{mthm_rmk1.3_formula}Z^Iu&=\sum_{|J|\leq |I|,\ k_J=0}c_{I,J}(t,x)Z^Ju}
where $k_J$ denotes the number of translations in $Z^J$ and $c_{I,J}(t,x)$ can be computed explicitly from the coefficients in \cite[(3.1)]{MR2382144}.   Combining this expression and \eqref{mthmcl}, we thus obtain a more accurate asymptotic formula for $Z^Iu$ when $Z^I$ contains translations. 

We also remark that $c_{I,J}=O((t-r)^{-k_I})$ whenever $r<t-1$ and $t\gg1$. In other words, by setting $\wt{c}_{I,J}=(t-r)^{k_I}c_{I,J}=O(1)$, we obtain
\eq{\label{mthmclext}(t-r)^{k_I}Z^Iu(t,x)&=-\frac{\eps}{4\pi}\sum_{|J|\leq |I|,\ k_J=0}\wt{c}_{I,J}(t,x)\int_{\mathbb{S}^2}A_J(x\cdot \theta-t,\theta)\ dS_\theta+O_{|I|}(\lra{t-|x|}^{-2}t^{C_{|I|}\eps})}
for all $t\geq e^{\delta/\eps}$ and $|x|< t-1$. When $Z^I$ contains translations, the estimate \eqref{mthmclext} is a more accurate asymptotic formula for $Z^Iu$. For simplicity, we do not present the formulas for $\wt{c}_{I,J}$ here, but we emphasize that they can be computed explicitly.
}

\rmkn{\rm We now have a full description of the asymptotic behaviors of the global solution $u$ to \eqref{qwe} and \eqref{init}. First, Theorem \ref{mthm} provides a good approximation for $Z^Iu$ for $t\gg_{\delta,\eps}1$ and $|x|<t-t^{1/2+}$. Moreover, we obtain from \cite{MR4772266} an approximation for $u$ whenever $t\gg_{\delta,\eps}1$ and $|x|>t-t^{1-}$.  For each multiindex $I$, we have 
\fm{Z^Iu(t,x)=Z^I\wt{u}(t,x)+O_{\gamma,|I|}(\eps t^{-2+C_{I}\eps}\lra{r-t})}
whenever $t\geq 2e^{\delta/\eps}$ and $|r-t|\lesssim t^\gamma$ for each fixed $\gamma\in(0,1)$. Also see \eqref{prevac:uapp} in Section \ref{secreview:secapp}.
Here $\wt{u}$ is defined in both Section \ref{sec1.2.3ac} and Section \ref{secreview:secapp}. We also recall that $u\equiv 0$ whenever $r-t\geq R$ where $R>0$ is a constant such that the initial data \eqref{init} vanishes for $|x|\geq R$. In other words, for all $(t,x)$ with $t\gg_{\delta,\eps}1$, we have at least one way to approximate $Z^Iu(t,x)$.
}

\rmkn{\label{mthm_rmk1.5}\rm The proof of \eqref{mthmcl} is robust and does not only work for \eqref{qwe}. Let $\phi=\phi(t,x)$ be a $C^2$ function defined for all $t>0$ and $|x|<t$. If $\Box \phi$ satisfies a good pointwise bound, and if we have a good way to approximate $(-\partial_t+\partial_r)(r\phi)$ near the light cone $r=t$ for all $t\gg1$, then we can prove an asymptotic formula for $\phi$ similar to \eqref{mthmcl}. In fact,  
approximately, one can apply  Proposition \ref{prop:inhomwe} to show that \eq{\label{rmk1.5f1}\phi(t,x)&\approx-\frac{1}{4\pi}\int_{\mathbb{S}^2}\mcl{A}_0(x\cdot\theta-t,\theta)\ dS_\theta,\qquad  \forall |x|<t,\ t\gg1.}
Here $\mcl{A}_0$ is the radiation field for $\phi_t-\phi_r$ in the sense that
\fm{\phi_t-\phi_r\approx r^{-1}\mcl{A}_0(r-t,\omega),\qquad \text{whenever }|r-t|\ll t.}
For a nonrigorous derivation of this approximate identity, we refer to the discussion above \eqref{intro:phifor} in Section \ref{sec:structure}. The form of the error term in \eqref{rmk1.5f1} depends on the pointwise bounds for $\phi$ and $\Box\phi$, and we will not discuss it in detail here for simplicity.}  

\rmkn{\rm We now discuss how \eqref{rmk1.5f1} is related to \cite[Theorem 1.5]{lindblad2023scattering}. Recall that the authors of \cite{lindblad2023scattering} constructed a global solution $\phi$ to $\Box \phi=0$ with a prescribed radiation field and then described the interior and exterior asymptotics of $\phi$. In particular, they mentioned that $\phi$ has the interior homogeneous asymptotics \cite[(1.11)]{lindblad2023scattering}:
\eq{\label{rmk1.5f2}\phi(t,x)\approx\frac{1}{4\pi}\int_{\mathbb{S}^2}\frac{N(\theta)}{t-x\cdot \theta}\ dS_\theta,\qquad \text{for } |x|/t<1\text{ fixed},\ t\to\infty.}
Interestingly, given the radiation field of $\phi_t$ in \cite[Theorem 1.5]{lindblad2023scattering}, one can nonrigorously derive \eqref{rmk1.5f2} to the leading order by \eqref{rmk1.5f1}. In \cite[Remark 1.3]{lindblad2023scattering}, it was mentioned that $\phi_t$ has the radiation field
\fm{\phi_t\approx r^{-1}\mcl{G}_0(r-t,\omega),\qquad \text{with }\mcl{G}_0=N_{01}(\omega)q\lra{q}^{-2}-\partial_q\mcl{F}_{01}(q,\omega).}
Both $N_{01}$ and $\mcl{F}_{01}$ are from the statement of \cite[Theorem 1.5]{lindblad2023scattering}. When $q<-1$, we have $\partial_q\mcl{F}_{01}=O(\lra{q}^{-2})$, which is a lower order term compared to $N_{01}(\omega)q\lra{q}^{-2}$. Since $\phi_t+\phi_r$ is a tangential derivative, we expect $\phi_r\approx-\phi_t$, so $\phi_t-\phi_r\approx 2\mcl{G}_0\approx 2N_{01}q^{-1}$ whenever $q<0$ and $|q|\gg1$. Plug this into \eqref{rmk1.5f1}, and we obtain
\fm{\phi\approx -\frac{1}{2\pi}\int_{\mathbb{S}^2}\frac{N_{01}(\theta)}{x\cdot\theta-t}\ dS_\theta.}}

\rmkn{\rm 
The formula \eqref{mthmcl} can also be extended to a scalar quasilinear wave equation satisfying the null condition:
\eq{\label{mthm:rmk:qwenc}g^{\alpha\beta}(u,\partial u)\partial_\alpha\partial_\beta u=f(u,\partial u).}
The Cauchy problem \eqref{mthm:rmk:qwenc} and \eqref{init} admits a global solution $u$. One can also show that the following limit exists and is a continuous function:
\fm{B(q,\omega):=\lim_{r\to\infty}(-\partial_t+\partial_r)(ru)(r-q,r\omega).}
By Proposition \ref{prop:inhomwe}, we can show not only \eqref{rmk1.5f1} but also a more precise estimate
\eq{\label{mthm_rmk1.5_2}u(t,x)&=\frac{1}{4\pi}\int_{\mathbb{S}^2}B(x\cdot\theta-t,\theta)\ dS_\theta+O(t^{-1}\lra{t-|x|}^{-3}),\qquad \forall t>1,|x|<t.}
This estimate is weaker than the results in \cite{dong2022generically,luk2024late,luk2024latenew}, but it provides a potentially different way to prove the results in those papers. To recover, e.g., \cite[(1.12)]{dong2022generically}, one needs an asymptotic formula for $B$ as follows:
\eq{\label{asymB}B(q,\omega)=c_0q^{-2}+O(|q|^{-3}),\qquad \forall q<-1, \text{for some constant }c_0.}
If it is true, then we have
\fm{u(t,x)&=\frac{1}{4\pi}\int_{\mathbb{S}^2}c_0(x\cdot\theta-t)^{-2}+O(|x\cdot\theta-t|^{-3})\ dS_\theta+O(t^{-1}\lra{t-|x|}^{-3})\\
&= c_0(t^2-|x|^2)^{-1}+O( t^{-1}\lra{t-|x|}^{-2})} for all $t>1$ and $|x|<t-1$.  If one can furthermore prove that $c_0=0$, then we obtain the better pointwise decay rate $u=O(t^{-1}\lra{|x|-t}^{-2})$ in \cite{dong2022generically,luk2024late,luk2024latenew}. It is however unclear to the author whether one can show the asymptotics \eqref{asymB} and $c_0=0$, so we do not claim in this paper that we have given different proofs of those results.
}

\rm 

\bigskip

By \eqref{mthmcl}, we can show several interesting properties of global solutions to \eqref{qwe} when the null condition is violated. All but one of them fail if the null condition is satisfied. 

Recall that the null condition was defined for a general system \eqref{nlw} in Section \ref{sec:lifespan:intro}. By \eqref{asypde1}, the equation \eqref{qwe} satisfies the null condition if and only if $G(\omega)\equiv0$ on $\mathbb{S}^2$, with  $G(\omega)$ defined by \eqref{gomedefn}. Since $g^{00}\equiv -1$, we claim that $G(\omega)\equiv 0$ on $\mathbb{S}^2$ if and only if $g^{\alpha\beta}_0=0$ for all $\alpha,\beta=0,1,2,3$.  In fact, if  $G(\omega)\equiv 0$ on $\mathbb{S}^2$, we  have a decomposition\fm{(g^{\alpha\beta}_0)=g^{00}_0\text{diag}(1,-1,-1,-1)+(\text{an asymmetric matrix}).} 
Since $(g^{\alpha\beta})$ is symmetric, so is $(g^{\alpha\beta}_0)$. Thus the asymmetric part in this decomposition vanishes. And since $g^{00}\equiv -1$, we have $g^{00}_0=0$ and thus $g^{\alpha\beta}_0=0$.

We first state a theorem related to the scattering data.

\thm{\label{sdltthm}Suppose that the equation \eqref{qwe} \emph{violates} the null condition. Fix $(u_0,u_1)\in C_c^\infty(\R^3)$. Let $u$ be the global solution to \eqref{qwe} with initial data \eqref{init} for $\eps\ll1$. Fix a small $\delta\in(0,1)$ and let $\wh{A}$ be the corresponding scattering data. 

Moreover, suppose that one of the following assumptions holds:
\begin{enumerate}[\rm a)]
    \item We have $\min_{(q,\omega)\in\mathbb{R}\times\mathbb{S}^2}\wh{A}(q,\omega)\geq 0$ or $\max_{(q,\omega)\in\mathbb{R}\times\mathbb{S}^2}\wh{A}(q,\omega)\leq 0$; 
    \item We have $\min_{ \omega\in \mathbb{S}^2}G(\omega)\geq 0$ and $\min\{\wh{A},0\}\in L^1_{q,\omega}(\R\times\mathbb{S}^2)$, or we have $\max_{ \omega\in \mathbb{S}^2}G(\omega)\leq 0$ and $\max\{\wh{A},0\}\in L^1_{q,\omega}(\R\times\mathbb{S}^2)$;
    \item Let $C_0>0$ be a constant uniform in all $\eps\ll1$ such that \fm{C_0\geq \sup_{(q,\omega)\in\R\times\mathbb{S}^2}|\frac{1}{2} G(\omega)\wh{A}(q,\omega)|.}
    Such a $C_0$ always exists because of \eqref{estwhA}.
    Suppose that there exists a  constant $B_0>C_0$ that is also uniform in all $\eps\ll1$, such that  \fm{\lim_{q\to-\infty}\lra{q}^{1+B_0\eps}\sup_{\omega\in\mathbb{S}^2}|\wh{A}(q,\omega)|=0.}
\end{enumerate}
Then, as long as $\eps\ll1$ (depending on the constant $B_0$ in part c)),  we have $\wh{A}\equiv 0$ and thus $u\equiv 0$.

}

\rmkn{\label{sdltthm_rmk}\rm The assumption a) implies $\wh{A}\equiv 0$ and $u\equiv 0$ no matter whether \eqref{qwe} satisfies the null condition. We refer to Section \ref{vanishingA:assumpa} and Footnote \ref{FN12}. However, for b) and c), it is necessary to assume that the null condition is violated. Suppose that \eqref{qwe} satisfies the null condition, and let $u$ be a global solution to \eqref{qwe} and \eqref{init}. One can check that 
\eq{\label{sdltthm_rmk_nc}\sum_{|I|\leq N}|Z^I u(t,x)|\lesssim\eps t^{-1}\lra{|x|-t}^{-1},\qquad \text{as long as }\eps\ll_N1.}
The case $N=0$ has been proved in Christodoulou \cite{MR820070}. For $N>0$, one can apply the estimates in \cite[Section 6.6]{MR1466700} and \cite[Proposition 3.1]{MR2244605}.
Following the proofs in \cite[Section 5.4]{MR4772266}, from these pointwise bounds we obtain $|\wh{A}(q,\omega)|\lesssim \lra{q}^{-2}$. In other words, the assumptions b) and c) hold in general.}

\rmkn{\label{sdltthm_rmk2}\rm Let us compare the assumption c) with the pointwise bound \eqref{estwhA} for $\wh{A}$. By  \eqref{estwhA}, we have $|\wh{A}(q,\omega)|\lesssim\lra{q}^{-1+C\eps}$ for all $(q,\omega)\in\R\times\mathbb{S}^2$. Part c) of Theorem \ref{sdltthm} indicates that if the null condition is violated and if $|\wh{A}(q,\omega)|\lesssim\lra{q}^{-1-B_0\eps}$ everywhere for a large constant $B_0>1$, then $u\equiv 0$ as long as $\eps\ll1$. In other words, if the null condition is violated, the scattering data cannot decay arbitrarily fast in all directions as $q\to-\infty$. Besides, the exponent $-1+C\eps$ in $|\wh{A}(q,\omega)|\lesssim\lra{q}^{-1+C\eps}$ is almost sharp.

We also compare the assumption c) with the assumption \eqref{weakassumA} in the modified wave operator problem. There we assume that $|\wh{A}|\lesssim \lra{q}^{-2-}$, so the global solutions constructed in \cite{yu2022nontrivial,MR4232783} cannot have compactly supported initial data if the null condition is violated. It is unclear whether one can further relax the assumption \eqref{weakassumA} in the modified wave operator problem. For example, can we assume $\gamma_->1$ instead of $\gamma_->2$?

The scattering $\wh{A}$ satisfies an extra decaying property. By Proposition~\ref{prop:decaysphmeans}, we have
\fm{|\int_{\mathbb{S}^2}\wh{U}(\eps\ln t-\delta,-2t,\omega )dS_\omega|\lesssim\eps^{-1}t^{-1/3+C\eps},\qquad \forall t\geq 2e^{\delta/\eps}}
where (recall that $R>0$ is chosen so that the data \eqref{init} vanishes for $|x|\geq R$)
\fm{\wh{U}(s,q,\omega)=-\int_{q}^R\wh{A}(p,\omega)\exp(\frac{1}{2}G(\omega)\wh{A}(p,\omega) s)\ dp.}
In contrast, from $\wh{A}=O(\lra{q}^{-1+C\eps})$, we only have $\wh{U}(s,q,\omega)=O(\eps^{-1}\lra{q}^{C\eps}e^{Cs})$ and $\wh{U}(\eps\ln t-\delta,-2t,\omega)=O(\eps^{-1}t^{C\eps})$. Because of part c), we cannot improve this bound for $\wh{U}$ if we only rely on the pointwise bounds for $\wh{A}$. Thus, to relax the assumption \eqref{weakassumA} in the modified wave operator problem, we may need to also take such a decaying property into account.
}
\rm
\bigskip

We now state a corollary of Theorem \ref{sdltthm}. Here we focus on the pointwise decay for the solutions and their first derivatives.

\thm{\label{ltthm_weak} Suppose that  \eqref{qwe} \emph{violates} the null condition. Fix $(u_0,u_1)\in C_c^\infty(\R^3)$. Let $u$ be the global solution to \eqref{qwe} with initial data \eqref{init} for $\eps\ll1$. Suppose that there exist constants $B_0,B_1>0$  depending on $u_0,u_1$ (but not on $\eps$), such that one of the following assumptions holds:
\begin{enumerate}[\rm i)]
    \item Fix  $0<\nu_0<1$. For all $t\geq 1$ and $\omega\in\mathbb{S}^2$, we have \fm{|(u_t-u_r)(t,(t-t^{\nu_0})\omega)|\lesssim \eps t^{-1-\nu_0-\nu_0B_0\eps};}
    \item  Fix  $0<\nu_0<1$. For all $t\geq 2$ and $x\in\R^3$ with $|x|\in[t-2t^{\nu_0},t-t^{\nu_0}/2]$, we have   \fm{ |u(t,x) |\lesssim\eps t^{-1-(\nu_0 B_0+B_1)\eps};}
    \item Suppose that $u$ is a spherically symmetric solution  with respect to $x$ and that \fm{|u(t,0)|\lesssim\eps t^{-1-B_0\eps} \qquad \forall t\geq 1.}
\end{enumerate}

Then, as long as $\eps\ll1$ which may depend on the constants in each of the assumptions, we have  $u
\equiv0$.}
\rmkn{\rm Note that the assumptions i) and ii) are related to the pointwise decays of $u$ and $\partial u$ for $|r-t|\lesssim t^{\nu_0}$. Recall from \cite{MR2382144} that $|\partial u|\lesssim \eps t^{-1}\lra{r-t}^{-1+C\eps}$ and $|u|\lesssim \eps t^{-1+C\eps}$. Thus, both the exponents $-1+C\eps$ in these estimates are almost sharp if the null condition is violated.

If the null condition is satisfied, then as shown in \eqref{sdltthm_rmk_nc} in Remark \ref{sdltthm_rmk}, we have $|\partial u|\lesssim \eps t^{-1}\lra{r-t}^{-2}$ and $|u|\lesssim \eps t^{-1}\lra{r-t}^{-1}$. That is, the assumptions i) and ii) hold in general.}

\rmkn{\rm\label{rmk:ltthm_weak:realasymp} We now prove the fact stated in Remark \ref{rmk:realasymp}. Suppose that the null condition is violated. Let $u$ be a nonzero global solution to \eqref{qwe} and \eqref{init} with $\eps\ll1$. If we take $\nu_0=3/4$ (or any number in $(1/2,1)$) in part ii) of Theorem \ref{ltthm_weak}, we obtain a sequence of points $\{(t_n,x_n)\}$ with \fm{t_n\uparrow\infty,\qquad 0<t_n-|x_n|\sim t^{3/4},\qquad |u(t_n,x_n)|\geq \eps t^{-1-C_{B_0,B_1}\eps}.} We also have
\fm{\lra{|x_n|-t_n}^{-2}t_n^{C\eps}\lesssim t_n^{-3/2+C\eps}\ll \eps t_n^{-1-C_{B_0,B_1}\eps},\qquad \forall t\geq e^{\delta/\eps}.}
In other words, the estimate \eqref{rmk:realasymp:est} must fail for this sequence of points.
}

\rmkn{\rm  Part iii) of Theorem \ref{ltthm_weak} suggests that, if the null condition is violated, and if the solution is spherically symmetric with respect to $x$, then $u(t,0)$ cannot decay arbitrarily fast as $t\to\infty$. Unlike the assumptions i) and ii), here we make an assumption on the decay rate of $u$ along the time-axis which is far inside the light cone. In Proposition \ref{ltthm} below), we can replace the spherical symmetry condition with an assumption on the decay rate of $\partial\Omega_{ij}u$ near the light cone.

It is however unclear to the author whether the spherical symmetry condition and the assumption on the decay rate of $\partial\Omega_{ij}u$  can be completely removed. 
For example, assuming
\fm{|u(t,x)|&\lesssim \eps t^{-1-},\qquad \forall t\geq 1,\ |x|\leq t/2,}
can we show $u\equiv 0$ for $\eps\ll1$?}

\rm

\subsection{Structure of this paper}
\label{sec:structure}
We now explain how this paper is organized. We also invite our readers to check the beginning of each of Sections \ref{secreview}--\ref{sec:corofmthm} for an overview of that specific section.

In Section \ref{sec:preliminaries}, we introduce the notation and several preliminary lemmas. We emphasize that a special convention involving $\eps$ is made in Section \ref{prel:depeps}.

In Section \ref{secreview}, we present an overview of the asymptotic completeness result for \eqref{qwe}. We follow the discussion in Section \ref{sec1.2.3ac} and provide more details. We will also prove several results that did not appear in \cite{MR4772266}. For example, in Proposition \ref{prop:trzwtu}, we provide an asymptotic formula for higher derivatives of $\wh{u}$ which was constructed from our scattering data in \cite[Section 7]{MR4772266}. In Section \ref{secreview:gauage}, we provide a new gauge independence result which is stated in a much cleaner way than \cite[Proposition 6.1]{MR4772266}. 

In Section \ref{sec:generalinhom}, we study a general inhomogeneous wave equation $\Box \phi=F$ inside the light cone $|x|=t$. Assuming that \eq{\label{sec:intro:prop4.1:asu}|F(t,x)|\lesssim \lra{t}^{-\gamma_1}\lra{r-t}^{-\gamma_2},\qquad \text{in }\Oc=\{t>0,|x|<t\},\ \text{with }\gamma_1,\gamma_2>0,} we derive an asymptotic formula for $\phi$: whenever $(t,x)\in\Oc$ and whenever  $T>2t$, we have \eq{\label{sec:intro:prop4.1}\phi(t,x)&=\frac{1}{4\pi}\int_{\mathbb{S}^2}\Phi(T,x-(T-t)\theta)\ dS_\theta+\text{(a remainder term)},\quad \text{with }\Phi:=(-\partial_t+\partial_r)(|x|\phi).} See Proposition \ref{prop:inhomwe}. In its proof, we view $\phi$ as a solution to the backward Cauchy problem $\Box \phi=F$ with initial data assigned at $t=T$, so we can express $\phi(t,x)$ explicitly as a sum of an integral of $F$ and an integral of $(\phi,\phi_t)|_{t=T}$. We use \eqref{sec:intro:prop4.1:asu} to control the integral of $F$, so $\phi$ is approximated well by the solution to $\Box w=0$ with data $(\phi,\phi_t)|_{t=T}$ at $t=T$. Direct computations yield \eqref{sec:intro:prop4.1}.

The formula \eqref{sec:intro:prop4.1} is useful in this paper because one can relate the integral of $\Phi$ with the asymptotics of $\phi$ at null infinity. In fact, we will prove in Lemma \ref{pfmthm:lem:basic} that
\fm{\lim_{T\to\infty}(|x-(T-t)\theta|-T,\frac{x-(T-t)\theta}{|x-(T-t)\theta|})=(-t-x\cdot \theta,-\theta),\qquad \forall t>|x|>0,\ \theta\in\mathbb{S}^2.}
In other words, the curve $T\mapsto (T,x-(T-t)\theta)$ converges to the straight line $T\mapsto(T,(x+t\cdot\theta-T)\theta)$ as $T\to\infty$. If the radiation field of $\phi_t-\phi_r$
\fm{\mcl{A}_0(q,\omega)=\lim_{r\to\infty}r(\phi_t-\phi_r)(r-q,r\omega)} exists for all $(q,\omega)\in\R\times\mathbb{S}^2$, then we expect
\fm{\lim_{T\to\infty}\int_{\mathbb{S}^2}\Phi(T,x-(T-t)\theta)\ dS_\theta=-\frac{1}{4\pi}\int_{\mathbb{S}^2}\mcl{A}_0(-t-x\cdot\theta,-\theta)\ dS_\theta.}
If one can also control the remainder term \eqref{sec:intro:prop4.1} when $T\to\infty$, then we have
\eq{\label{intro:phifor}\phi(t,x)=-\frac{1}{4\pi}\int_{\mathbb{S}^2}\mcl{A}_0(x\cdot\theta-t,\theta)\ dS_\theta+\text{(a remainder term)},\qquad \forall (t,x)\in\mcl{O}.}
Unfortunately, if $u$ solves \eqref{qwe}, $u_t-u_r$ does not necessarily admit a radiation field in the usual sense. In other words, we cannot apply \eqref{intro:phifor} directly to finish the proof of Theorem~\ref{mthm}. We thus require a more delicate analysis of the asymptotic behaviors of $u$ near the light cone.

In Section \ref{sec:pfmthm}, we prove Theorem \ref{mthm} by combining Proposition \ref{prop:inhomwe} with the results in Section~\ref{secreview}. For each multiindex $I$, we first prove a pointwise estimate for $\Box Z^Iu$. Then, by Proposition \ref{prop:inhomwe}, we notice that it remains to estimate the integral
\fm{\int_{\mathbb{S}^2}[(-\partial_t+\frac{y}{|y|}\cdot \partial_y)(|y|Z^Iu)](T,y(t,x,T,\theta))\ dS_\theta.}
Here $y(t,x,T,\theta)=x-(T-t)\theta$. In our proof, we set $T=t^3$ (or a sufficiently large power of $t$) and take $t\gg_{\delta,\eps}1$ and $|x|<t$. Under this setting, we have $||y|-T|\lesssim T^{1/3}$, so one  can apply Proposition \ref{prop:trzwtu} to approximate the integrand.

In Section \ref{sec:corofmthm}, we prove Theorems \ref{sdltthm} and \ref{ltthm_weak}. We first notice that Theorem \ref{ltthm_weak} follows from part c) of Theorem \ref{sdltthm}. That is, we need to show $\wh{A}=O(\lra{q}^{-1-B_0\eps})$ under one of the assumptions i)--iii). To achieve this goal, we apply either \eqref{prevac:uapp} or \eqref{mthmcl} to transfer the pointwise bounds for $u$ or $\partial u$ to the bounds for $\wh{A}$. For example, in the case when $u$ is spherically symmetric in $x$, the scattering data $\wh{A}$ is independent of $\omega$ (i.e.\ $\wh{A}=\wh{A}(q)$). By \eqref{mthmcl} and the bound $|u(t,0)|\lesssim\eps t^{-1-B_0\eps}$, we obtain immediately that $|\wh{A}(-t)|\lesssim t^{-1-B_0\eps}$. This is exactly what we want, so we conclude the proof of part iii). We also refer our readers to Proposition \ref{ltthm} which is slightly stronger than Theorem \ref{ltthm_weak}.

We now explain how to prove Theorem \ref{sdltthm}. Let us focus on parts b) and c), since the proof of part a) is similar and simpler. The proof relies on the explicit formula \eqref{asyeqn:sol:2A} of solutions to the geometric reduced system \eqref{asyeqn}. Recall that in Section \ref{sec1.2.3ac}, we defined $\wh{U}$ by 
\fm{\wh{U}_q(s,q,\omega)&=\wh{A}(q,\omega)\exp(\frac{1}{2}G(\omega)\wh{A}(q,\omega)s),\qquad \lim_{q\to\infty}\wh{U}(s,q,\omega)=0.}
We observe that  $|\wh{A}|e^{-Cs}\lesssim |\wh{U}_q|\lesssim |\wh{A}|e^{Cs}$ and that $\lim_{s\to\infty}\wh{U}_q\cdot 1_{G(\omega)\wh{A}<0}=0$.
Besides, we can show that a certain integral of $\wh{U}_q$ decays to $0$ as $t\to\infty$. For example, we have a decaying property in Remark \ref{sdltthm_rmk2} or Proposition \ref{prop:decaysphmeans}. In part c), we also have a pointwise decay for $\wh{U}$ in Lemma \ref{lem:est:wtu3}. Now we rewrite the integral of $\wh{U}_q$ as the sum of $\mcl{I}_++\mcl{I}_-+\mcl{I}_0$, where $\mcl{I}_\pm$ denotes the integral of $\wh{U}_q\cdot 1_{\pm G(\omega)\wh{A}>0}$ and $\mcl{I}_0$ denotes the integral of $\wh{U}_q\cdot 1_{G(\omega)\wh{A}=0}$. We send $t\to \infty$ and study the asymptotics of these three integrals. Under the assumptions of Theorem \ref{sdltthm}, we can show that $\mcl{I}_0$ can be neglected (because the null condition is violated), that $\mcl{I}_+\to0$ (because $\mcl{I}_-\to 0$ by the Lebesgue dominated convergence theorem), and that $e^{Cs}\mcl{I}_-\to0$ (because of the good decay rate of the integral of $\wh{U}_q$). In summary, we can show that $\wh{A}\equiv0$ and thus $u\equiv0$. We refer to Section \ref{sec:prop:vanishingA} for more details.

\subsection{Acknowledgement}
The author would like to thank Sung-Jin Oh and Daniel Tataru for many helpful discussions. The author would also like to thank the anonymous reviewer for valuable comments and suggestions. The author has been funded by a Simons Investigator grant of Daniel Tataru from the Simons Foundation and a VandyGRAF Fellowship from Vanderbilt University.

\section{Preliminaries}\label{sec:preliminaries}

\subsection{Notation}

We use $C$ to denote universal positive constants. We write $A\lesssim B$, $B\gtrsim A$, or $A=O(B)$ if $|A|\leq CB$ for some  $C>1$. The values of all constants in this paper may vary from line to line. We write $A\sim B$ if $A\lesssim B$ and $B\lesssim A$. We use $C_{v}$, $\gtrsim_v$, or $\lesssim_v$ if we want to emphasize that the constant depends on a parameter $v$. Moreover, if we make a statement such as ``$|A|\leq CB$ for all $v$ in a certain set $V$'', then it implies that the constant $C$ can be chosen to be uniform for all $v\in V$.

Unless specified otherwise, we always assume that the Latin indices $i,j,k$ take values in $\{1,2,3\}$ and the Greek indices $\alpha,\beta$ take values in $\{0,1,2,3\}$.   We use subscript to denote partial derivatives unless specified otherwise. For example, $u_{\alpha\beta}=\partial_\alpha\partial_\beta u$, $q_r=\partial_rq$, $A_q=\partial_qA$, etc. For a fixed integer $k\geq 0$, we use $\partial^k$ to denote either a specific partial derivative of order $k$, or the collection of partial derivatives of order $k$.

Given a function $f=f(\omega)$ defined on $\mathbb{S}^2$, in order to define angular derivatives $\partial_\omega$, we first extend $f$ to $\R^3\setminus 0$ by setting $f(x):=f(x/|x|)$ and then set $\partial_{\omega_i}f:=\partial_{x_i}f|_{\mathbb{S}^2}$. To prevent confusion, we will only use $\partial_\omega$ to denote angular derivatives under the coordinate system $(s,q,\omega)$, and will never use it under the coordinate system $(t,r,\omega)$. For a fixed integer $k\geq 0$, we will use $\partial_\omega^k$ to denote either a specific angular derivative of order $k$, or the collection of all angular derivatives of order $k$.

\subsection{Dependence on $\eps$}\label{prel:depeps}
In this paper, the parameter $\eps$ always denotes the size of initial data in \eqref{init}. We always assume that $\eps\ll 1$ which means $0<\eps<\eps_0$ for some sufficiently small constant $0<\eps_0<1$. We say that $\eps \ll1$ depends on $v$, or $\eps\ll_v1$, if we want to emphasize that this constant $\eps_0$ depends on a parameter $v$. 

The convention in the previous subsection allows us to write a constant depending on $v$ as $C$ without the subscript $v$. However, we make a special exception for $\eps$. In this paper, if a constant depends on  $\eps$, then we must write it as $C_{\eps}$. In other words, any constant without the subscript $\eps$ must be independent of the choice of $\eps$. Similarly for $\lesssim_\eps$, $\gtrsim_\eps$, or $\ll_\eps$. 

Note that most statements in this paper follow this format: For each integer $N\geq 1$, we have a certain property $(P_{N,\eps})$ as long as $\eps\ll_N1$. For simplicity, later we will only say that we have $(P_{N,\eps})$ for $\eps\ll1$ without explicitly stating the dependence of $\eps$ on $N$. For example, we can say
\fm{(P_{N,\eps}):\ t^{-1+C_N\eps}\leq 1,\qquad \forall t\geq 1}
as long as $\eps\ll 1$.
Note that $(P_{N,\eps})$ holds for all $\eps\in(0,\frac{1}{1+C_N})$.

\subsection{Commuting vector fields}

Let $Z$ be any of the following vector fields:
\eq{\label{c1vf} \partial_\alpha,\ \alpha=0,1,2,3;\ S=t\partial t+r\partial_r=t\partial_t+\sum_{i=1}^3x_i\partial_i;\\
\Omega_{ij}=x_i\partial_j-x_j\partial_i,\ 1\leq i<j\leq 3;\ \Omega_{0i}=x_i\partial_t+t\partial_i,\ i=1,2,3.}
We write these vector fields as $Z_1,Z_2,\dots,Z_{11}$, respectively. For any multiindex $I=(i_1,\dots,i_m)$ with length $m=|I|$ such that $1\leq i_*\leq 11$, we set $Z^I=Z_{i_1}Z_{i_2}\cdots Z_{i_m}$. Then we have Leibniz's rule
\begin{equation}Z^I(fg)=\sum_{|J|+|K|=|I|}C^I_{JK}Z^Jf\cdot Z^Kg,\qquad\text{where $C_{JK}^I$ are constants.}\end{equation}

We have the following commutation properties.
\begin{equation}[S,\square]=-2\square,\qquad[Z,\square]=0\text{ for other $Z$};\end{equation}
\begin{equation}[Z,\wt{Z}]=\sum_{|I|=1} C_{Z,\wt{Z},I}Z^I,\qquad\text{where $C_{Z,\wt{Z},I}$ are constants};\end{equation}
\begin{equation}[Z,\partial_\alpha]=\sum_\beta C_{Z,\alpha\beta}\partial_\beta,\hspace{2em} \text{where $C_{Z,\alpha\beta}$ are constants}.\end{equation}

We now recall several useful estimates involving the commuting vector fields. Their proofs are standard and can be found in, e.g., \cite{MR2455195,MR1466700,MR2382144}. We start with the following pointwise estimates.
\begin{lem}\label{c1l2.1}
For any function $\phi=\phi(t,x)$, we have 
\begin{equation}\label{c1l2.1f1}|\partial^k\phi|\leq C\lra{t-r}^{-k}\sum_{|I|\leq k}|Z^I\phi|,\qquad\forall k\geq 0;\end{equation}
\begin{equation}\sum_{i=1}^3|(\partial_i+\omega_i\partial_t)\phi|\leq C\lra{t+r}^{-1}|Z \phi|.\end{equation}
For each $x\in\R^3$, we set $r=|x|$ and $\omega=(\omega_i)_{i=1,2,3}=x/|x|$. For each $s\in\R$, we define the Japanese bracket $\lra{s}:=\sqrt{1+|s|^2}$. We also define  $|Z\phi|:=\sum_{|I|=1}|Z^I\phi|$.
\end{lem}\rm

In addition, we have the Klainerman-Sobolev inequality.  
\begin{prop} For $\phi\in C^\infty(\R^{1+3})$ which vanishes for large $|x|$, we have
\begin{equation}\lra{t+|x|}\lra{t-|x|}^{1/2}|\phi(t,x)|\leq C\sum_{|I|\leq 2}\norm{Z^I\phi(t,\cdot)}_{L^2(\R^3)}.\end{equation}
\end{prop}\rm

Moreover, we state a corollary of \cite[Lemma 2.6]{yu2022nontrivial}. It can be viewed as the estimates for Taylor's series
adapted to the commuting vector fields. Note that the original lemma is stronger but unnecessary for this paper.
\lem{\label{lem:prel:taylor}Fix two positive integer $N,l_0$. Assume that $\phi=\phi(t,x)$ is an $\R^k$-valued function satisfying $\sum_{|J|\leq N}|Z^J\phi|\leq 1$ at a certain point $(t_0,x_0)$. Also assume that $f\in C^\infty(\R^k)$ and that $f(s)=O(|s|^{l_0})$ whenever $s$ is near the origin. Then, at $(t_0,x_0)$ we have
\fm{\sum_{|J|\leq N}|Z^J(f(\phi))|\lesssim_{f,N} \sum_{|J|\leq N}|Z^J\phi|^{l_0}.}
The implicit constant is independent of the choice of $(t_0,x_0)$.}\rm

\bigskip

We also recall \cite[Lemma 2.1]{yu2022nontrivial} on the commutator $[\partial_t-\partial_r,Z^I]$.
\lem{\label{sec:prel:commutetrz}For each multiindex $I$, we have
\fm{\ [\partial_t-\partial_r,Z^I]F&=\sum_{|J|<|I|}\kh{f_0\cdot Z^J(F_t-F_r)+\sum_{i=1,2,3}f_0\cdot (\partial_i+\omega_i\partial_t)Z^JF}.}
Here $f_0$ denotes an arbitrary finite sum of terms of the form $C\cdot \prod_{j=1}^p Z^{I_j}\omega_{i_j}$ with $i_*\in\{1,2,3\}$, and we allow $f_0$ to vary from line to line.}
\rm

\subsection{A function space}
To simplify our computations, we make the following definition. Fix a small parameter $\delta\in(0,1)$. For each $\eps\in(0,1)$, we suppose that there is a nonempty open set $\D_\eps\subset\{(t,x):\ t\geq e^{\delta/\eps},\ |x|/t\in[1/2,2]\}$.
\defn{\label{prel:defn:SD}\rm Fix $n,m,p\in\R$. Suppose that we have a function $F=F(t,x;\eps)$ that is defined for all $(t,x)\in\D_\eps$ as long as $\eps\ll 1$. We say that $F\in \eps^nS^{m,p}_\D$ if for each integer $N\geq1$, there exists an $\eps_N\in(0,1)$ depending on $\{\D_*\},\delta,N$, such that for all $\eps\in(0,\eps_N)$, we have $F(\cdot;\eps)\in C^N(\D_\eps)$  and 
\eq{\label{prel:fs:est}\sum_{|I|\leq N}|Z^IF(t,x;\eps)|\lesssim_N \eps^nt^{m+C_N\eps}\lra{r-t}^{p},\qquad \forall (t,x)\in\D_\eps.}
The constants in \eqref{prel:fs:est} depend on $\{\D_*\}$, $N$, and $\delta$, but they are independent of $\eps\in(0,\eps_N)$.

Moreover, we say that $F\in \eps^nS^{m,p}_{\ln,\D}$ if we add a factor $(1+\ln \lra{r-t})$ on the right hand side of \eqref{prel:fs:est}. Such a definition is useful in the proof of Proposition \ref{prop:trzwtu}.

Finally, we will omit the subscript $\D$ for simplicity if it is clear what $\D_\eps$ is according to the context. That is, we can write $\eps^n S^{m,p}$ and $\eps^n S^{m,p}_{\ln}$. If we have $n=0$, then we will also omit $\eps^0$ and write $S^{m,p}$ and $S^{m,p}_{\ln}$.
}
\rm

\bigskip

We summarize some useful properties in the next lemma. This lemma is essentially the same as \cite[Lemma 2.8]{yu2022nontrivial}, so we omit the proof here.
\lem{\label{c1lem6}
\begin{enumerate}[\rm (a)]
\item We have $r-t\in S^{0,1}$, $r^{m},t^{m},(t+r)^{m}\in S^{m,0}$, and $\omega\in S^{0,0}$.

\item For any $F_1\in\eps^{n_1} S^{m_1,p_1}$ and $F_2\in \eps^{n_2}S^{m_2,p_2}$, we have \fm{F_1+ F_2\in \eps^{\min\{n_1,n_2\}}S^{\max\{m_1,m_2\},\max\{p_1,p_2\}},\hspace{2em}F_1F_2\in\eps^{n_1+n_2} S^{m_1+m_2,p_1+p_2}.}

\item For any $F_1\in\eps^{n_1} S^{m_1,p_1}_{\ln}$ and $F_2\in \eps^{n_2}S^{m_2,p_2}$, we have \fm{F_1+ F_2\in \eps^{\min\{n_1,n_2\}}S_{\ln}^{\max\{m_1,m_2\},\max\{p_1,p_2\}},\hspace{2em}F_1F_2\in\eps^{n_1+n_2} S_{\ln}^{m_1+m_2,p_1+p_2}.}

\item For any $F\in\eps^n S^{m,p}$,  we have $ZF\in \eps^nS^{m,p}$, $\partial F\in \eps^nS^{m,p-1}$ and $(\partial_i+\omega_i\partial_t)F\in \eps^nS^{m-1,p}$.
\end{enumerate}}\rm

\bigskip

We also present a key lemma that will be applied repeatedly in this paper.

\lem{\label{c1lem6.1}Let $\rho=\rho(t,x;\eps)\in S^{0,1}$ such that $\lra{r-t}\lesssim \lra{\rho(t,x;\eps)}t^{C\eps}$ in $\D_\eps$ for $\eps\ll1$. Let $\wt{\D}_\eps\subset [0,\infty)\times\R\times\mathbb{S}^2$ be an open set such that \fm{\wt{\D}_\eps\supset\{(\eps\ln t-\delta,\rho(t,x;\eps),x/|x|):\ (t,x)\in\D_\eps\}.} Suppose that $G=G(s,q,\omega;\eps)$ is a function such that for each integer $N\geq 0$, as long as $\eps\ll_N1$, we have $G(\cdot;\eps)\in C^N(\wt{D}_\eps)$ and  
\eq{\label{c1lem6.1:asu}\sum_{a+b+c\leq N}\lra{q}^{a}\eps^{b}|\partial_q^a\partial_s^b\partial_\omega^cG(s,q,\omega;\eps)|\lesssim_N \eps^n\lra{q}^{p}e^{C_Ns+m\cdot \frac{s+\delta}{\eps}},\qquad (s,q,\omega)\in\wt{\D}_\eps.}
Then we have $G(\eps\ln t-\delta,\rho(t,x;\eps),x/|x|;\eps)\in \eps^nS^{m,p}$.}
\begin{proof}
We first notice that \fm{\lra{r-t}t^{-C\eps}\lesssim \lra{\rho}\lesssim t^{C\eps}\lra{r-t},\qquad \forall (t,x)\in\D\cap\{t\geq e^{\delta/\eps}\}.}
The first half is in the assumption of this lemma, and the second half follows from $\rho\in S^{0,1}$.

Fix an integer $N\geq 0$. Choose $\eps\ll_N1$ so that $G(s,q,\omega;\eps)\in C^{N}(\wt{\D}_\eps)$, $\rho(t,x;\eps)\in C^N(\D_\eps)$, and that both \eqref{prel:fs:est} (with $F$ replaced by $\rho$ and $(n,m,p)=(0,0,1)$) and \eqref{c1lem6.1:asu} hold for this $N$. We need to check \eqref{prel:fs:est} with $F$ replaced by $G(t,x;\eps):=G(\eps\ln t-\delta,\rho(t,x;\eps),x/|x|;\eps)$ and $m=0$.
Fix a multiindex $I$ with $|I|\leq N$. If $|I|=0$, then \eqref{c1lem6.1:asu} implies that \fm{|G|\lesssim \eps^n\lra{\rho}^p t^{m+C\eps}\lesssim \eps^n\lra{r-t}^p t^{m+C\eps}\qquad \forall (t,x)\in \D_\eps,} so we are done. If $|I|>0$, we write $Z^I(G(t,x;\eps))$ as a linear combination (with real constant coefficients independent of $\eps$) of terms of the form
\fm{\partial_q^a\partial_s^b\partial_\omega^cG(\eps\ln t-\delta,\rho(t,x;\eps),x/|x|;\eps)\cdot\prod_{j=1}^bZ^{I_j}(\eps\ln t-\delta)\cdot\prod_{j=1}^aZ^{J_j}\rho\cdot\prod_{j=1}^cZ^{K_j}\omega.}
Here $a+b+c\leq N$, each of $I_*,J_*,K_*$ is nonzero, and the sum of these multiindices is $|I|$. Since $|I_*|>0$, we have $Z^{I_*}(\eps\ln t-\delta)=O(\eps)$. Since $\rho\in S^{0,1}$ and $\omega\in S^{0,0}$, we have $Z^{J_*}\rho=O(t^{C\eps}\lra{r-t}) $ and $Z^{K_*}\omega=O(1)$ (we can show that there is no $t^{C\eps}$ factor here). We also apply \eqref{c1lem6.1:asu} to control the derivatives of $G$. In summary, each term of the form above is controlled by 
\fm{\eps^{n-b}\lra{\rho}^{p-a}t^{m+C\eps}\cdot \eps^b\cdot \lra{r-t}^{a}t^{C\eps}\cdot 1\lesssim\eps^n\lra{r-t}^{p}t^{m+C\eps}. }
Thus, we have $Z^IG=O(\eps^n\lra{r-t}^{p}t^{m+C\eps})$ in $\D_\eps$.
\end{proof}
\rm

\section{A review on the asymptotic completeness result}\label{secreview}
The proofs of the main theorems rely heavily on the author's previous paper \cite{MR4772266} on the asymptotic completeness for \eqref{qwe}.
In this section, we summarize the results in \cite{MR4772266} that are necessary for the proofs in this paper. In Section \ref{secreview:globalsol}, we quickly review the global existence result in \cite{MR2382144}. In Section \ref{sec:prevac:opticalfunction}, we discuss the construction of the optical function $q$. In Sections \ref{prevac:grs} and \ref{secreview:secapp}, we define the scattering data $\wh{A}$ for the asymptotic completeness problem and show that one can construct a matching approximate solution $\wt{u}$ to \eqref{qwe}. Finally, in Section \ref{secreview:gauage}, we discuss a gauge independence result.

We also need some new results whose proofs rely on the computations in \cite{MR4772266}. In Section \ref{secreview:furthercomputations}, we prove Proposition \ref{prop:trzwtu} on the higher derivatives of the approximate solution $\wt{u}$ constructed in Section \ref{sec1.2.3ac}. The proof of Proposition \ref{prop:trzwtu} is lengthy, so we add a brief overview of its proof after its statement. In Section \ref{secreview:gauage}, we present a new version of the gauge independence result. By changing the definition of the scattering data, we show Corollary \ref{secreview:gauge:prop:cor} which is of a much cleaner form than the one in \cite[Section 6]{MR4772266}.

\subsection{Global existence}\label{secreview:globalsol}
Fix $u_0,u_1\in C_c^\infty(\R^3)$. Choose $R>0$ such that $\supp(u_0,u_1)\subset\{|x|\leq R\}$, and we remark that the meaning of $R$ will remain unchanged in the rest of this paper. For $\eps\ll1$,  the Cauchy problem \eqref{qwe} and \eqref{init} admits a unique global $C^\infty$ solution $u=u(t,x;\eps)$ for all $t\geq 0$. Sometimes we omit $\eps$ and write $u(t,x)$ for simplicity. Here we first apply the global existence result for \eqref{qwe} in Lindblad \cite{MR2382144} to get a global $C^{14}$ solution, and then apply \cite[Theorems I.4.1 and I.4.2]{MR2455195} to show that this $C^{14}$ solution is indeed the unique $C^\infty$ solution.

Moreover, for each nonnegative integer $N$, we have the following pointwise bounds as long as $\eps\ll_N1$:
\eq{\label{secreview:uptb}\sum_{|I|\leq N}|Z^Iu(t,x;\eps)|\lesssim_N \eps \lra{t}^{-1+C_N\eps},\qquad \forall t\geq 0,\ x\in\R^3.}
There is a better bound for $\partial u$: $|\partial u|\lesssim \eps \lra{t}^{-1}$.
By the finite speed of propagation, we also have $u\equiv 0$ for $r-t\geq R$.

\subsection{The optical function}\label{sec:prevac:opticalfunction}
Fix a small parameter $\delta\in(0,1)$. Define 
\eq{\label{secreview:defnOmega}\Omega=\Omega_{\delta,R,\eps}:=\{(t,x)\in(e^{\delta/\eps},\infty)\times\R^3:\ |x|>\frac{1}{2}(t+e^{\delta/\eps})+2R\}.}
Heuristically, one can view $\Omega$ as $\{t>e^{\delta/\eps},\ |x|>t/2\}$. Here we use \eqref{secreview:defnOmega} for technical reasons; see the discussion below \eqref{eik}. From now on, if we use the notation $\eps^nS^{m,p}$ or $\eps^nS^{m,p}_{\ln}$, we always take $\D_{\eps}=\Omega_{\delta,R,\eps}\cap \{r-t<2R\}$ in Definition \ref{prel:defn:SD}. Thus, by \eqref{secreview:uptb}, we have $u\in \eps S^{-1,0}$.

In \cite[Propositions 3.1 and 4.1]{MR4772266}, we applied the method of characteristics and proved that the eikonal equation
\eq{\label{eik}g^{\alpha\beta}(u)\partial_\alpha q\partial_\beta q=0\quad\text{in }\Omega;\qquad q=|x|-t\quad \text{on }\partial\Omega}
admits a unique global $C^2$ solution $q=q(t,x;\eps)$ in $\Omega$. Because $\partial\Omega$ is not a smooth surface, it seems that we need to carefully choose the boundary condition near the ``corner'' of $\partial\Omega$
\fm{\{t=e^{\delta/\eps},|x|=e^{\delta/\eps}+2R\}\subset\partial\Omega}
to make $q$ a $C^2$ solution. However, it is not necessary here because of our choice \eqref{secreview:defnOmega} of $\Omega$. We have $u|_{r-t\geq R}\equiv 0$ and $q|_{r-t\geq R}=r-t$, so we only need to consider those characteristics emanating from a smooth subset of $\partial\Omega$ (i.e.\ $\partial\Omega\cap\{t>e^{\delta/\eps}\}$).  Besides, we have $q\in S^{0,1}$ and $\lra{r-t}\lesssim \lra{q}t^{C\eps}$ in $\Omega$, so we can take $\rho=q$ in Lemma \ref{c1lem6.1}. To see this, we apply \cite[Lemma 3.8]{MR4772266} to obtain $|q-(r-t)|\lesssim t^{C\eps}$ and $t^{-C\eps}\lesssim \lra{q}/\lra{r-t}\lesssim t^{C\eps}$ in $\Omega$. Next, we apply \cite[Proposition 4.1]{MR4772266} to conclude that for each fixed integer $N\geq 2$, the $C^2$ solution $q$ is $C^N$ as long as $\eps\ll_N1$. We also have $\sum_{|I|\leq N}|Z^Iq|\lesssim \lra{q}t^{C\eps}\lesssim \lra{r-t}t^{C\eps}$ in $\Omega$.

Next, if we set 
\eq{\label{secreview:defnOmega'}\Omega'&:=\{(s,q,\omega)\in[0,\infty)\times\R\times\mathbb{S}^2:\ q>\frac{1}{2}(1-e^{s/\eps})e^{\delta/\eps}+2R\},}
then  the map
\fm{\Omega\to\Omega':\ (t,x)\mapsto (\eps\ln t-\delta,q(t,x;\eps),x/|x|)}
is bijective and has an inverse map. As a result, any function defined for all $(t,x)\in\Omega$ induces a function of $(s,q,\omega)$. Again, for each fixed integer $N\geq 2$, both the map $\Omega\to\Omega'$ and its inverse are $C^N$ as long as $\eps\ll_N1$. 

\subsection{The geometric reduced system}\label{prevac:grs}
In $\Omega$, we set $(\mu,U)(t,x;\eps):=(q_t-q_r,\eps^{-1}ru)(t,x;\eps)$ and obtain an induced function $(\mu,U)(s,q,\omega;\eps)$ in $\Omega'$. In \cite[Section 5]{MR4772266}, we proved that $(\mu,U)$ is an approximate solution to the geometric reduced system \eqref{asyeqn}. In addition, the following three limits exist for all $(q,\omega)\in\R\times\mathbb{S}^2$:
\eq{\label{prevac:limita}\left\{\begin{array}{l}
\displaystyle A(q,\omega;\eps):=-\frac{1}{2}\lim_{s\to\infty}(\mu U_q)(s,q,\omega;\eps),\\
\displaystyle A_1(q,\omega;\eps):=\lim_{s\to\infty}\exp(\frac{1}{2}G(\omega)A(q,\omega;\eps)s)\mu(s,q,\omega;\eps),\\
\displaystyle A_2(q,\omega;\eps):=\lim_{s\to\infty}\exp(-\frac{1}{2}G(\omega)A(q,\omega;\eps)s)U_q(s,q,\omega;\eps).
\end{array}\right.}
Note that $A_1A_2\equiv -2A$. Since $u\equiv 0$ whenever $r\geq t+R$, we have $(A,A_2)\equiv 0$ and $A_1\equiv -2$ for $q\geq R$. By setting
\eq{\label{prevac:defnwtmuU}\left\{\begin{array}{l}
\displaystyle\wt{\mu}(s,q,\omega;\eps):=A_1\exp(-\frac{1}{2}G(\omega)As),\\[1em]
\displaystyle \wt{U}_q(s,q,\omega;\eps):=A_2\exp(\frac{1}{2}G(\omega)As),\\[1em]\displaystyle 
\lim_{q\to\infty}\wt{U}(s,q,\omega;\eps)=0,
\end{array}\right.}
we obtain an exact solution $(\wt{\mu},\wt{U})$ to the reduced system \eqref{asyeqn}. We refer our readers to \cite[Proposition 5.1]{MR4772266} for more details. Since $A,A_1,A_2$ are defined for all $(q,\omega)\in\R\times\mathbb{S}^2$, we emphasize that $(\wt{\mu},\wt{U})$ is defined for all $(s,q,\omega)\in\R\times\R\times\mathbb{S}^2$. On the other hand, $(\mu,U)$ is only defined in $\Omega'$.

In \cite[Proposition 5.1]{MR4772266}, we also proved several bounds for $A,A_1,A_2$ in the coordinate set $(q,\omega)$: 
\eq{\label{prevac:esta}\partial_q^a\partial_\omega^c(A,A_2)=O_{a,c}(\lra{q}^{-1-a+C_{a,c}\eps}),\quad \partial_q^a\partial_\omega^cA_1=O_{a,c}(\lra{q}^{-a+C_{a,c}\eps}),\qquad \forall a,c\geq 0.}
Here we remind our readers of the convention introduced in Section \ref{prel:depeps}. First, the estimates \eqref{prevac:esta} should be understood in the following way: for each integer $N\geq 0$, as long as $\eps\ll_N1$, the functions $A,A_1,A_2$ are $C^N$ and satisfy \eqref{prevac:esta} for all $a+c\leq N$. Moreover, we emphasize that the constants in \eqref{prevac:esta} can be chosen to be uniform for all $\eps\ll_N1$.

In \cite[Lemma 5.8]{MR4772266}, we proved that
\eq{|A_1+2|\leq \lra{q}^{-1+C\eps},\qquad -3\leq A_1\leq -1<0.}

Finally, we recall from \cite[Proposition 5.1, (5.11), and (5.12)]{MR4772266} that 
\eq{\label{prevac:estmore}\partial_q^a\partial_s^b\partial_\omega^c(\mu-\wt{\mu})&=O_{a,b,c}(\eps^{-b}\lra{q}^{-a}t^{-1+C_{a,b,c}\eps}),\qquad &\forall a,b,c\geq 0,\ (s,q,\omega)\in\Omega';\\
\partial_q^a\partial_s^b\partial_\omega^c(U-\wt{U})&=O_{a,b,c}(\eps^{-b}\lra{q}^{1-a}t^{-1+C_{a,b,c}\eps}),\qquad& \forall a,b,c\geq 0,\ (s,q,\omega)\in\Omega';}
\eq{\label{prevac:estmore2}
\partial_q^a\partial_s^b\partial_\omega^c(\wt{\mu},\wt{U})&=O_{a,b,c}(\eps^{-b}\lra{q}^{-a}t^{C_{a,b,c}\eps}),\qquad &\forall a,b,c\geq 0,\ (s,q,\omega)\in\Omega'.}
All the functions above are defined for $(s,q,\omega)\in \Omega'\subset[0,\infty)\times\R\times\mathbb{S}^2$ where $\Omega'$ is defined by \eqref{secreview:defnOmega'}. Besides, we set $t=e^{\frac{s+\delta}{\eps}}$. Again, all the estimates above should be understood in the following way: for each integer $N\geq0$, we can choose $\eps\ll_N1$ so that these estimates hold for all $a+b+c\leq N$.

As discussed in Section \ref{sec:prevac:opticalfunction}, every function of $(s,q,\omega)$ induces a function of $(t,x)$ via the map $(t,x)\mapsto (\eps\ln t-\delta,q(t,x;\eps),x/|x|)$.  Below \eqref{eik}, we mentioned that $q\in S^{0,1}$ and $\lra{r-t}\lesssim \lra{q}t^{C\eps}$. Thus, we can apply Lemma \ref{c1lem6.1} with $\rho=q$ to obtain the following lemma.

\lem{\label{prevac:estmore:lem} We have  $u-\eps r^{-1}\wt{U}\in\eps S^{-2,1}$, and $\partial_\alpha u-\eps\wh{\omega}_\alpha r^{-1}A\in \eps S^{-2,0}$ for $\alpha=0,1,2,3$. Here $\wh{\omega}_0=-1$ and $\wh{\omega}_j=\omega_j$ for $j=1,2,3$. That is, for each fixed integer $N\geq 0$, as long as $\eps\ll_N1$,   we have
\fm{\sum_{|I|\leq N}|Z^I(u(t,x)-\eps r^{-1}\wt{U}(\eps\ln t-\delta,q(t,x),x/|x|))|&\lesssim \eps t^{-2+C\eps}\lra{r-t},\\
\sum_{|I|\leq N}|Z^I(\partial_\alpha u(t,x)-\eps \wh{\omega}_\alpha r^{-1}A(q(t,x),x/|x|))|&\lesssim \eps t^{-2+C\eps}}
for all $(t,x)\in\Omega$. Note that the estimates hold trivially when $r-t\geq R$ since $u,\wt{U},A$ vanish there.

}
\begin{proof}
We apply Lemma \ref{c1lem6.1} with $\rho=q$, $G=U-\wt{U}$, and $(n,m,p)=(0,-1,1)$. Thus, we have $U-\wt{U}\in \ S^{-1,1}$. By Lemma \ref{c1lem6}, we have $\eps r^{-1}\in \eps S^{-1,0}$ and thus $u-\eps r^{-1}\wt{U}=\eps r^{-1}(U-\wt{U})\in \eps S^{-1,0}\cdot S^{-1,1}\subset\eps S^{-2,1}$.

Next, by part (d) of Lemma \ref{c1lem6}, we have $\partial_{\alpha}(u-\eps r^{-1}\wt{U})\in \eps S^{-2,0}$. It  suffices to prove that $\partial_{\alpha}(\eps r^{-1}\wt{U})-\eps \wh{\omega}_\alpha r^{-1}A\in \eps S^{-2,0}$. By the chain rule and since $A=-\frac{1}{2}\wt{\mu}\wt{U}_q$, we have
\fm{&\partial_\alpha(\eps r^{-1}\wt{U})-\eps r^{-1}\wh{\omega}_\alpha A\\
&=\eps r^{-1}(\wt{U}_qq_\alpha-\wh{\omega}_\alpha A)-\eps r^{-2}(\partial r)\cdot \wt{U}+\eps r^{-1}\wt{U}_s\cdot\partial(\eps\ln t)+\eps r^{-1}\wt{U}_{\omega}\cdot\partial\omega\\
&=\eps r^{-1}(q_\alpha+\frac{1}{2}\wh{\omega}_\alpha \wt{\mu})\wt{U}_q+\eps S^{-2,0}.}
To obtain the last estimate, we first apply \eqref{prevac:estmore2} and Lemma \ref{c1lem6.1} to obtain $\wt{U},\wt{U}_\omega\in S^{0,0}$ and $\wt{U}_s\in \eps^{-1}S^{0,0}$. By Lemma \ref{c1lem6}, we have $\eps r^{-2}\partial r,\eps r^{-1}\partial\omega\in \eps S^{-2,0}$ and $\eps r^{-1}\partial(\eps\ln t)\in \eps^2S^{-2,0}$.

Finally, we need to show $\eps r^{-1}(q_\alpha+\frac{1}{2}\wh{\omega}_\alpha \wt{\mu})\wt{U}_q\in \eps S^{2,0}$. We have $\eps r^{-1}\in \eps S^{-1,0}$ and $\wt{U}_q\in S^{0,-1}$ by \eqref{prevac:estmore2} and Lemma \ref{c1lem6.1}. Moreover,  we have
\fm{q_t&=\frac{1}{2}\wt{\mu}+\frac{1}{2}(\mu-\wt{\mu})+\frac{1}{2}(q_t+q_r),\\
q_i&=-\frac{1}{2}\wt{\mu}\omega_i-\frac{1}{2}(\mu-\wt{\mu})\omega_i+(q_i-\omega_iq_r).}
By \eqref{prevac:estmore} and Lemma \ref{c1lem6.1}, we have $\mu-\wt{\mu}\in S^{-1,0}$. Since $q\in S^{0,1}$, we apply part (d) of Lemma \ref{c1lem6} to obtain $q_t+q_r,q_i-\omega_iq_r\in S^{-1,1}$. As a result, we have $q_\alpha+\frac{1}{2}\wt{\mu}\wh{\omega}_\alpha\in S^{-1,1}$. In summary, we have
\fm{\eps r^{-1}(q_\alpha+\frac{1}{2}\wh{\omega}_\alpha \wt{\mu})\wt{U}_q\in \eps S^{-1,0}\cdot S^{-1,1}\cdot S^{0,-1}\subset \eps S^{-2,0}.}
\end{proof}\rm

\subsection{Approximation}\label{secreview:secapp} In \cite[Section 7]{MR4772266}, we used the limits in \eqref{prevac:limita} to generate an approximate solution to \eqref{qwe} that is sufficiently close to the exact solution $u$ in some sense. The construction is as follows. Recall from \eqref{prevac:defnwtmuU}  that we have an exact solution $(\wt{\mu},\wt{U})(s,q,\omega;\eps)$ to the geometric reduced system \eqref{asyeqn}. Now, we define a new function $\wt{q}=\wt{q}(t,x;\eps)$ in $\Omega$ by solving
\eq{\wt{q}_t-\wt{q}_r=\wt{\mu}(\eps\ln t-\delta,\wt{q}(t,x;\eps),x/|x|;\eps)\quad\text{in }\Omega;\qquad \wt{q}=|x|-t\text{ for }|x|-t\geq 2R.}
We then define
\eq{\wt{u}(t,x;\eps)&=\eps |x|^{-1}\wt{U}(\eps\ln t-\delta,\wt{q}(t,x;\eps),x/|x|;\eps),\qquad \forall (t,x)\in\Omega.}
Once again, both $\wt{q}$ and $\wt{u}$ are $C^N$ as long as  $\eps\ll_N1$. In \cite[Proposition 7.1 and Lemma 7.7]{MR4772266}, we proved that this $\wt{u}$ is an approximate solution to \eqref{qwe}, that for each integer $N\geq 0$, as long as $\eps\ll_N1$,
\eq{\label{prevac:wtu}\sum_{|I|\leq N}|Z^I\wt{u}(t,x;\eps)|\lesssim \eps t^{-1+C\eps},\qquad \forall (t,x)\in\Omega,} and that for each $\gamma\in(0,1)$ and each integer $N\geq 0$, as long as $\eps\ll_{\gamma,N}1$,
\eq{\label{prevac:uapp}\sum_{|I|\leq N}|Z^I(u-\wt{u})(t,x;\eps)|\lesssim_{\gamma,N}\eps t^{-2+C_N\eps}\lra{r-t},\qquad \forall (t,x)\in\Omega\cap\{|r-t|\lesssim t^\gamma\}.}
We ask our readers to compare \eqref{prevac:uapp} with Lemma \ref{prevac:estmore:lem}. First, the estimates \eqref{prevac:uapp} hold in $\Omega\cap\{|r-t|\lesssim t^\gamma\}$ and the implicit constants depend on $\gamma$. From this point of view, the estimates in Lemma \ref{prevac:estmore:lem} are better since they hold in $\Omega$. However, we can construct $\wt{u}$ from the functions $A,A_1,A_2$ without knowing what the exact solution $u$ to \eqref{qwe} and the exact optical function $q$ are. Recall that we construct an approximate optical function $\wt{q}$ by solving a transport equation which does not depend on $u$ and $q$. However, to apply Lemma \ref{prevac:estmore:lem}, we need the exact optical function $q$ which is determined by the eikonal equation \eqref{eik} involving~$u$.

Here we recall an alternative way to define $\wt{u}$ from \cite[Section 7.2]{MR4772266} which will simplify our computations. Define\footnote{To motivate this definition, we  recall that $\wt{q}_t-\wt{q}_r=\wt{\mu}$. If we set $\wh{q}=F(\wt{q},\omega)$ and set $\wh{\mu}=\wh{q}_t-\wh{q}_r$, then $\wh{\mu}=F_q\cdot (\wt{q}_t-\wt{q}_r)=F_q\wt{\mu}$. Recall that $\wt{\mu}|_{s=0}=A_1$ and that our goal is to make $\wh{\mu}|_{s=0}=-2$. This suggests that we should choose $F$ so that $A_1F_q=-2$.}
\eq{\label{prevac:defF}F(q,\omega;\eps):=2R-\int_{2R}^q \frac{2}{A_1(p,\omega;\eps)}\ dp,\qquad\forall (q,\omega)\in\R\times\mathbb{S}^2.}
Recall that $A_1\in[-3,-1]$ everywhere, so $F$ is well defined. In \cite[Section 7.2]{MR4772266}, we proved that $\lra{F(q,\omega;\eps)}\sim \lra{q}$ for all $(q,\omega)\in\R\times\mathbb{S}^2$.
For each fixed $\omega\in\mathbb{S}^2$, the map $q\mapsto F(q,\omega;\eps)$ has an inverse $\wh{F}(q,\omega;\eps)$. Set
\eq{\label{prevac:defhata}
\left\{\begin{array}{l}
\displaystyle \wh{A}(q,\omega;\eps):=A(\wh{F}(q,\omega;\eps),\omega;\eps),\\
\displaystyle \wh{\mu}(s,q,\omega;\eps):=-2\exp(-\frac{1}{2}G(\omega)\wh{A}(q,\omega;\eps)s),\\
\displaystyle \wh{U}(s,q,\omega;\eps):=-\int_q^\infty \wh{A}(p,\omega;\eps)\exp(\frac{1}{2}G(\omega)\wh{A}(p,\omega;\eps)s) \ dp.\end{array}\right.}
This $\wh{A}$ is called the \emph{scattering data} for the asymptotic completeness problem.
Note that $(\wh{\mu},\wh{U})$ is again a solution to the geometric reduced system \eqref{asyeqn}.

We summarize several estimates for $\wh{A},\wh{\mu},\wh{U}$ in the next lemma. We again remind our readers of the convention in Section \ref{prel:depeps}.
`
\lem{\label{prevac:lem:esthat}For all $a,c\geq 0$, we have
\eq{\label{prevac:estahat} \partial_q^a\partial_\omega^c\wh{A}=O_{a,c}(\lra{q}^{-1-a+C_{a,c}\eps}),\qquad \forall (q,\omega)\in\R\times\mathbb{S}^2.}
For all $a,b,c\geq 0$, we have
\eq{\label{prevac:estmuhat} \partial_q^a\partial_s^b\partial_\omega^c(\wh{\mu}+2)&=O_{a,b,c}(\lra{q}^{-\max\{1,b\}-a+C_{a,b,c}\eps}e^{C_{a,b,c}s});\\
\partial_q^a\partial_s^b \partial_\omega^c\wh{U}_q&=O_{a,b,c}(\lra{q}^{-1-a-b+C_{a,b,c}\eps}e^{C_{a,b,c}s});\\
\partial_s^b \partial_\omega^c\wh{U}&=O_{b,c}((\eps^{-1}\lra{q}^{C_{b,c}\eps}1_{b=0}+1_{b>0})e^{C_{b,c}s})}
for all $(s,q,\omega)\in[0,\infty)\times\R\times\mathbb{S}^2$.
Finally, for all $c\geq 0$, $s\geq 0$, $\omega\in \mathbb{S}^2$, and $q>C_0e^{C_0s}(e^{\delta/\eps}-e^{(\delta+s)/\eps}-1)$ for some constant $C_0>1$, we have
\eq{\label{prevac:estUhat} \partial_\omega^c\wh{U}&=O_{c,C_0}(e^{C_{c,C_0}s}).}}
\begin{proof}
The estimates  \eqref{prevac:estahat} have been proved in \cite[Proposition  7.3] {MR4772266}. We now prove \eqref{prevac:estmuhat} for all $(s,q,\omega)\in[0,\infty)\times\R\times\mathbb{S}^2$. Since $|e^x-1|\leq |x|e^{|x|}$, we have $|\wh{\mu}+2|\lesssim |G(\omega)\wh{A}s|e^{Cs}\lesssim \lra{q}^{-1+C\eps}e^{Cs}$. Next, for $a+c>0$, we write  $\partial_q^a\partial_\omega^c\wh{\mu}$ as a linear combination (with real constant coefficients depending only on $a$ and $c$) of terms of the form
\fm{\exp(-\frac{1}{2}G \wh{A}s)\cdot\prod_{j=1}^k \partial_q^{a_j}\partial_\omega^{c_j}(-\frac{1}{2}G \wh{A}s),\qquad \sum a_*=a,\ \sum c_*=c,\ a_j+c_j>0\text{ for each }j.}
The estimates for $\wh{A}$ imply that \fm{|\exp(-\frac{1}{2}G \wh{A}s)\cdot\prod_{j=1}^k \partial_q^{a_j}\partial_\omega^{c_j}(-\frac{1}{2}G \wh{A}s)|\lesssim e^{Cs}\cdot \lra{q}^{-1-\sum a_j+C\eps}s^k\lesssim \lra{q}^{-1-a+C\eps}e^{Cs}.}
In summary, we have
$|\partial_q^a\partial_\omega^c\wh{\mu}|\lesssim_{a,c} \lra{q}^{-1-a+C_{a,c}\eps}e^{C_{a,c}s}$. That is, we obtain the first estimate in \eqref{prevac:estmuhat} with $b=0$.

Now, by \eqref{prevac:defhata}, we notice that
\eq{\label{prevac:estmuhat:pf}\partial_q^a\partial_s^b\partial_\omega^c\wh{U}_q&=\partial_q^a\partial_\omega^c((\frac{1}{2}G\wh{A})^b\wh{A}\cdot\exp(\frac{1}{2}G\wh{A}s)),\qquad a,b,c\geq 0;\\
\partial_q^a\partial_s^b\partial_\omega^c\wh{\mu}&=\partial_q^a\partial_\omega^c((-\frac{1}{2}G\wh{A})^b\cdot\exp(-\frac{1}{2}G\wh{A}s)),\qquad a,c\geq 0,b>0.}
Both of these rows are of the form 
\fm{\partial_q^a\partial_\omega^c[F_1(q,\omega)\cdot \exp(F_2(q,\omega)s)].}
Here for all $a,c\geq 0$ and $j=1,2$, we have $\partial_q^a\partial_\omega^cF_j=O_{a,c}(\lra{q}^{-b_j+C_{a,c}\eps})$ for some fixed constants $b_1>0$ and $b_2=1$. Following the same proof in the previous paragraph for $\wh{\mu}$, we obtain that for all $a,c\geq 0$ and $(s,q,\omega)$,
\fm{\partial_q^a\partial_\omega^c(\exp(F_2(q,\omega)s)-1)=O_{a,c}(\lra{q}^{-1-a+C_{a,c}\eps}e^{C_{a,c}s}).}
Thus, by Leibniz's rule, for all $a,c\geq 0$ and $(s,q,\omega)$, we have
\fm{&\partial_q^a\partial_\omega^c[F_1(q,\omega)\cdot \exp(F_2(q,\omega)s)]=\partial_q^a\partial_\omega^cF_1+\partial_q^a\partial_\omega^c[F_1 \cdot (\exp(F_2 s)-1)]\\
&=O_{a,c}(\lra{q}^{-b_1-a+C_{a,c}\eps}e^{C_{a,c}s}+\lra{q}^{-b_1-1-a+C_{a,c}\eps}e^{C_{a,c}s})=O_{a,c}(\lra{q}^{-b_1-a+C_{a,c}\eps}e^{C_{a,c}s}).}
To estimate the derivatives of $\wh{U}_q$, we use \eqref{prevac:estmuhat:pf} to set $F_1=(\frac{1}{2}G\wh{A})^b\wh{A}$ and $F_2=\frac{1}{2}G\wh{A}$. In this case, we have $b_1=b+1$, so we obtain the second row of \eqref{prevac:estmuhat}. To estimate the derivatives of $\wh{\mu}$, we use \eqref{prevac:estmuhat:pf} to set $F_1=(\frac{1}{2}G\wh{A})^b$ and $F_2=-\frac{1}{2}G\wh{A}$. In this case, we have $b_1=b$, so we obtain the first row of \eqref{prevac:estmuhat} with $b>0$.

To end the proof of \eqref{prevac:estmuhat}, we prove  the bounds for $\partial_s^b\partial_\omega^c\wh{U}$.
Note that $\wh{U}\equiv 0$ for $q\geq R$ and that for $q<R$,
\fm{|\partial_s^b\partial_\omega^c\wh{U}(s,q,\omega)|&\lesssim \int_{q}^R \lra{\rho}^{-1-b+C\eps}e^{Cs}\ d\rho.}
When $b=0$, the right side is $O(\eps^{-1}\lra{q}^{C\eps}e^{Cs})$. If $b\geq 1$, the right side is $O(e^{Cs})$. 

Finally, we check \eqref{prevac:estUhat}. Recall from \cite[Proposition 7.3]{MR4772266} that $\partial_\omega^c\wh{U}=O(e^{Cs})$ for $(s,q,\omega)\in\Omega'$. That is, $(s,\omega)\in [0,\infty)\times\mathbb{S}^2$ and $q>q_{(s,\eps)}:=\frac{1}{2}(e^{\delta/\eps}-e^{(s+\delta)/\eps})+2R$. Since $e^{\delta/\eps}-e^{(\delta+s)/\eps}-1=2q_{(s,\eps)}-4R-1<0$ and since $q_{(s,\eps)}\leq 2R$, we have \fm{\lra{q}/\lra{q_{(s,\eps)}}+\lra{q_{(s,\eps)}}/\lra{q}\lesssim_{R,C_0}e^{C_0s},\qquad \text{whenever } C_0e^{C_0s}(e^{\delta/\eps}-e^{(\delta+s)/\eps}-1)<q\leq q_{(s,\eps)}.} Since $\partial_\omega^c\wh{U}_q=O(\lra{q}^{-1+C\eps}e^{Cs})$ by \eqref{prevac:estmuhat}, by the mean value theorem we have
\fm{\partial_\omega^c\wh{U}(s,q,\omega;\eps)&=\partial_\omega^c\wh{U}(s,q_{(s,\eps)},\omega;\eps)+O(|q-q_{(s,\eps)}|\cdot \lra{q_{(s,\eps)}}^{-1+C\eps}e^{Cs})\\
&=O(e^{Cs})+O(\lra{q_{(s,\eps)}}\cdot \lra{q_{(s,\eps)}}^{-1+C\eps}\cdot e^{Cs})=O(e^{Cs})}
for all $C_0e^{C_0s}(e^{\delta/\eps}-e^{(\delta+s)/\eps}-1)<q\leq q_{(s,\eps)}$. This finishes the proof.
\end{proof}\rm

\rmk{\rm 
The bounds for $\wh{\mu}$ in \eqref{prevac:estmuhat}  are better than those for $\wt{\mu}$. For example, we have $\wh{\mu}_q=O(\lra{q}^{-2+C\eps}e^{Cs})$ but $\wt{\mu}_q=O(\lra{q}^{-1+C\eps}e^{Cs})$. Such a difference occurs because we have $\partial_q(-2)=0$ and $\partial_qA_1=O(\lra{q}^{-1+C\eps})$. This explains why we introduce this new map~$F$. }
\rm
 
\bigskip

It follows from \cite[Lemma 7.2]{MR4772266} that $\wh{q}(t,x;\eps):=F(\wt{q}(t,x;\eps),x/|x|;\eps)$ solves
\eq{\label{qhateqn}\wh{q}_t-\wh{q}_r=\wh{\mu}(\eps\ln t-\delta,\wh{q}(t,x;\eps),x/|x|;\eps),\qquad\text{in }\Omega,}
that $\wh{q}=r-t$ for $r-t\geq R$, and that 
\eq{\label{wtuwhuequal}\wt{u}(t,x;\eps)=\wh{u}(t,x;\eps):=\eps |x|^{-1}\wh{U}(\eps\ln t-\delta,\wh{q}(t,x;\eps),x/|x|;\eps),\qquad\text{in }\Omega.}
In the rest of this subsection, every function of $(s,q,\omega)$ induces a function of $(t,x)$ via the map $(t,x)\mapsto (\eps\ln t-\delta,\wh{q}(t,x;\eps),x/|x|)$. Note that this is different from the setting in Section~\ref{prevac:grs}. By \cite[Lemmas 7.4 and 7.7]{MR4772266}, we have $\wh{q}\in S^{0,1}$ and $\lra{r-t}\lesssim \lra{\wh{q}}t^{C\eps}$. Thus, we can apply Lemma \ref{c1lem6.1} with $\rho=\wh{q}$. For convenience, we list some results from \eqref{prevac:estmuhat} and \eqref{prevac:estUhat}:
\eq{\label{prevac:uapp:collection}
\wh{U},\eps\wh{U}_s\in S^{0,0};\ \wh{u}=\eps r^{-1}\wh{U}\in \eps S^{-1,0};\ \wh{\mu}+2\in S^{0,-1};
}
All the functions here are of $(t,x)$ via the map $(t,x)\mapsto (\eps\ln t-\delta,\wh{q}(t,x;\eps),x/|x|)=(s,q,\omega)$. We will use them in the next subsection. Note that $(t,x)\in\Omega$ does not imply that $(\eps\ln t-\delta,\wh{q},x/|x|)\in\Omega'$. Instead, by \cite[Lemma 7.4]{MR4772266}, we have $\lra{\wh{q}}/\lra{r-t}+\lra{r-t}/\lra{\wh{q}}\lesssim t^{C\eps}$ in $\Omega$. That is, in $\Omega\cap\{r-t\leq 2R\}$ we have $\wh{q}\leq 2R$ and
\fm{0<(2R+1-\wh{q})\sim\lra{\wh{q}}\lesssim\lra{r-t}t^{C\eps}\lesssim (t-e^{\delta/\eps}+1)t^{C\eps}. }
We hope that the estimate \eqref{prevac:estUhat} holds for such a $\wh{q}$, so we assume  $q>C_0e^{C_0s}(e^{\delta/\eps}-e^{(\delta+s)/\eps}-1)$ in Lemma \ref{prevac:lem:esthat}.

\subsection{Further computations}\label{secreview:furthercomputations}
In the next proposition, we present asymptotic formulas for $Z^I\wh{u}$ and $(\partial_t-\partial_r)Z^I\wh{u}$ in $\Omega$. We remind our readers of the definitions of $\eps^*S^{*,*}$ and $\eps^*S^{*,*}_{\ln}$ in Definition \ref{prel:defn:SD}.
\prop{\label{prop:trzwtu}
For each multiindex $I$, there exist  functions $U_I=U_I(s,q,\omega;\eps)$ and $A_I=A_I(q,\omega;\eps)$ defined for all $(s,q,\omega)\in[0,\infty)\times\R\times\mathbb{S}^2$ and $\eps\ll_I1$, such that 
\eq{\label{prop:trzwtu:c:est0}
r Z^I\wh{u}&=\eps U_I(\eps\ln t-\delta,|x|-t,x/|x|;\eps)+\eps S^{0,-1}_{\ln}+\eps S^{-1,1},\\
r(\partial_t-\partial_r)Z^I\wh{u}&=\eps A_I(|x|-t,x/|x|;\eps)+\eps S^{0,-2}_{\ln}+\eps S^{-1,0}.}
Here  $U_I$ and $A_I$ are defined inductively on $|I|$. First, we have $U_0=\wh{U}$ and $A_0=-2\wh{A}$. Moreover, for each $|I|>0$, by writing $Z^I=ZZ^{I'}$ with $|I'|=|I|-1$, we define $U_I$ and $A_I$ inductively by
\eq{\label{prop:trzwtu:inductformula:u}
U_I(s,q,\omega;\eps)=\left\{
\begin{array}{ll}
    \eps\partial_sU_{I'}+q\partial_qU_{I'}-U_{I'}, & Z^I=SZ^{I'};  \\
    (\omega_i\partial_{\omega_j}-\omega_j\partial_{\omega_i})U_{I'}, & Z^I=\Omega_{ij}Z^{I'},\ 1\leq i<j\leq 3;\\
     \eps\omega_i\partial_sU_{I'}-q\omega_i\partial_qU_{I'}+\partial_{\omega_i}U_{I'}-\omega_iU_{I'}, & Z^I=\Omega_{0i}Z^{I'},\ 1\leq i\leq 3;\\
   0, & Z^I=\partial Z^{I'};
\end{array}\right.
}
\eq{\label{prop:trzwtu:inductformula:a}
A_I(q,\omega;\eps)=\left\{
\begin{array}{ll}
    q\partial_qA_{I'}, & Z^I=SZ^{I'};  \\
    (\omega_i\partial_{\omega_j}-\omega_j\partial_{\omega_i})A_{I'}, & Z^I=\Omega_{ij}Z^{I'},\ 1\leq i<j\leq 3;\\
     -q\omega_i\partial_qA_{I'}+\partial_{\omega_i}A_{I'}-2\omega_iA_{I'}, & Z^I=\Omega_{0i}Z^{I'},\ 1\leq i\leq 3;\\
   0, & Z^I=\partial Z^{I'};
\end{array}\right.
}
Note that $\wh{U}\equiv \wh{A}\equiv 0$ for $q\geq R$, so by \eqref{prop:trzwtu:inductformula:u} and \eqref{prop:trzwtu:inductformula:a}, we have $U_I\equiv A_I\equiv 0$ for all $q\geq R$. For all $a,b,c\geq 0$, we have\eq{\label{trzwtu:esthat} \partial_q^a\partial_\omega^cA_I&=O_{a,c}(\lra{q}^{-1-a+C_{a,c}\eps});\\
\partial_q^a\partial_s^b \partial_\omega^c\partial_qU_I&=O_{a,b,c}(\lra{q}^{-1-a-b+C_{a,b,c}\eps}e^{C_{a,b,c}s});\\
\partial_s^b \partial_\omega^cU_I&=O_{b,c}((\eps^{-1}\lra{q}^{C_{b,c}\eps}1_{b=0}+1_{b>0})e^{C_{b,c}s})}
for all $(s,q,\omega)\in[0,\infty)\times\R\times\mathbb{S}^2$. For all $c\geq 0$, we have \eq{\label{trzwtu:estUhat} \partial_\omega^cU_I=O_{c,C_0}(e^{C_{c,C_0}s})}
for all $s>0$, $\omega\in\mathbb{S}^2$, and $q>C_0e^{C_0s}(e^{\delta/\eps}-e^{(\delta+s)/\eps}-1)$ for some constant $C_0>1$. Finally, for all $a,b,c\geq 0$, we have
\eq{\label{prop:trzwtu:derivative} \partial_q^a\partial_s^b\partial_\omega^c(2\partial_qU_I+A_I)(s,q,\omega;\eps)
&=O_{a,b,c}(\lra{q}^{-\max\{1,b\}-1-a+C_{a,b,c}\eps}e^{C_{a,b,c}s})}
for all $(s,q,\omega)\in[0,\infty)\times\R\times\mathbb{S}^2$.
}

\rmk{\rm 
Since $\wh{u}=\eps r^{-1}\wh{U}(s,\wh{q},\omega)$ is defined explicitly by \eqref{wtuwhuequal}, one can compute $Z^I\wh{u}$ by the chain rule and Leibniz's rule. This has been done in, e.g., the proof of Lemma \ref{c1lem6.1}. Thus, we can obtain an identity for $Z^I\wh{u}$. Similarly for $(\partial_t-\partial_r)Z^I\wh{u}$. In this paper, we prefer the asymptotic formulas   \eqref{prop:trzwtu:c:est0} over these identities because \eqref{prop:trzwtu:c:est0} are much simpler.

Let us compute $Z\wh{u}$ as an example. It suffices to compute $Z(r\wh{u})=Z(\eps\wh{U}(\eps\ln t-\delta,\wh{q},\omega))$ because of the commutator $[Z,r]=Zr$.
By the chain rule, we obtain the following term in the expansion of $Z(r\wh{u})$:
\eq{\label{whuqepslnt}\eps\wh{U}_q(\eps\ln t-\delta,\wh{q},\omega)\cdot Z\wh{q}.}
Recall that $\wh{q}$ is a solution to the transport equation \eqref{qhateqn}. There is no explicit formula for $Z\wh{q}$, so we cannot further simplify \eqref{whuqepslnt}.  In Proposition \ref{prop:trzwtu}, we use an alternative way. We first show that 
$Z(r\wh{u})$ equals $Z(\eps\wh{U}(\eps\ln t-\delta,r-t,\omega))$ up to a remainder term. This remainder term is good because we expect $\wh{q}\approx r-t$ to some extent; see Lemma \ref{prevac:lem:qrt} below. By the chain rule, we replace \eqref{whuqepslnt} with
\fm{\eps\wh{U}_q(\eps\ln t-\delta,r-t,\omega)\cdot Z(r-t).}
We can compute this expression explicitly because we have $S(r-t)=r-t$, $\Omega_{0i}(r-t)=(t-r)\omega_i$, and $Z(r-t)=0$ for other $Z$. As a result, we obtain a simpler asymptotic formula for $Z\wh{u}$.}

\rmk{\rm As commented in Remark \ref{mthm_rmk1.3}, if $Z^I$ contains a translation, we have $A_I\equiv 0$. In this case, we can still approximate $Z^I\wh{u}$. We achieve this goal by  following Remark \ref{mthm_rmk1.3}, and we remark that the key is the identity \eqref{mthm_rmk1.3_formula}.}\rm

\bigskip

The key step in the proof of Proposition \ref{prop:trzwtu} is to show the following estimates:
\eq{\label{pf:trzwtu:keyest}r\wh{u}-\eps\wh{U}(\eps\ln t-\delta,r-t,\omega;\eps)&\in\eps S^{0,-1}_{\ln}+\eps S^{-1,1},\\
(\partial_t-\partial_r)(r\wh{u})+2\eps\wh{A}(r-t,\omega;\eps)&\in\eps S^{0,-2}_{\ln}+\eps S^{-1,0}.}
Recall that $r=|x|$ and $\omega=x/|x|$.
Proving these estimates requires several estimates for the approximate optical function $\wh{q}$; see Section \ref{sec:trzwtu:step1}.

In Section \ref{sec:trzwtu:step2}, we extend \eqref{pf:trzwtu:keyest} by showing that for each multiindex $I$, we have
\eq{\label{pf:trzwtu:keyest2}rZ^I\wh{u}&=\eps U_I(\eps\ln t-\delta,r-t,\omega;\eps)+\eps S^{0,-1}_{\ln}+\eps S^{-1,1},\\(\partial_t-\partial_r)(rZ^I\wh{u})&=\eps A_I(r-t,\omega;\eps)+\eps S^{0,-2}_{\ln}+\eps S^{-1,0}.}
Here the $U_I$ and $A_I$ are defined in the statement of Proposition \ref{prop:trzwtu}: we set $(U_0,A_0)=(\wh{U},-2\wh{A})$ and define $(U_I,A_I)$ inductively by \eqref{prop:trzwtu:inductformula:u} and \eqref{prop:trzwtu:inductformula:a} for $|I|>0$. We can induct on $|I|$ to check  \eqref{trzwtu:esthat} and \eqref{trzwtu:estUhat}; the details of the proofs are omitted. We also notice that
$[\partial_t-\partial_r,r]Z^I\wh{u}=-Z^I\wh{u}=O(\eps t^{-1+C\eps})$ by \eqref{prevac:wtu}, so we obtain \eqref{prop:trzwtu:c:est0}.

Finally, we prove \eqref{prop:trzwtu:derivative}. By setting $B_I=2\partial_qU_I+A_I$ for each $I$, we note that the $B_I$'s satisfy the induction formulas
\fm{
B_I=\left\{
\begin{array}{ll}
    \eps\partial_sB_{I'}+q\partial_qB_{I'}, & Z^I=SZ^{I'};  \\
    (\omega_i\partial_{\omega_j}-\omega_j\partial_{\omega_i})B_{I'}, & Z^I=\Omega_{ij}Z^{I'},\ 1\leq i<j\leq 3;\\
     \eps\omega_i\partial_sB_{I'}-q\omega_i\partial_qB_{I'}+\partial_{\omega_i}B_{I'}-2\omega_iB_{I'}, & Z^I=\Omega_{0i}Z^{I'},\ 1\leq i\leq 3;\\
   0, & Z^I=\partial Z^{I'}.
\end{array}\right.
}
We also have $B_0=2\wh{A}(\exp(\frac{1}{2}G(\omega)\wh{A}s)-1)$. Since  $\exp(\frac{1}{2}G(\omega)\wh{A}s)-1$ has essentially the same form as $\wh{\mu}+2$, by following the proof of \eqref{prevac:estmuhat}, we obtain for all $a,b,c\geq 0$ and all $(s,q,\omega)\in[0,\infty)\times\R\times\mathbb{S}^2$,
\fm{\partial_q^a\partial_s^b\partial_\omega^c (\exp(\frac{1}{2}G(\omega)\wh{A}s)-1)=O_{a,b,c}(\lra{q}^{-\max\{1,b\}-a+C_{a,b,c}\eps}e^{C_{a,b,c}s}).}
It follows from Leibniz's rule that
\fm{\partial_q^a\partial_s^b\partial_\omega^cB_0&=O_{a,b,c}(\lra{q}^{-\max\{1,b\}-1-a+C_{a,b,c}\eps}e^{C_{a,b,c}s}).}
Now, \eqref{prop:trzwtu:derivative} follows by induction on $|I|$. We again omit the details of the proof.

\subsubsection{Proof of \eqref{pf:trzwtu:keyest}}\label{sec:trzwtu:step1}
By the definition of $\wh{u}$, in $\Omega$ we have 
\fm{r\wh{u}&=\eps\wh{U},\\
(\partial_t-\partial_r)(r\wh{u})&=(\partial_t-\partial_r)(\eps\wh{U})=\eps^2t^{-1}\wh{U}_s+\eps\wh{U}_q\wh{\mu}=\eps^2t^{-1}\wh{U}_s-2\eps\wh{A}.}
Here $\wh{U},\wh{U}_s,\wh{A}$ are functions of $(t,x)$ via the map $(t,x)\mapsto (\eps\ln t-\delta,\wh{q}(t,x),x/|x|)$. By \eqref{prevac:uapp:collection} at the end of Section \ref{secreview:secapp}, we have $\eps \wh{U}_s\in S^{0,0}$. Thus, we have
\fm{(\partial_t-\partial_r)(r\wh{u})&=-2\eps\wh{A}+\eps S^{-1,0}.}
It remains to show that \eq{\label{pf:trzwtu:keyest3}\eps \wh{U}(\eps\ln t-\delta,\wh{q},\omega;\eps)-\eps \wh{U}(\eps\ln t-\delta,r-t,\omega;\eps)&\in \eps S^{0,-1}_{\ln}+\eps S^{-1,1},\\
\eps \wh{A}(\wh{q},\omega;\eps)-\eps \wh{A}(r-t,\omega;\eps)&\in \eps S^{0,-2}_{\ln}+\eps S^{-1,0}.}

We first prove a lemma for $\wh{q}$. This lemma both relies on and improves \cite[Lemma 7.7]{MR4772266}. 
\lem{\label{prevac:lem:qrt}We have $\wh{q}-r+t\in S^{0,0}_{\ln}$. That is, for each $I$, we have $Z^I(\wh{q}-r+t)=O((1+\ln\lra{r-t})t^{C\eps})$ in $\Omega$.}
\begin{proof}
For $r-t\geq R$, we have $q=r-t$. It remains to prove $Z^I(\wh{q}-r+t)=O((1+\ln\lra{r-t})t^{C\eps})$ for each $I$ in $\Omega\cap\{r-t<2R\}$.
Before we start the proof, we recall from \cite[Lemma 7.4]{MR4772266}  that $t^{-C\eps}\lesssim \lra{\wh{q}}/\lra{r-t}\lesssim t^{C\eps}$ in $\Omega$. Thus, every power of $\lra{\wh{q}}$ can be controlled by $\lra{r-t}t^{C\eps}$ and vice versa.

We first apply Lemma \ref{sec:prel:commutetrz} to compute $(\partial_t-\partial_r)Z^I(\wh{q}-r+t)$ for each multiindex $I$. For each integer $N\geq 0$, in $\Omega$ we have
\eq{\label{lem:qrt:fff}&\sum_{|I|\leq N}|(\partial_t-\partial_r)Z^I(\wh{q}-r+t)|\\
&\lesssim \sum_{|I|\leq N}|Z^I(\wh{\mu}+2)|+\sum_{|J|<N}\kh{|Z^J(\wh{\mu}+2)|+\sum_{j=1}^3|(\partial_j+\omega_j\partial_t)Z^J(\wh{q}-r+t)|}\\
&\lesssim \sum_{|I|\leq N}|Z^I(\wh{\mu}+2)|+\lra{r+t}^{-1}\sum_{|I|\leq N}|Z^I(\wh{q}-r+t)|.}
Here we recall that $\wh{\mu}=\wh{q}_t-\wh{q}_r$. In the last estimate, we apply Lemma \ref{c1l2.1} to get the factor $\lra{r+t}^{-1}$. Meanwhile, we recall  that $\wh{\mu}+2\in S^{0,-1}$. In fact, the results in \eqref{prevac:uapp:collection} hold because it is already known that $\wh{q}\in S^{0,1}$ from \cite[Lemma 7.7]{MR4772266}. Thus, we have $\sum_{|I|\leq N}|Z^I(\wh{\mu}+2)|\lesssim \lra{r-t}^{-1}t^{C\eps}$ in $\Omega$.

We can now use Gronwall's inequality to finish the proof. Fix an integer $N\geq 0$ and a point $(t_0,x_0)\in \Omega\cap\{|x_0|-t_0<2R\}$. For each $\tau\in[\frac{t_0+|x_0|}{2}-R,t_0]$, we set
\fm{z_N(\tau)&=\sum_{|I|\leq N}|Z^I(\wh{q}-r+t)(\tau,(t_0+|x_0|-\tau)\frac{x_0}{|x_0|};\eps)|.}
Note that $(\tau,(t_0+|x_0|-2\tau)\frac{x_0}{|x_0|})\in \Omega\cap\{r-t<2R\}$ for all $\tau\in[\frac{t_0+|x_0|}{2}-R,t_0]$. Moreover, given a function $h(t,x)$, for $\wt{z}(\tau)=h(\tau,(|x_0|+t_0-\tau)\cdot \frac{x_0}{|x_0|})$, we have $\frac{d}{d\tau}\wt{z}=(h_t-h_r)(\tau,(|x_0|+t_0-\tau)\cdot \frac{x_0}{|x_0|})$. Thus, the inequality \eqref{lem:qrt:fff}  implies that for each $\tau_0\in[\frac{t_0+|x_0|}{2}-R,t_0]$,
\fm{z_N(\tau_0)&\lesssim z_N(\frac{t_0+|x_0|}{2}-R)+\int_{\frac{t_0+|x_0|}{2}-R}^{\tau_0}\lra{|x_0|+t_0-2\tau}^{-1}\tau^{C\eps}+\lra{|x_0|+t_0}^{-1}z_N(\tau)\ d\tau\\
&\lesssim\lra{|x_0|+t_0}^{-1}\int_{\frac{t_0+|x_0|}{2}-R}^{\tau_0}z_N(\tau)\ d\tau+ t_0^{C\eps}(1+\ln\lra{|x_0|-t_0}).}
To obtain the second estimate, we notice that $(\frac{t_0+|x_0|}{2}-R,(\frac{t_0+|x_0|}{2}-R)\frac{x_0}{|x_0|})\in\{(t,x):\ |x|-t>R\}$. Thus,  we have $\wh{q}=r-t$  near this point and thus $z_N(\frac{t_0+|x_0|}{2}-R)=0$. Moreover, for $\frac{t_0+|x_0|}{2}-R\leq \tau\leq t_0$, we have $\tau^{C\eps}\leq t_0^{C\eps}$ and $t+x_0-2\tau\leq 2R$. It follows that
\fm{\int_{\frac{t_0+|x_0|}{2}-R}^{\tau_0}\lra{|x_0|+t_0-2\tau}^{-1}\tau^{C\eps}\ d\tau&\lesssim t_0^{C\eps}\int_{\frac{t_0+|x_0|}{2}-R}^{\tau_0}(2\tau-|x_0|-t_0+2R+1)^{-1}\ d\tau\\
&\lesssim t_0^{C\eps}\ln(t_0-|x_0|+2R+1)\lesssim t_0^{C\eps}(1+\ln\lra{t_0-|x_0|}).}
Thus, by Gronwall's inequality, we have
\fm{z_N(t_0)&\lesssim t_0^{C\eps}(1+\ln\lra{t_0-|x_0|})\exp(\int_{\frac{t_0+|x_0|}{2}-R}^{t_0}\lra{|x_0|+t_0}^{-1}\ d\tau)\\
&\lesssim t_0^{C\eps}(1+\ln\lra{t_0-|x_0|})\exp(C\cdot \frac{\frac{t_0-|x_0|}{2}+R}{|x_0|+t_0+1})\lesssim t_0^{C\eps}(1+\ln\lra{t_0-|x_0|}).}
\end{proof}
\rm

\bigskip

Using Lemma \ref{prevac:lem:qrt} and the bounds for $\wh{U},\wh{A}$, we now show \eqref{pf:trzwtu:keyest3} which implies \eqref{pf:trzwtu:keyest}.

\lem{\label{prevac:lem:ahatqrt}Fix an integer $N\geq 0$. Then, for $\eps\ll_N1$ and all $(t,x)\in\Omega$, we have \fm{\sum_{|I|\leq N}|Z^I(\wh{U}(\eps\ln t-\delta,\wh{q},\omega;\eps)-\wh{U}(\eps\ln t-\delta,r-t,\omega;\eps))|\lesssim (1+\ln\lra{r-t})\lra{r-t}^{-1}t^{C\eps},\\
\sum_{|I|\leq N}|Z^I(\wh{A}(\wh{q},\omega)-\wh{A}(r-t,\omega))|\lesssim (1+\ln\lra{r-t})\lra{r-t}^{-2}t^{C\eps}.}}
\begin{proof}
Since $\wh{q}=r-t$ for $r-t\geq R$, it remains to prove this lemma in $\Omega\cap\{r-t<2R\}$. We will prove by induction on $|I|$.

When $|I|=0$, by the mean value theorem, we have (here we set $s=\eps\ln t-\delta$ and omit the parameter $\eps$)
\fm{&|\wh{U}(s,\wh{q},\omega)-\wh{U}(s,r-t,\omega)|\lesssim |\wh{q}-r+t|\cdot \sup_{\rho\text{ lies between }\wh{q},r-t}|\wh{U}_q(s,\rho,\omega)|\\
&\lesssim (1+\ln\lra{r-t})t^{C\eps}\cdot \lra{r-t}^{-1+C\eps}t^{C\eps}\lesssim (1+\ln\lra{r-t})\lra{r-t}^{-1}t^{C\eps},\\
&|\wh{A}(\wh{q},\omega)-\wh{A}(r-t,\omega)|\lesssim |\wh{q}-r+t|\cdot \sup_{\rho\text{ lies between }\wh{q},r-t}|\wh{A}_q(\rho,\omega)|\\
&\lesssim (1+\ln\lra{r-t})t^{C\eps}\cdot \lra{r-t}^{-2+C\eps}t^{C\eps}\lesssim (1+\ln\lra{r-t})\lra{r-t}^{-2}t^{C\eps}.}
Here we use \eqref{prevac:estahat}, \eqref{prevac:estmuhat}, and Lemma \ref{prevac:lem:qrt}. Besides, we also use the fact that if $\rho$ lies between $\wh{q}$ and $r-t$, then $\lra{r-t}\lesssim\lra{\rho}t^{C\eps}$. In fact, since $\wh{q}<2R$ and $r-t<2R$ in $\Omega\cap\{r-t<2R\}$, we also have $\rho<2R$. And since $t^{-C\eps}\lesssim\lra{\wh{q}}/\lra{r-t}\lesssim t^{C\eps}$, we have
\fm{\lra{\rho}^{-1}&\sim (1+2R-\rho)^{-1}\leq (1+2R-\wh{q})^{-1}+(1+2R-r+t)^{-1}\\
&\lesssim \lra{\wh{q}}^{-1}+\lra{r-t}^{-1}\lesssim \lra{r-t}^{-1}t^{C\eps}.}
We thus finish the proof when $|I|=0$.

In general, we fix a multiindex $I$ with $|I|>0$. By the chain rule and Leibniz's rule, we can write $Z^I(\wh{U}(s,\wh{q},\omega)-\wh{U}(s,r-t,\omega))$ as a linear combination (with real constant coefficients depending on $I$) of terms of the form
\fm{&(\partial_q^a\partial_s^b\partial_\omega^c \wh{U})(s,\wh{q},\omega)\cdot \prod_{j=1}^b Z^{I_j}(\eps\ln t-\delta)\cdot\prod_{j=1}^a Z^{J_j}\wh{q}\cdot\prod_{j=1}^cZ^{K_j}\omega
\\&-(\partial_q^a\partial_s^b\partial_\omega^c \wh{U})(s,r-t,\omega)\cdot \prod_{j=1}^b Z^{I_j}(\eps\ln t-\delta)\cdot\prod_{j=1}^a Z^{J_j}(r-t)\cdot\prod_{j=1}^cZ^{K_j}\omega.}
Here $a+b+c>0$, the sum of all the multiindices is $|I|$, and each of these multiindices is nonzero. 
This term equals 
\fm{&[(\partial_q^a\partial_s^b\partial_\omega^c \wh{U})(s,\wh{q},\omega)-(\partial_q^a\partial_s^b\partial_\omega^c \wh{U})(s,r-t,\omega)]\cdot\prod_{j=1}^bZ^{I_j}(\eps\ln t)\cdot\prod_{j=1}^a Z^{J_j}\wh{q}\cdot\prod_{j=1}^cZ^{K_j}\omega\\
& +(\partial_s^b\partial_q^a\partial_\omega^c \wh{U})(s,r-t,\omega)\cdot \kh{\sum_{k=1}^a\prod_{j=1}^k Z^{J_j}\wh{q}\cdot Z^{J_k}(\wh{q}-r+t)\cdot \prod_{j=k+1}^a Z^{J_j}(r-t)}\\&\hspace{25em}\cdot\prod_{j=1}^bZ^{I_j}(\eps\ln t)\cdot\prod_{j=1}^cZ^{K_j}\omega.}
If $a=0$, then the second row vanishes. We claim that this term is $O((1+\ln\lra{r-t})\lra{r-t}^{-1}t^{C\eps})$.
In fact, by \eqref{prevac:estmuhat}, we follow a similar proof in the case $|I|=0$ to show that \fm{(\partial_q^a\partial_s^b\partial_\omega^c \wh{U})(s,\wh{q},\omega)-(\partial_q^a\partial_s^b\partial_\omega^c \wh{U})(s,r-t,\omega)=O((1+\ln\lra{r-t})\lra{r-t}^{-1-a}t^{C\eps}).} Moreover, by \eqref{prevac:estmuhat} and \eqref{prevac:estUhat}, we have $(\partial_q^a\partial_s^b\partial_\omega^c \wh{U})(s,r-t,\omega)=O(\lra{r-t}^{-a}t^{C\eps})$ in $\Omega\cap\{r-t<2R\}$. We also have $Z^J\omega=O(1)$,  $Z^J(\wh{q},r-t)=O(\lra{r-t}t^{C\eps})$, and $Z^J(\wh{q}-r+t)=O((1+\ln\lra{r-t})t^{C\eps})$. Combine all these estimates and we obtain the estimates for $\wh{U}$ in the lemma. The estimates for $\wh{A}$ can be proved similarly.
\end{proof}\rm

\subsubsection{Proof of \eqref{pf:trzwtu:keyest2}}\label{sec:trzwtu:step2}

We  induct on $|I|$. The case when $|I|=0$ has been proved; see \eqref{pf:trzwtu:keyest}. Now, fix a multiindex $I$ with $|I|>0$ and write $Z^I=ZZ^{I'}$.  We have
\eq{\label{trzwtu:step2:est0}rZ^I\wh{u}&= Z(rZ^{I'}\wh{u})+[r,Z]Z^{I'}\wh{u}\\
&=\eps (\partial_sU_{I'})\cdot Z(\eps\ln t)+\eps (\partial_qU_{I'})\cdot Z(r-t)\\
&\quad +\sum_{k=1}^3\eps (\partial_{\omega_k}U_{I'})\cdot Z\omega_k-Zr\cdot Z^{I'}\wh{u}+\eps S^{0,-1}_{\ln}+\eps S^{-1,1},}
\eq{\label{trzwtu:step2:est}(\partial_t-\partial_r)(rZ^I\wh{u})&=Z[(\partial_t-\partial_r)(rZ^{I'}\wh{u})]+[\partial_t-\partial_r,Z](rZ^{I'}\wh{u})+(\partial_t-\partial_r)([r,Z]Z^{I'}\wh{u})\\
&=\eps (\partial_qA_{I'})\cdot Z(r-t)+\sum_{k=1}^3\eps (\partial_{\omega_k}A_{I'})\cdot Z\omega_k\\
&\quad+[\partial_t-\partial_r,Z](rZ^{I'}\wh{u})-(\partial_t-\partial_r)(Zr\cdot Z^{I'}\wh{u})+\eps S^{0,-2}_{\ln}+\eps S^{-1,0}.}
On the right sides, every function of $(s,q,\omega)$ is evaluated at $(s,q,\omega)=(\eps\ln t-\delta,r-t,\omega)$. For example, $\partial_sU_{I'}=(\partial_sU_{I'})(\eps\ln t-\delta,r-t,\omega;\eps)$; similarly for $U_{I'},A_{I'},\partial_q(U_{I'},A_{I'}),\partial_\omega (U_{I'}, A_{I'})$.
We apply the induction hypotheses and the chain rule to compute $Z(rZ^{I'}\wh{u})$ and $Z((\partial_t-\partial_r)(rZ^{I'}\wh{u}))$.
Now we can prove \eqref{pf:trzwtu:keyest2}.

\lem{Fix a multiindex $I$ with $|I|>0$ and $Z^I=ZZ^{I'}$. Assuming that \eqref{pf:trzwtu:keyest2} holds with $I$ replaced by $I'$, we have \eqref{pf:trzwtu:keyest2} for $I$. Here the functions $U_I$ and $A_I$ are defined in the statement of Proposition \ref{prop:trzwtu}.}
\begin{proof}
In this proof, every function of $(s,q,\omega)$ (such as $U_{I'}$ and $A_{I'}$) is evaluated at $(s,q,\omega)=(\eps\ln t-\delta,r-t,\omega)$.

Since the $U_I$'s and $A_I$'s are defined inductively by \eqref{prop:trzwtu:inductformula:u} and \eqref{prop:trzwtu:inductformula:a}, we can prove \eqref{trzwtu:esthat} and \eqref{trzwtu:estUhat} without proving \eqref{prop:trzwtu:c:est0}.
Thus, for each multiindex $I$, we have \eq{\label{step2:collection}U_I,\eps \partial_sU_I,\partial_\omega U_I\in S^{0,0};\quad \partial_qU_I,\partial_\omega A_I\in S^{0,-1};\quad \partial_qA_I\in S^{0,-2}.}
Here we apply Lemma \ref{c1lem6.1} (with $\rho=r-t$), \eqref{trzwtu:esthat}, and \eqref{trzwtu:estUhat}. 

(a) Suppose that $Z^I=\partial Z^{I'}$. We have $\partial (\eps\ln t),\partial\omega\in \eps S^{-1,0}$ and $\partial (r-t)\in S^{0,0}$. By \eqref{step2:collection} and Lemma \ref{c1lem6}, we have
\fm{& \eps (\partial_sU_{I'})\cdot \partial(\eps\ln t)+\eps (\partial_qU_{I'}) \cdot \partial(r-t)+\sum_{k=1}^3\eps (\partial_{\omega_k}U_{I'}) \cdot \partial\omega_k\\
&\in S^{0,0}\cdot\eps S^{-1,0}+\eps S^{0,-1}\cdot S^{0,0}+\eps S^{0,0}\cdot S^{-1,0}\subset \eps S^{0,-1},\\
& \eps (\partial_qA_{I'}) \cdot \partial(r-t)+\sum_{k=1}^3\eps (\partial_{\omega_k}A_{I'}) \cdot \partial\omega_k\in \eps S^{0,-2}\cdot S^{0,0}+\eps S^{0,-1}\cdot S^{-1,0}\subset \eps S^{0,-2}.}
Moreover, we have $[\partial_t-\partial_r,\partial_i]=r^{-1}(\partial_i-\omega_i\partial_r)=S^{-1,0}\cdot \partial$ (by $S^{-1,0}\cdot \partial$ we mean $f^\alpha(t,x) \partial_\alpha$ where $f^\alpha\in S^{0,0}$), $\partial_ir=\omega_i\in S^{0,0}$ and $[\partial_t-\partial_r,\partial_t]=\partial_tr=0$. Since $\wh{u}\in \eps S^{-1,0}$, by Lemma \ref{c1lem6} we have \fm{\partial r\cdot Z^{I'}\wh{u}&\in S^{0,0}\cdot \eps S^{-1,0}\subset \eps S^{-1,0},\\ [\partial_t-\partial_r,\partial](rZ^{I'}\wh{u})&\in S^{-1,0}\cdot \partial(\eps S^{0,0})\subset \eps S^{-1,-1},\\
(\partial_t-\partial_r)(\partial r\cdot Z^{I'}\wh{u})&\in \partial (S^{0,0}\cdot \eps S^{-1,0})\subset \eps S^{-1,-1}.}
In summary, by \eqref{trzwtu:step2:est0} and \eqref{trzwtu:step2:est}, we have $ r\partial Z^{I'}\wh{u}\in \eps S^{0,-1}_{\ln}+\eps S^{-1,1}$ and $(\partial_t-\partial_r)(r\partial Z^{I'}\wh{u})\in \eps S^{0,-2}_{\ln}+\eps S^{-1,0}$, i.e.\ $U_I\equiv A_I\equiv 0$.

(b) Suppose that $Z^I=SZ^{I'}$. In this case, we have $S(\eps\ln t)=\eps$, $S(r-t)=r-t$, $S\omega=0$, $[\partial_t-\partial_r,S]=\partial_t-\partial_r$, and $Sr=r$. By \eqref{trzwtu:step2:est0} and \eqref{trzwtu:step2:est}, we have
\fm{rSZ^{I'}\wh{u}&=\eps^2 (\partial_sU_{I'}) +\eps (\partial_qU_{I'}) \cdot (r-t) -r\cdot Z^{I'}\wh{u}+\eps S^{0,-1}_{\ln}+\eps S^{-1,1}\\
&=\eps (\eps \partial_sU_{I'}+q\partial_qU_{I'}-U_{I'}) +\eps S^{0,-1}_{\ln}+\eps S^{-1,1},\\
(\partial_t-\partial_r)(rSZ^{I'}\wh{u})
&=\eps (\partial_qA_{I'}) \cdot (r-t)+\eps S^{0,-2}_{\ln}+\eps S^{-1,0}=\eps(q\partial_qA_{I'})+\eps S^{0,-2}_{\ln}+\eps S^{-1,0}.}
That is, $U_I=\eps\partial_sU_{I'}+q\partial_qU_{I'}-U_{I'}$ and $A_I=q\partial_qA_{I'}$. In these  estimates, we apply the induction hypotheses to compute $rZ^{I'}\wh{u}$. We also notice that $[S,r(\partial_t-\partial_r)]=0$.

(c) Suppose that $Z^I=\Omega_{ij}Z^{I'}$. In this case, we have $\Omega_{ij}(\eps\ln t)=\Omega_{ij}(r-t)=\Omega_{ij}r=0$, $\Omega_{ij}\omega_k=\omega_i\delta_{jk}-\omega_j\delta_{ik}$, and $[\partial_t-\partial_r,\Omega_{ij}]=0$. By \eqref{trzwtu:step2:est0} and \eqref{trzwtu:step2:est}, we have
\fm{r\Omega_{ij}Z^{I'}\wh{u}&=\sum_{k=1}^3\eps (\partial_{\omega_k}U_{I'}) \cdot (\omega_{i}\delta_{jk}-\omega_j\delta_{ik})+\eps S^{0,-1}_{\ln}+\eps S^{-1,1}\\
&=\eps ((\omega_i\partial_{\omega_j}-\omega_j\partial_{\omega_i})U_{I'}) +\eps S^{0,-1}_{\ln}+\eps S^{-1,1},\\
(\partial_t-\partial_r)(r\Omega_{ij}Z^{I'}\wh{u})&=\sum_{k=1}^3\eps (\partial_{\omega_k}A_{I'}) \cdot (\omega_i\delta_{jk}-\omega_j\delta_{ik})+\eps S^{0,-2}_{\ln}+\eps S^{-1,0}\\
&=\eps ((\omega_i\partial_{\omega_j}-\omega_j\partial_{\omega_i})A_{I'}) +\eps S^{0,-2}_{\ln}+\eps S^{-1,0}.}
That is, $U_I=(\omega_i\partial_{\omega_j}-\omega_j\partial_{\omega_i})U_{I'}$ and  $A_I=(\omega_i\partial_{\omega_j}-\omega_j\partial_{\omega_i})A_{I'}$.

(d) Suppose that $Z^I=\Omega_{0i}Z^{I'}$. In this case, we have $\Omega_{0i}(\eps\ln t)=\eps \omega_i(1+\frac{r-t}{t}) $, $\Omega_{0i}(r-t)=-(r-t)\omega_i$, $\Omega_{0i}\omega_k=(1-\frac{r-t}{r})(\delta_{ik}-\omega_i\omega_k)$, $[\partial_t-\partial_r,\Omega_{0i}]=-\omega_i(\partial_t-\partial_r)-\frac{r+t}{r^2}\sum_{j=1}^3\omega_j\Omega_{ij}$, and $\Omega_{0i}r=t\omega_i=(1-\frac{r-t}{r})r\omega_i$. By \eqref{trzwtu:step2:est0} and \eqref{trzwtu:step2:est}, we have
\fm{&r\Omega_{0i}Z^{I'}\wh{u}\\
&=\eps (\partial_sU_{I'}) \cdot \eps\omega_i(1+\frac{r-t}{t})-\eps (\partial_qU_{I'}) \cdot (r-t)\omega_i\\
&\quad +\sum_{k=1}^3\eps (\partial_{\omega_k}U_{I'}) \cdot (1-\frac{r-t}{r})(\delta_{ik}-\omega_i\omega_k)-(1-\frac{r-t}{r})r\omega_i\cdot Z^{I'}\wh{u}+\eps S^{0,-1}_{\ln}+\eps S^{-1,1}\\
&=\eps(\eps\omega_i\partial_sU_{I'}-q\omega_i\partial_qU_{I'}+\partial_{\omega_i}U_{I'}-\omega_iU_{I'})\\
&\quad+\eps^2 (\partial_sU_{I'})\omega_i \cdot \frac{r-t}{t}-\eps (\partial_{\omega_i}U_{I'}) \cdot \frac{r-t}{r}+(r-t)\omega_iZ^{I'}\wh{u}+\eps S^{0,-1}_{\ln}+\eps S^{-1,1},}
\fm{&(\partial_t-\partial_r)(r\Omega_{0i}Z^{I'}\wh{u})\\
&=-\eps (\partial_qA_{I'}) \cdot (r-t)\omega_i+\sum_{k=1}^3\eps (\partial_{\omega_k}A_{I'}) \cdot (1-\frac{r-t}{r})(\delta_{ik}-\omega_i\omega_k)\\
&\quad-\omega_i(\partial_t-\partial_r)((r+t)Z^{I'}\wh{u})-\frac{r+t}{r^2}\sum_{j=1}^3\omega_j\Omega_{ij}(rZ^{I'}\wh{u})+\eps S^{0,-2}_{\ln}+\eps S^{-1,0}\\
&=\eps (-q\omega_i\partial_qA_{I'}+\partial_{\omega_i}A_{I'}-2\omega_iA_{I'}) -\eps (\partial_{\omega_i}A_{I'}) \cdot  \frac{r-t}{r}\\
&\quad+\omega_i(\partial_t-\partial_r)((r-t)Z^{I'}\wh{u})-\frac{r+t}{r^2}\sum_{j=1}^3x_j\Omega_{ij}Z^{I'}\wh{u}+\eps S^{0,-2}_{\ln}+\eps S^{-1,0}.} By \eqref{step2:collection}, we have
\fm{\eps^2 (\partial_sU_{I'})\omega_i \cdot \frac{r-t}{t}-\eps (\partial_{\omega_i}U_{I'}) \cdot \frac{r-t}{r}+(r-t)\omega_iZ^{I'}\wh{u}&\in \eps S^{-1,1},\\
-\eps (\partial_{\omega_i}A_{I'}) \cdot  \frac{r-t}{r}+\omega_i(\partial_t-\partial_r)((r-t)Z^{I'}\wh{u})-\frac{r+t}{r^2}\sum_{j=1}^3x_j\Omega_{ij}Z^{I'}\wh{u}&\in \eps S^{-1,0}.}
In summary, we have
\fm{r\Omega_{0i}Z^{I'}\wh{u}&=\eps(\eps\omega_i\partial_sU_{I'}-q\omega_i\partial_qU_{I'}+\partial_{\omega_i}U_{I'}-\omega_iU_{I'})+\eps S^{0,-1}_{\ln}+\eps S^{-1,1},\\
(\partial_t-\partial_r)(r\Omega_{0i}Z^{I'}\wh{u})
&=\eps (-q\omega_i\partial_qA_{I'}+\partial_{\omega_i}A_{I'}-2\omega_iA_{I'})+\eps S^{0,-2}_{\ln}+\eps S^{-1,0}.}
That is, $U_I=\eps\omega_i\partial_sU_{I'}-q\omega_i\partial_qU_{I'}+\partial_{\omega_i}U_{I'}-\omega_iU_{I'}$ and $A_I=-q\omega_i\partial_qA_{I'}+\partial_{\omega_i}A_{I'}-2\omega_iA_{I'}$.
\end{proof}\rm

\subsection{Gauge independence}\label{secreview:gauage}
In Section \ref{sec:prevac:opticalfunction}, we construct the optical function $q=q(t,x;\eps)$ by solving the eikonal equation \eqref{eik}. There, we consider the region $\Omega=\Omega_{\delta,R,\eps}$ defined by \eqref{secreview:defnOmega} and assign the boundary condition $q|_{\partial\Omega}=r-t$. Meanwhile, we will obtain different optical functions by choosing different regions $\Omega$ and by assigning different boundary conditions. For example, for $\overline{\delta},\kappa\in(0,1)$, we can set
\fm{\overline{\Omega}_{\overline{\delta},R,\eps,\kappa}=\{(t,x)\in\R^{1+3}:\ t>e^{\overline{\delta}/\eps},\ |x|-e^{\overline{\delta}/\eps}-2R>\kappa(t-e^{\overline{\delta}/\eps})\}}
and solve the eikonal equation \eqref{eik} with $\Omega$ replaced by $\overline{\Omega}$. Note that $\Omega_{\delta,R,\eps}=\overline{\Omega}_{\delta,R,\eps,1/2}$. All the computations and proofs in Sections \ref{sec:prevac:opticalfunction} --  \ref{secreview:furthercomputations} will work, so we will again obtain functions such as $\wh{\overline{A}}$ (from \eqref{prevac:defhata}) and $\overline{A}_I$ (from Proposition \ref{prop:trzwtu}), etc. It is natural to ask how these new functions are related to the old ones. We seek to obtain a gauge independence result; see \cite[Section 6]{MR4772266}.

For convenience, we  make the following definition.

\defn{\rm Fix $\delta,\kappa\in(0,1)$. Let $q=q(t,x;\eps)$ be a $C^2$ optical function defined for all $\eps\ll_{\kappa,\delta}1$ and for all $(t,x)\in\R^{1+3}$ with $t>e^{\delta/\eps}$ and $|x|>\kappa t$. That is, we have
\fm{g^{\alpha\beta}(u)\partial_\alpha q\partial_\beta q=0,\qquad \forall t>e^{\delta/\eps}, |x|>\kappa t.}
Also suppose that $q|_{r-t\geq R}=r-t$. We say that the optical function $q$ is \emph{$(\delta,\kappa)$-admissible} if it satisfies the following assumptions:
\begin{enumerate}[a)]
\item For all $t>e^{\delta/\eps}$ and $|x|>\kappa t$, we have
\fm{\sum_{|I|\leq 1}|Z^I(q_t-q_r,q_r^{-1},(q_t-q_r)^{-1})|&\lesssim t^{C\eps};\\
\sum_{|I|\leq 1}\sum_{i=1,2,3}|Z^I(q_i+\omega_iq_t)|&\lesssim t^{-1+C\eps};}
\item The map $(s,q,\omega)=(\eps\ln t-\delta,q(t,x;\eps),x/|x|)$ induces a $C^1$ diffeomorphism from $\{t>e^{\delta/\eps},\ |x|>\kappa t\}$ to a subset of $[0,\infty)\times\R\times\mathbb{S}^2$, so every $C^1$ function of $(t,x)$ induces a $C^1$ function of $(s,q,\omega)$;
\item Define $(\mu,U)(t,x)=(q_t-q_r,\eps^{-1}ru)(t,x)$ and consider the induced functions  $(\mu,U)(s,q,\omega;\eps)$. The limits in \eqref{prevac:limita} exist, i.e.\ 
\fm{\left\{\begin{array}{l}
\displaystyle A(q,\omega;\eps):=-\frac{1}{2}\lim_{s\to\infty}(\mu U_q)(s,q,\omega;\eps),\\
\displaystyle A_1(q,\omega;\eps):=\lim_{s\to\infty}\exp(\frac{1}{2}G(\omega)A(q,\omega;\eps)s)\mu(s,q,\omega;\eps),\\
\displaystyle A_2(q,\omega;\eps):=\lim_{s\to\infty}\exp(-\frac{1}{2}G(\omega)A(q,\omega;\eps)s)U_q(s,q,\omega;\eps).
\end{array}\right.}
\end{enumerate}
}
\rmk{\rm For each $\kappa\in(1/2,1)$ and $\overline{\delta}\in[\delta,1)$, the optical function $q$ constructed in Section \ref{sec:prevac:opticalfunction} is $(\overline{\delta},\kappa)$-admissible. We cannot take $\kappa=1/2$ because $\Omega_{\delta,R,\eps}\subsetneq\{t>e^{\delta,\eps},|x|>t/2\}$.}
\rm

We now state the gauge independence result. It is indeed stronger than the one in \cite{MR4772266}.

\prop{\label{secreview:gauge:prop}Fix $\delta,\kappa,\overline{\delta},\overline{\kappa}\in(0,1)$. Let $q,\overline{q}$ be a $(\delta,\kappa)$-admissible optical function and a $(\overline{\delta},\overline{\kappa})$-admissible optical function, respectively. Define $\delta_0=\max\{\delta,\overline{\delta}\}$ and $\kappa_0=\max\{\kappa,\overline{\kappa}\}$. Then, whenever $\eps\ll_{\delta,\kappa,\overline{\delta},\overline{\kappa}}1$, we have 
\begin{enumerate}[\rm i)]
\item  There exists $Q_\infty=Q_\infty(q,\omega;\eps)\in C^1(\R\times\mathbb{S}^2)$ such that 
\fm{Q_\infty(q,\omega;\eps)=\lim_{s\to\infty}\overline{q}(s,q,\omega;\eps)}
where $\overline{q}(s,q,\omega;\eps)$ is the function induced by $\overline{q}(t,x;\eps)$ via the coordinate changes $(s,q,\omega)=(\eps\ln t-\delta,q(t,x;\eps),x/|x|)$. Besides,  we have
\fm{\partial_qQ_\infty=\lim_{s\to\infty}\partial_q\overline{q},\qquad \partial_{\omega_i}Q_\infty=\lim_{s\to\infty}\partial_{\omega_i}\overline{q}.}
All the convergences are here uniform in $(q,\omega)\in\R\times\mathbb{S}^2$.
\item
One can exchange the roles of $q$ and $\overline{q}$ to obtain $\overline{Q}_\infty(\overline{q},\omega;\eps)\in C^1(\R\times\mathbb{S}^2)$ such that
\fm{\overline{Q}_\infty(\overline{q},\omega;\eps)=\lim_{\overline{s}\to\infty}q(\overline{s},\overline{q},\omega;\eps),\\
\partial_q\overline{Q}_\infty=\lim_{\overline{s}\to\infty}\partial_{\overline{q}}q,\quad \partial_{\omega_i}\overline{Q}_\infty=\lim_{\overline{s}\to\infty}\partial_{\omega_i}q.}
Moreover,  we have \fm{Q_\infty(\overline{Q}_\infty(\overline{q},\omega;\eps),\omega;\eps)&=\overline{q},\\\overline{Q}_\infty(Q_\infty(q,\omega;\eps),\omega;\eps)&=q.}

\item Let $(A,A_1,A_2)(q,\omega;\eps)$ and $(\overline{A},\overline{A}_1,\overline{A}_2)(\overline{q},\omega;\eps)$ be the limits defined by \eqref{prevac:limita}, respectively. Then we have
\fm{\overline{A}(Q_\infty(q,\omega;\eps),\omega;\eps)&=A(q,\omega;\eps),\\
\overline{A}_1(Q_\infty(q,\omega;\eps),\omega;\eps)&=(A_1\cdot \partial_qQ_\infty)(q,\omega;\eps)\cdot\exp(\frac{1}{2}G(\omega)A(q,\omega;\eps)\cdot (\delta-\overline{\delta})).}
\end{enumerate}
}
\begin{proof}
The existence of $Q_\infty$ has been proved in \cite[Proposition 6.1]{MR4772266}. We recall that the key step in the proof is to write 
\fm{\partial_s\overline{q}&=\eps^{-1}t\overline{q}_r(\frac{\overline{\nu}}{\overline{q}_r}-\frac{\nu}{q_r}).}Here $\nu=q_t+q_r$ and $\overline{\nu}=\overline{q}_t+\overline{q}_r$. Using the eikonal equation, we write $\nu/q_r$ as the sum of a quantity not involving $q$ and a remainder term. We thus obtain $\partial_s\overline{q}=O(t^{-1+C\eps})$ and show that the limit $Q_\infty$ exists. Following a  similar proof, we can show $\partial_qQ_\infty=\lim_{s\to\infty}\partial_q\overline{q}$ and $ \partial_{\omega_i}Q_\infty=\lim_{s\to\infty}\partial_{\omega_i}\overline{q}$. The details are skipped. This finishes the proof of part i).

In part ii), we recall from the definition that both  $(s,q,\omega)=(\eps\ln t-\delta,q(t,x;\eps),x/|x|)$ and $(\overline{s},\overline{q},\omega)=(\eps\ln t-\overline{\delta},\overline{q}(t,x;\eps),x/|x|)$ are $C^1$ diffeomorphisms. 
These diffeomorphisms are defined whenever $t>e^{\delta_0/\eps}$ and $|x|>\kappa_0t$. We thus have two maps that are the inverses of each other:
\fm{(s,q,\omega)\mapsto (s+\delta-\overline{\delta},\overline{q}(s,q,\omega;\eps),\omega);\\(\overline{s},\overline{q},\omega)\mapsto (\overline{s}-\delta+\overline{\delta},q(\overline{s},\overline{q},\omega;\eps),\omega).}
It follows that for each $(s,q,\omega)$, we have
\fm{q&=q(s,\overline{q}(s+\delta-\overline{\delta},q,\omega;\eps),\omega;\eps).}
By the continuity and the uniform convergence, we obtain
\fm{q=\overline{Q}_\infty(Q_\infty(q,\omega;\eps),\omega;\eps).}
The other limit can be proved similarly.

In part iii),  we have proved $\overline{A}(Q_\infty(q,\omega;\eps),\omega;\eps)=A(q,\omega;\eps)$ in \cite[Proposition 6.1]{MR4772266}. To show the other limit, we recall that 
\fm{\overline{A}_1(\overline{q},\omega;\eps)=\lim_{\overline{s}\to\infty}\exp(\frac{1}{2}G(\omega)\overline{A}(\overline{q},\omega;\eps)\overline{s})\overline{\mu}(\overline{s},\overline{q},\omega;\eps).}
Here $\overline{s}=\eps\ln t-\overline{\delta}=s+\delta-\overline{\delta}$. Set $\overline{q}=Q_\infty(q,\omega;\eps)$, and we obtain 
\fm{\overline{A}_1(Q_\infty(q,\omega;\eps),\omega;\eps)=\lim_{s\to\infty}\exp(\frac{1}{2}G(\omega)A(q,\omega;\eps)\cdot (s+\delta-\overline{\delta}))\overline{\mu}(s+\delta-\overline{\delta},Q_\infty(q,\omega;\eps),\omega;\eps).}
Using the $C^1$ diffeomorphisms mentioned in the proof of part ii), we identify the following three points:
\eq{\label{secreview:gauge:threepoints}(s+\delta-\overline{\delta},Q_\infty(q,\omega;\eps),\omega)\qquad\text{ in }(\overline{s},\overline{q},\omega)\text{-coordinates},\\
(s,q(s+\delta-\overline{\delta},Q_\infty(q,\omega;\eps),\omega),\omega)\qquad\text{ in }(s,q,\omega)\text{-coordinates},\\
(e^{(s+\delta)/\eps},x(s,q,\omega;\eps))\qquad\text{ in }(t,x)\text{-coordinates}.}
Here $x(s,q,\omega;\eps)$ is chosen so that $x(s,q,\omega;\eps)=|x(s,q,\omega;\eps)|\omega$ and $\overline{q}(e^{(s+\delta)/\eps},x(s,q,\omega;\eps))=Q_\infty(q,\omega)$. Meanwhile, we have
\fm{\overline{\mu}=\overline{\nu}-2\overline{q}_r=\overline{\nu}-2q_r\partial_q\overline{q}=(\mu-\nu)\partial_q\overline{q}+O(t^{-1+C\eps}).}
At the point specified by \eqref{secreview:gauge:threepoints}, we have (using $(s,q,\omega)$-coordinates)
\fm{(\partial_q\overline{q})(s,\underline{q(s+\delta-\overline{\delta},Q_\infty(q,\omega;\eps),\omega;\eps)},\omega;\eps)\to (\partial_qQ_\infty)(q,\omega;\eps),\qquad s\to\infty.}
This limit holds because the underlined part converges to $\overline{Q}_\infty(Q_\infty(q,\omega;\eps),\omega;\eps)=q$.  As a result, we can skip the term $-\nu\partial_q\overline{q}+O(t^{-1+C\eps})$ as it will not affect the limit. In summary, we have
\fm{&\lim_{s\to\infty}\exp(\frac{1}{2}G(\omega)A(q,\omega;\eps)\cdot (s+\delta-\overline{\delta}))\overline{\mu}(s+\delta-\overline{\delta},Q_\infty(q,\omega;\eps),\omega;\eps)\\
&=\lim_{s\to\infty}\exp(\frac{1}{2}G(\omega)A(q,\omega;\eps)\cdot (s+\delta-\overline{\delta}))(\mu\cdot \partial_q\overline{q})(s,q(s+\delta-\overline{\delta},Q_\infty(q,\omega;\eps),\omega),\omega;\eps)\\
&=(A_1\cdot \partial_qQ_\infty)(q,\omega;\eps)\cdot \exp(\frac{1}{2}G(\omega)A(q,\omega;\eps)\cdot (\delta-\overline{\delta})).}
\end{proof}
\rm

We now discuss two important corollaries of this gauge independence result. First, suppose that $\delta<\overline{\delta}$, $\kappa<\overline{\kappa}$, and that $\overline{q}$ is simply a restriction of $q$. In this case, we have $Q_\infty=\overline{Q}_\infty=q$ everywhere, so part iii) of Proposition \ref{secreview:gauge:prop} suggests that
\fm{A\equiv\overline{A},\qquad \overline{A}_1\equiv A_1\exp(\frac{\delta-\overline{\delta}}{2}G(\omega)A).}
If we define $(\wt{\mu},\wt{U})$ and $(\wt{\overline{\mu}},\wt{\overline{U}})$ by \eqref{prevac:defnwtmuU} correspondingly, it is easy to check that 
\fm{(\wt{\overline{\mu}},\wt{\overline{U}})(s+\delta-\overline{\delta},q,\omega;\eps)=(\wt{\mu},\wt{U})(s,q,\omega;\eps).}
Such a time translation follows from the fact that the time $t$ corresponds to $s=\eps\ln t-\delta$ in $(s,q,\omega)$-coordinates and $\overline{s}=\eps\ln t-\overline{\delta}$ in $(\overline{s},\overline{q},\omega)$-coordinates.

The second corollary states that, if $\delta=\overline{\delta}$, then the outcome $\wh{A}$ defined by \eqref{prevac:defhata} should be the same.

\cor{\label{secreview:gauge:prop:cor}Fix $\delta,\kappa,\overline{\kappa}\in(0,1)$. Let $q,\overline{q}$ be a $(\delta,\kappa)$-admissible optical function and a $(\delta,\overline{\kappa})$-admissible optical function, respectively. Let $\wh{A},\wh{\overline{A}}$ be the functions defined by \eqref{prevac:defhata}, respectively. Then we have $\wh{A}\equiv \wh{\overline{A}}$.
}
\begin{proof}
Fix $(\rho,\omega)\in\R\times\mathbb{S}^2$. We seek to prove $\wh{A}(\rho,\omega;\eps)=\wh{\overline{A}}(\rho,\omega;\eps)$, or equivalently,
\fm{A(\wh{F}(\rho,\omega;\eps),\omega;\eps)= \overline{A}(\wh{\overline{F}}(\rho,\omega;\eps),\omega;\eps).}
The definitions of $\wh{F},\wh{\overline{F}}$ are given above \eqref{prevac:defhata}.
Since $\overline{A}(Q_\infty(q,\omega;\eps),\omega;\eps)=A(q,\omega;\eps)$, it suffices to prove
\fm{\wh{\overline{F}}(\rho,\omega;\eps)&=Q_\infty(\wh{F}(\rho,\omega;\eps),\omega;\eps).}

Set $q=\wh{F}(\rho,\omega;\eps)$ and $\overline{q}=\wh{\overline{F}}(\rho,\omega;\eps)$. By the definitions, we have
\fm{\rho=F(q,\omega;\eps)=\overline{F}(\overline{q},\omega;\eps).}
We seek to prove $\overline{q}=Q_\infty(q,\omega;\eps)$. Note that
\fm{F(q,\omega;\eps)&=2R-\int_{2R}^{q}\frac{2}{A_1(p,\omega;\eps)}\ dp=2R-\int_{2R}^{q}\frac{2\partial_qQ_\infty(p,\omega;\eps)}{\overline{A}_1(Q_\infty(p,\omega;\eps),\omega;\eps)}\ dp\\
&=2R-\int_{2R}^{Q_\infty(q,\omega;\eps)}\frac{2}{\overline{A}_1(p,\omega;\eps)}\ dp=\overline{F}(Q_\infty(q,\omega;\eps),\omega;\eps).}
In the second identity, we use part iii) of Proposition \ref{secreview:gauge:prop}. Since $Q_\infty(\cdot;\omega;\eps)$ has an inverse, its derivative with respect to $q$ cannot vanish. Moreover, since $q,\overline{q}=r-t$ for $r-t\geq R$, we have $Q_\infty(2R,\omega;\eps)=2R$. It follows that
\fm{\rho=F(q,\omega;\eps)=\overline{F}(Q_\infty(q,\omega;\eps),\omega;\eps)=\overline{F}(\overline{q},\omega;\eps).}
Since $\partial_q\overline{F}=2/A_1\neq 0$, we must have $\overline{q}=Q_\infty(q,\omega;\eps)$. This finishes our proof.
\end{proof}
\rm

\section{A general result for inhomogeneous wave equations}\label{sec:generalinhom}
We now consider the following  general inhomogeneous wave equation: 
\eq{\label{inhomwe} \Box \phi=F(t,x)\qquad \text{in }\Oc.}
Here the region $\Oc\subset\R^{1+3}$ is defined by 
\eq{\Oc=\{(t,x)\in\R^{1+3}:\ t>0,\ |x|<t\}.}
The main proposition in this section is a representation formula for a solution $\phi$ in $\Oc$.

\prop{\label{prop:inhomwe} Suppose that $\phi\in C^2(\Oc)$ is a solution to \eqref{inhomwe}. Suppose that  for some constants $M,\gamma_1,\gamma_2>0$, we have \fm{|F(t,x)|\leq M\lra{t}^{-\gamma_1}\lra{|x|-t}^{-\gamma_2},\qquad \forall (t,x)\in\Oc.}Then, by setting $\Phi(t,x)=(-\partial_t+\partial_r)(|x|\phi)$, for all $(t,x)\in\Oc$ and $T> 2t$, we have 
\eq{\label{prop:inhomwe:clde}\phi(t,x)&=\frac{1}{4\pi}\int_{\mathbb{S}^2}\Phi(T,x-(T-t)\theta)\ dS_\theta\\
&\quad+O(|x|\norm{\partial\phi(T)}_{L^\infty}+M\lra{t-|x|}^{-\gamma_2}\cdot(\lra{t}^{2-\gamma_1}+\int_{t^2}^{(T-t)^2} (1+\rho)^{-\gamma_1/2}\ d\rho)).}
The implicit constant in $O(\dots)$ is universal and does not depend on $T,t,x,M,\gamma_1,\gamma_2$.
}
\rmk{\rm Note that the integral can be estimated explicitly: 
\fm{\int_{t^2}^{\wt{T}^2}(1+\rho)^{-\gamma_1/2}\ d\rho\leq \left\{
\begin{array}{ll}
  \frac{\lra{\wt{T}}^{2-\gamma_1}}{1-\gamma_1/2},  & \gamma_1<2; \\[1em]
    2\ln(\frac{\lra{\wt{T}}}{\lra{t}}),  & \gamma_1=2;\\[1em]
    \frac{\lra{t}^{2-\gamma_1}}{\gamma_1/2-1},  & \gamma_1>2.
\end{array}
\right.
}}\rm

\bigskip

It suffices to prove this proposition for $M=1$. For general $M>0$, we simply replace $(\phi,F)$ with $(\phi/M,F/M)$ in \eqref{inhomwe}. In the proof, we view $\phi$ as a solution to \eqref{inhomwe} with initial data $(\phi_0,\phi_1)=(\phi,\phi_t)(T)$ at $t=T$. This allows us to write $\phi=\phi_{\rm lin}+\phi_{\rm inh}$. Here we have $\phi_{\rm lin}$ solves $\Box \phi_{\rm lin}=0$ with data $(\phi_0,\phi_1)$ at $t=T$, and $\phi_{\rm inh}$ solves \eqref{inhomwe} with zero data at $t=T$. See Lemma \ref{lem:inhomsol}. Note that $\phi_{\rm lin}$ can be expressed explicitly in terms of $(\phi_0,\phi_1)$ and that $\phi_{\rm inh}$ can be expressed explicitly in terms of $F$. In Lemmas \ref{lemestphilin} and \ref{lemestphiinh}, we estimate $\phi_{\rm lin}$ and $\phi_{\rm inh}$ by applying these explicit formulas.

\subsection{A representation formula}
\lem{\label{lem:inhomsol} Let $\phi\in C^2(\Oc)$ be a solution to \eqref{inhomwe}. Then, for all $(t,x)\in\Oc$ and $T> t$, we have
\eq{\label{eq:inhomsol}4\pi \phi(t,x)&=\int_{\mathbb{S}^2} \phi(T,x-(T-t)\theta)-(T-t)(\phi_t+\theta\cdot\nabla_x\phi)(T,x-(T-t)\theta)\ dS_\theta\\
&\quad -\int_{y\in\R^3:\ |y|<T-t}F(t+|y|,x+y) |y|^{-1}\ dy.}}
\begin{proof}
Fix $T>0$ and set $w(t,x)=\phi(T-t,-x)$. Note that $w\in C^2(\{0<t<T,\ |x|<T-t\})$ and that $w$  solves $(\Box w)(t,x)=F(T-t,-x)$ in its domain. Using formulas from, e.g., \cite{MR2455195,MR1466700}, we can explicitly express $w$ in terms of $(w,w_t)|_{t=0}$ and $F$. Whenever $t\in(0,T)$ and $|x|<T-t$, we have
\fm{4\pi w(t,x)&=\int_{\mathbb{S}^2}\partial_t(tw(0,x+t\theta))+tw_t(0,x+t\theta)\ dS_\theta-\int_{y\in\R^3:\ |y|<t}(\Box w)(t-|y|,x-y)\  \frac{dy}{|y|} .}
Since $(w,w_t)|_{t=0}=(\phi(T,-x),-\phi_t(T,-x))$, by the chain rule, we have \fm{&\partial_t(tw(0,x+t\theta))+tw_t(0,x+t\theta)=w(0,x+t\theta)+t(w_t+\theta\cdot \nabla_xw)(0,x+t\theta)\\
&=\phi(T,-x-t\theta)-t(\phi_t+\theta\cdot \nabla_x \phi)(T,-x-t\theta).}
And since $(\Box w)(t-|y|,x-y)=F(T-t+|y|,y-x)$, we have 
\fm{4\pi w(t,x)&=\int_{\mathbb{S}^2}\phi(T,-x-t\theta)-t(\phi_t+\theta\cdot \nabla_x \phi)(T,-x-t\theta)\ dS_\theta\\
&\quad-\int_{y\in\R^3:\ |y|<t}|y|^{-1}F(T-t+|y|,y-x)\  dy. }
Finally, by replacing $(t,x)$ with $(T-t,-x)$, we obtain the formula for $\phi$. This is because $\phi(t,x)=w(T-t,-x)$.
\end{proof}\rm

\bigskip

For each $(t,x)\in\Oc$ and $T>t$, we rewrite \eqref{eq:inhomsol} as $\phi=\phi_{\rm lin}+\phi_{\rm inh}$ where 
\eq{\label{defnphilin}\phi_{\rm lin}(t,x;T)=\frac{1}{4\pi}\int_{\mathbb{S}^2} \phi(T,x-(T-t)\theta)-(T-t)(\phi_t+\theta\cdot\nabla_x\phi)(T,x-(T-t)\theta)\ dS_\theta,}
\eq{\label{defnphiinh}\phi_{\rm inh}(t,x;T)=-\frac{1}{4\pi}\int_{y\in\R^3:\ |y|<T-t}F(t+|y|,x+y) |y|^{-1}\ dy.}
When $T$ is fixed, we note that $\phi_{\rm lin}$ is the solution to the linear wave equation $\Box \phi_{\rm lin}=0$ in $\Oc\cap\{0<t<T\}$ with data $(\phi_{\rm lin},\partial_t\phi_{\rm lin})|_{t=T}=(\phi(T),\phi_t(T))$ and that $\phi_{\rm inh}$ is the solution to the inhomogeneous equation $\Box \phi_{\rm inh}=F$ in $\Oc\cap\{0<t<T\}$ with zero data at $t=T$.

In the next two subsections, we will estimate $\phi_{\rm lin}$ and $\phi_{\rm inh}$ respectively. Proposition \ref{prop:inhomwe} will then follow directly from Lemmas \ref{lemestphilin} and \ref{lemestphiinh} below. 

\subsection{Estimates for $\phi_{\rm lin}$}
Let us start with $\phi_{\rm lin}$ defined by \eqref{defnphilin}. 
\lem{\label{lemestphilin}Set $\Phi(t,x)=(-\partial_t+\partial_r)(|x|\phi)$. Then, for all $(t,x)\in\Oc$ with $t<T$,  we have
\fm{\phi_{\rm lin}(t,x;T)&=\frac{1}{4\pi}\int_{\mathbb{S}^2} \Phi(T,x-(T-t)\theta)\ dS_\theta+O(|x|\norm{\partial\phi(T)}_{L^\infty(\{x\in\R^3:\ |x|<T\})}).}
The implicit constant in $O(\dots)$ is universal and does not depend on $T,t,x$.
}
\begin{proof}
For simplicity, we write $y=y(t,x,\theta,T)=x-(T-t)\theta$. Then,  we have
\fm{\Phi(T,y)&=\phi(T,y)-|y|(\phi_t-\phi_r)(T,y)=\phi(T,y)-|y|\phi_t(T,y)+\sum_{j=1}^3y_j\phi_{y_j}(T,y),}
while the integrand on the right side of \eqref{defnphilin} is
\fm{\phi(T,y)-(T-t)(\phi_t+\theta\cdot\nabla_y\phi)(T,y)=\phi(T,y)-(T-t)\phi_t(T,y)-\sum_{j=1}^3\theta_j(T-t)\phi_{y_j}(T,y).}
The difference between these two formulas is $(T-t-|y|)\phi_t(T,y)+\sum_{j=1}^3(y_j+(T-t)\theta_j)\phi_{y_j}(T,y)$. By the definition of $y$, we have $y_j+(T-t)\theta_j=x_j$ and
\fm{T-t-|x|&\leq |y|=|x-(T-t)\theta|\leq |x|+T-t\Longrightarrow |T-t-|y||\leq |x|.}
Here we use $\theta\in\mathbb{S}^2$ and $T\geq t$. This finishes the proof.
\end{proof}
\rm

\subsection{Estimates for $\phi_{\rm inh}$}

We now estimate $\phi_{\rm inh}$ defined by \eqref{defnphiinh}.

\lem{\label{lemestphiinh}Suppose that $|F(t,x)|\leq \lra{t}^{-\gamma_1}\lra{|x|-t}^{-\gamma_2}$ in $\Oc$ for some constants $\gamma_1,\gamma_2>0$. Then, for all $t,x,\wt{T}$ with $(t,x)\in\Oc$ and $\wt{T}>t$ (later we will set $\wt{T}=T-t$), we have
\fm{\abs{\int_{y\in\R^3:\ |y|<\wt{T}}F(t+|y|,x+y)|y|^{-1}\ dy} \leq 2\pi \lra{t-|x|}^{-\gamma_2}\cdot(\lra{t}^{2-\gamma_1}+\int_{t^2}^{\wt{T}^2} (1+\rho)^{-\gamma_1/2}\ d\rho).}
\begin{proof}
Using the bound for $F$, we have
\fm{&\abs{\int_{y\in\R^3:\ |y|<\wt{T}}F(t+|y|,x+y)|y|^{-1}\ dy}\\&\leq \int_{y\in\R^3:\ |y|<\wt{T}}\lra{t+|y|}^{-\gamma_1}\lra{|x+y|-t-|y|}^{-\gamma_2}|y|^{-1}\ dy\\
&\leq \lra{t-|x|}^{-\gamma_2}\int_{y\in\R^3:\ |y|<\wt{T}}\lra{t+|y|}^{-\gamma_1}|y|^{-1}\ dy\\
&\leq 4\pi \lra{t-|x|}^{-\gamma_2}\int_{0}^{\wt{T}}\lra{t+\rho}^{-\gamma_1}\rho\ d\rho .}
To get the second estimate, we recall  $|x|<t$ and $|x+y|-|y|\leq |x|$. Since $s\mapsto \lra{s}$ is increasing for $s\geq 0$, we have
$\lra{t-|x+y|+|y|}\geq \lra{t-|x|}$. In the last step, we write the integral in polar coordinates.

Next, we notice that 
\fm{\int_{0}^{\wt{T}}\lra{t+\rho}^{-\gamma_1}\rho\ d\rho&\leq \int_{0}^{t}\lra{t}^{-\gamma_1}\rho\ d\rho+\int_{t}^{\wt{T}}(1+\rho^2)^{-\gamma_1/2}\rho\ d\rho\\
&\leq \frac{1}{2}\lra{t}^{2-\gamma_1}+\frac{1}{2}\int_{t^2}^{\wt{T}^2}(1+\rho)^{-\gamma_1/2}\ d\rho.}
This finishes the proof.
\end{proof}
\rm
\bigskip

It is now obvious that Proposition \ref{prop:inhomwe}  follows directly from Lemmas \ref{lemestphilin} and \ref{lemestphiinh}.  We also note that we assume $T>2t$ in Proposition \ref{prop:inhomwe} because $T-t=\wt{T}>t$ in Lemma \ref{lemestphiinh}.

\section{Proof of Theorem \ref{mthm}}\label{sec:pfmthm}

We now prove Theorem \ref{mthm}. Let $u$ be the global solution to \eqref{qwe} in the statement of the theorem. In Section \ref{sec:pfmthm:setup}, we apply Proposition \ref{prop:inhomwe} to $\phi=Z^Iu$ for a fixed multiindex $I$. Because of \eqref{sec:pfmthm:keyest2}, it remains to compute the integral \fm{\int_{\mathbb{S}^2}\Phi_I(T,x-(T-t)\theta;\eps)\ dS_\theta}
where  $\Phi_I=(-\partial_t+\partial_r)(rZ^Iu)$. This computation is then done in Section \ref{sec:pfmthm:5.2}. By setting\footnote{The proof will still work if we set $T=t^\lambda$ with a fixed parameter $\lambda>2$. However, we  need the estimate $T^{C\eps}\lesssim t^{C\eps}$ in our proof, so we cannot set this $T$ arbitrarily. For example, we cannot set $T=e^t$. } $T=t^3$ and $y=y(t,x,T,\theta)=x-(T-t)\theta$, we have $(T,y)\in\Omega$ and  $||y|-T|\lesssim T^{1/3}$ whenever $t>(2e^{\delta/\eps})^{1/3}$, $|x|<t$, and $\theta\in\mathbb{S}^2$. Here we recall that $\Omega=\Omega_{\delta,R,\eps}$ was defined by \eqref{secreview:defnOmega}. As a result, we can apply \eqref{prevac:uapp} to show \fm{\Phi_I(T,y;\eps)\approx -r(\partial_t-\partial_r)Z^I\wt{u}(T,y;\eps).} Moreover, by Proposition \ref{prop:trzwtu}, we also have 
\fm{-r(\partial_t-\partial_r)Z^I\wt{u}(T,y;\eps)=-r(\partial_t-\partial_r)Z^I\wh{u}(T,y;\eps)\approx -\eps A_I(|y|-T,y/|y|;\eps).}
Finally, by Lemma \ref{pfmthm:lem:basic}, we have $|y|-T\approx -t-x\cdot\theta$ and $y/|y|\approx-\theta$. We thus have
\fm{\int_{\mathbb{S}^2}\Phi_I(T,x-(T-t)\theta;\eps)\ dS_\theta\approx -\eps\int_{\mathbb{S}^2}A_I(-t-x\cdot\theta,-\theta;\eps)\ dS_\theta.}
Thus, we obtain the integral on the right hand side of \eqref{mthmcl}. We still need to estimate the error terms in these computations. The details will be given in the proof below.

\subsection{Setup}\label{sec:pfmthm:setup}
Fix $(u_0,u_1)\in C_c^\infty(\R^3)$.  Let $u$ be the global $C^\infty$ solution to \eqref{qwe} with data $(\eps u_0,\eps u_1)$ for $\eps\ll 1$. 
Now fix a multiindex $I$. Set $N=|I|+2$ and choose $\eps\ll_N1$ so that \eqref{secreview:uptb} holds for $N$.

In order to apply Proposition \ref{prop:inhomwe} to $Z^Iu$, we need to  estimate $\Box Z^I u$. In fact, since $g^{\alpha\beta}(u)=m^{\alpha\beta}+O(|u|)$, it follows from  \eqref{secreview:uptb} and Lemma \ref{lem:prel:taylor} (with $l_0=1$) that
\fm{|Z^I\Box u|&=|Z^I\kh{(g^{\alpha\beta}(u)-m^{\alpha\beta})\partial_\alpha\partial_\beta u}|\lesssim \sum_{|J_1|+|J_2|=|I|}|Z^{J_1}(g^{\alpha\beta}(u)-m^{\alpha\beta})||Z^{J_2}\partial^2u|\\
&\lesssim \sum_{|J|\leq |I|}|Z^Ju| \cdot\lra{r-t}^{-2}\sum_{|J|\leq |I|+2}|Z^Ju|\lesssim \lra{r-t}^{-2}(\eps \lra{t}^{-1+C\eps})^2\lesssim \eps^2 \lra{t}^{-2+C\eps}\lra{r-t}^{-2}.}

Now we apply Proposition \ref{prop:inhomwe} to $Z^Iu$ (with $M=C\eps^2$, $\gamma_1=2-C\eps$ and $\gamma_2=2$). By setting $\Phi_I(t,x;\eps)=(-\partial_t+\partial_r)(rZ^Iu)(t,x;\eps)$, for all $(t,x)$ with $t>0$ and $|x|<t$, and all $T>2t$, we have
\eq{\label{sec:pfmthm:keyest1}&\abs{(Z^Iu)(t,x;\eps)-\frac{1}{4\pi}\int_{\mathbb{S}^2}\Phi_I(T,x-(T-t)\theta;\eps)\ dS_\theta}\\
&\lesssim|x|\norm{\partial Z^Iu(T)}_{L^\infty}+\eps^2\lra{t-|x|}^{-2}\cdot(\lra{t}^{C\eps}+\int_{t^2}^{(T-t)^2} (1+\rho)^{-1+C\eps}\ d\rho)\\
&\lesssim \eps tT^{-1+C\eps}+\eps^2\lra{t-|x|}^{-2}\lra{t}^{C\eps}+\eps\lra{t-|x|}^{-2}T^{C\eps}.}
If $t>2$ and if we set $T=t^3>2t$, then $T^{C\eps}\lesssim t^{C\eps}$ and $tT^{-1+C\eps}\lesssim t^{-2+C\eps}\lesssim \lra{t-|x|}^{-2}t^{C\eps}$ (since $0<t-|x|<t$). That is, we have
\eq{\label{sec:pfmthm:keyest2}\abs{(Z^Iu)(t,x;\eps)-\frac{1}{4\pi}\int_{\mathbb{S}^2}\Phi_I(T,x-(T-t)\theta;\eps)\ dS_\theta}\lesssim \eps  \lra{t-|x|}^{-2}t^{C\eps}.}
It remains to compute the integral of $\Phi_I$.

\subsection{Computing the integral of $\Phi_I$}\label{sec:pfmthm:5.2}
We now apply the results in Section \ref{secreview}. Since most estimates there hold in $\Omega$ defined by \eqref{secreview:defnOmega}, our first goal is to show that
\eq{\label{sec:pfmthm:est3}(T,x-(T-t)\theta)\in\Omega,\qquad \forall t> (2e^{\delta/\eps})^{1/3}, T=t^3, |x|<t, \theta\in\mathbb{S}^2.}
Here we recall that $\delta\in(0,1)$ is a small constant chosen at the beginning of Section \ref{sec:prevac:opticalfunction}. By choosing $\eps\ll_{\delta,R}1$, we have $T>2e^{\delta/\eps}$ and 
\fm{&|x-(T-t)\theta|-(\frac{1}{2}(T+e^{\delta/\eps})+2R)\geq (T-t)-|x|-\frac{1}{2}(T+e^{\delta/\eps})-2R\\
&\geq \frac{1}{2}T-2t-\frac{1}{2}e^{\delta/\eps}-2R\geq \frac{1}{2}T(\frac{1}{2}-4T^{-2/3}-2RT^{-1})>0.}
Thus, we have \eqref{sec:pfmthm:est3}.

Next, we set $y=y(t,x,T,\theta)=x-(T-t)\theta$. Since $T\gg t>|x|$, we expect $y\approx -T\theta$. This expectation is made rigorous in the next lemma.

\lem{\label{pfmthm:lem:basic}For all $t>(2e^{\delta/\eps})^{1/3}$, $T=t^3$, $|x|<t$, and $\theta\in\mathbb{S}^2$, we have $|y|\sim T$, $T-|y|=t+x\cdot\theta+O(t^{-1})\in[t-|x|-1,2t+1]$, and $y/|y|+\theta=O(t^{-2})$.}
\begin{proof}
We have
$||y|-(T-t)|\leq |x|$. It  follows that \fm{|T-|y||&\leq ||y|-(T-t)|+t\leq t+|x|\leq 2t,\\
|\frac{|y|}{T}-1|&\leq \frac{2t}{T}\leq 2T^{-2/3}\leq 2(e^{\delta/\eps})^{-2/3}.}
For $\eps\ll_\delta1$, we have $(e^{\delta/\eps})^{-2/3}=e^{-\frac{2\delta}{3\eps}}<1/4$, so $|y|\sim T$.

Next, we have
\fm{T^2-|y|^2&=T^2-(|x|^2+(T-t)^2-2(T-t)x\cdot\theta)=2tT-t^2-|x|^2+2(T-t)x\cdot\theta,\\
T-|y|&=\frac{2tT-t^2-|x|^2+2(T-t)x\cdot\theta}{T+|y|}\\
&=t+x\cdot\theta+\frac{t(T-|y|)-t^2-|x|^2+(T-|y|-2t)x\cdot\theta}{T+|y|}\\
&=t+x\cdot\theta+O(t^2/T)=t+x\cdot\theta+O(t^{-1}).}
In the second last estimate, we use $|x|<t$, $|T-|y||\leq 2t$, and $|y|\sim T$. It is also clear that $t-|x|-1\leq T-|y|\leq 2t+1$.

Finally, we have
\fm{\abs{\frac{y}{|y|}+\theta}&=\frac{|x-(T-t-|y|)\theta|}{|y|}\lesssim T^{-1}(|x|+t)\lesssim t/T=t^{-2}.}
\end{proof}\rm

We can now compute $\Phi_I(T,y;\eps)$. For simplicity, in the remaining computations of the current subsection, we omit the parameter $\eps$ in the functions $\Phi_I,u,\wt{u},A_I$, etc. That is, we will only write $\Phi_I(T,y)$.  Since $||y|-T|\lesssim t=T^{1/3}$, by applying \eqref{prevac:uapp} with $N=|I|+2$, for $\eps\ll_N1$ we have
\fm{\Phi_I(T,y)&=-[r(\partial_t-\partial_r)(Z^Iu)](T,y)+Z^Iu(T,y)\\
&=-[r(\partial_t-\partial_r)(Z^I\wt{u})](T,y)-[r(\partial_t-\partial_r)(Z^I(u-\wt{u}))](T,y)+O(\eps T^{-1+C\eps})\\
&=-[r(\partial_t-\partial_r)(Z^I\wt{u})](T,y)+O(\eps T^{-1+C\eps}).}
Here we use $|\partial Z^I(u-\wt{u})|\lesssim \lra{r-t}^{-1}|Z Z^I(u-\wt{u})|$ to remove the $\lra{r-t}$ factor in \eqref{prevac:uapp}.
By Proposition \ref{prop:trzwtu}, we have \fm{&-[r(\partial_t-\partial_r)(Z^I\wt{u})](T,y)\\
&=-\eps A_I(|y|-T,y/|y|)+O(\eps (1+\ln\lra{|y|-T})\lra{|y|-T}^{-2}T^{C\eps}+\eps T^{-1+C\eps})\\
&=-\eps A_I(|y|-T,y/|y|)+O(\lra{t-|x|}^{-2}t^{C\eps}+\eps t^{-3+C\eps}).}
In the last step, we use that $\eps\ln \lra{T-|y|}\lesssim \eps+\eps\ln T\lesssim \eps+\eps\ln t\lesssim t^{C\eps}$ and that $\lra{|y|-T}^{-2}\lesssim (1+t-|x|)^{-2}\lesssim \lra{t-|x|}^{-2}$ by Lemma \ref{pfmthm:lem:basic}. Moreover, for $\rho$ lies between $|y|-T$ and $-t-x\cdot\theta$, 
\fm{-\rho\geq \min\{T-|y|,t+x\cdot\theta\}\geq t+x\cdot\theta-Ct^{-1}\geq t-|x|-1>-1.} Thus, we have $\lra{\rho}\sim (2-\rho)\geq (t-|x|+1)\sim \lra{t-|x|}$. By the mean value theorem, we have
\fm{&|-\eps A_I(|y|-T,y/|y|)+\eps A_I(-t-x\cdot\theta,-\theta)|\\
&\lesssim \eps|A_I(|y|-T,y/|y|)- A_I(-t-x\cdot\theta,y/|y|)|+\eps|A_I(-t-x\cdot \theta,y/|y|)- A_I(-t-x\cdot\theta,-\theta)|\\
&\lesssim\eps||y|-T+t+x\cdot\theta|\cdot \sup_{\rho\text{lies between }|y|-T,\ -t-x\cdot\theta}|\partial_qA_I(\rho,y/|y|)|\\
&\quad+\eps|\frac{y}{|y|}+\theta|\cdot \sup_{\eta\in\mathbb{S}^2}|\partial_\omega A_I(-t-x\cdot\theta,\eta)|\\
&\lesssim\eps t^{-1}\sup_{\rho\text{ lies between }|y|-T,\ -t-x\cdot\theta}\lra{\rho}^{-2+C\eps}+\eps t^{-2}\lra{t+x\cdot\theta}^{-1+C\eps}\\
&\lesssim \eps t^{-1}\lra{t-|x|}^{-2+C\eps}+\eps t^{-2}\lra{t-|x|}^{-1+C\eps}\lesssim \eps t^{-1+C\eps}\lra{|x|-t}^{-2}.}
In summary, we have
\fm{\Phi_I(T,y;\eps)&=-\eps A_I(-t-x\cdot\theta,-\theta;\eps)+O(\lra{r-t}^{-2}t^{C\eps}).}
Plugging this back to \eqref{sec:pfmthm:keyest2}, we obtain
\eq{\label{pfmthm:final}Z^Iu(t,x;\eps)&=-\frac{\eps }{4\pi}\int_{\mathbb{S}^2} A_I(-t+x\cdot\theta,\theta;\eps)\ dS_\theta+O(\lra{t-|x|}^{-2}t^{C\eps}),\  \forall t>(2e^{\delta/\eps})^{1/3},\ |x|<t.}

\section{Corollaries of Theorem \ref{mthm}}\label{sec:corofmthm}

Combining Theorem \ref{mthm} with the asymptotic completeness result for \eqref{qwe}, we can prove several interesting corollaries.

In Section \ref{sec:prop:decaysphmeans}, we show a special decaying property for $\wh{A}$ and $\wh{U}$; see \eqref{prop:decaysphmeans:c2} in Proposition \ref{prop:decaysphmeans}.  In the proof, we compute the integral  $\int_{\mathbb{S}^2}u(t,(t-t^{2/3})\omega)\ dS_\omega$ in two ways. First, we can use the estimates from the asymptotic completeness result. By \eqref{prevac:uapp} and Proposition \ref{prop:trzwtu}, we have
\fm{\int_{\mathbb{S}^2}u(t,(t-t^{2/3})\omega)\ dS_\omega\approx \eps (t-t^{2/3})^{-1}\cdot\int_{\mathbb{S}^2}\wh{U}(\eps\ln t-\delta,-t^{2/3},\omega)\ dS_\omega,\qquad \forall t\gg_{\delta,\eps}1.}
Meanwhile, we can apply the main theorem for this paper. By \eqref{mthmcl}, we have
\fm{&\int_{\mathbb{S}^2}u(t,(t-t^{2/3})\omega)\ dS_\omega\approx \frac{\eps}{2\pi} \int_{\mathbb{S}^2}\int_{\mathbb{S}^2}\wh{A}(-t+(t-t^{2/3})\omega\cdot\theta,\theta;\eps)\ dS_\theta dS_\omega\\
&=\frac{\eps}{t-t^{2/3}} \int_{\mathbb{S}^2}\int_{-2t+t^{2/3}}^{-t^{2/3}}\wh{A}(\rho,\theta;\eps)\ d\rho dS_\theta\\
&\approx\frac{\eps}{t-t^{2/3}} \int_{\mathbb{S}^2}[\wh{U}(\eps\ln t-\delta,-t^{2/3},\omega;\eps)-\wh{U}(\eps\ln t-\delta,-2t+t^{2/3},\omega;\eps)]\  dS_\theta,\qquad \forall t\gg_{\delta,\eps}1.}
In the second step, we integrate with respect to $\omega$. In the last step, we notice that $\wh{U}_q\approx \wh{A}$. Comparing these two approximate identities, we notice that there is a cancellation. We thus have
\fm{\frac{\eps}{t-t^{2/3}} \int_{\mathbb{S}^2} \wh{U}(\eps\ln t-\delta,-2t+t^{2/3},\omega;\eps)\  dS_\theta\approx 0}
to some extent. For simplicity, we skip the discussion on how to control the error terms in these computations, and we refer to the proof of Proposition \ref{prop:decaysphmeans} below. 
Following the same idea, one can derive similar decay properties for the $U_I$'s defined in Proposition \ref{prop:trzwtu}; see~\eqref{prop:decaysphmeans:c1}.

In Sections \ref{sec:prop:vanishingA} and \ref{sec:pfltthm}, we prove Theorems \ref{sdltthm} and \ref{ltthm_weak}. Here Proposition \ref{prop:vanishingA} is the same as Theorem \ref{sdltthm}, while Proposition \ref{ltthm} is slightly stronger than Theorem \ref{ltthm_weak}.
Note that the ideas of the proofs of Propositions \ref{prop:vanishingA} and \ref{ltthm} have been explained in the introduction, so we will not repeat them here.

\subsection{Decaying properties of the $U_I$'s}\label{sec:prop:decaysphmeans}
We now discuss the decaying properties mentioned at the beginning of this section.
\prop{\label{prop:decaysphmeans}For each multiindex $I$, we let $U_I$ be the function of $(s,q,\omega)$ defined in the statement of Proposition \ref{prop:trzwtu}. Then, for all $t\geq 2e^{\delta/\eps}$, we have
\eq{\label{prop:decaysphmeans:c1}|\int_{\mathbb{S}^2}U_I(\eps\ln t-\delta,-2t,\omega;\eps)\ dS_\omega|\lesssim\eps^{-1}t^{-1/3+C\eps}.}
In particular, since $U_0=\wh{U}$ which was defined in \eqref{prevac:defhata}, we have
\eq{\label{prop:decaysphmeans:c2}|\int_{\mathbb{S}^2}\int_{-2t}^R\wh{A}(q,\omega;\eps)\exp(\frac{1}{2}G(\omega)\wh{A}(q,\omega;\eps)\cdot(\eps\ln t-\delta))\ dqdS_\omega|\lesssim\eps^{-1}t^{-1/3+C\eps}.}}
\begin{proof}
Fix a constant $\gamma\in(0,1)$. For each $(t,x)\in\Omega\cap\{-2t^\gamma<r-t<-t^\gamma/2\}$, we have two ways to approximate $Z^Iu$. By \eqref{prevac:uapp} and Proposition \ref{prop:trzwtu}, we have
\fm{Z^Iu&=Z^I\wh{u}+O_\gamma(\eps t^{-2+C\eps}\lra{r-t})\\
&=\eps r^{-1}U_I(\eps\ln t-\delta,r-t,\omega;\eps)+O_\gamma(t^{-1-\gamma+ C\eps}+\eps t^{-2+\gamma+C\eps}).}
Here we use $|r-t|\sim t^{\gamma}$.
Meanwhile, by Theorem \ref{mthm}, we have
\fm{Z^Iu&=-\frac{\eps}{4\pi}\int_{\mathbb{S}^2}A_I(-t+x\cdot\theta,\theta;\eps)\ dS_\theta+O(\lra{r-t}^{-2}t^{C\eps})\\
&=-\frac{\eps}{4\pi}\int_{\mathbb{S}^2}A_I(-t+x\cdot\theta,\theta;\eps)\ dS_\theta+O(t^{-2\gamma+C\eps}).}
By combining these two estimates, we have
\eq{\label{prop:decaysphmeans:est}&U_I(\eps\ln t-\delta,r-t,\omega;\eps)+\frac{r}{4\pi}\int_{\mathbb{S}^2}A_I(-t+x\cdot\theta,\theta;\eps)\ dS_\theta\\
&=O_\gamma(\eps^{-1}t^{1-2\gamma+ C\eps}+ t^{-1+\gamma+C\eps})}
whenever $(t,x)\in\Omega$ and $-2t^\gamma<r-t<-t^\gamma/2$.

To continue, we let $t>2e^{\delta/\eps}$ and $r=t-t^{2/3}$ in \eqref{prop:decaysphmeans:est}. Note that $\max\{1-2\gamma,-1+\gamma\}$ from the power of $t$ in \eqref{prop:decaysphmeans:est} is minimized at $\gamma=2/3$. Moreover, we have
\fm{r-(t+e^{\delta/\eps})/2-2R=(t-e^{\delta/\eps})/2-t^{2/3}-2R>t/4-t^{2/3}-2R>0,}
so $\{t>2e^{\delta/\eps},r=t-t^{2/3}\}\subset\Omega$. This gives us 
\fm{U_I(\eps\ln t-\delta,-t^{2/3},\omega;\eps)+\frac{r}{4\pi}\int_{\mathbb{S}^2}A_I(-t+r\omega\cdot\theta,\theta;\eps)\ dS_\theta=O(\eps^{-1}t^{-1/3+ C\eps}).}
Integrate both sides with respect to $\omega$ on $\mathbb{S}^2$, and we get
\fm{\int_{\mathbb{S}^2}U_I(\eps\ln t-\delta,-t^{2/3},\omega;\eps)\ dS_\omega&=-\frac{r}{4\pi}\int_{\mathbb{S}^2\times\mathbb{S}^2}A_I(-t+r\omega\cdot\theta,\theta;\eps)\ dS_\omega dS_\theta+O(\eps^{-1}t^{-1/3+ C\eps})\\
&=-\frac{1}{2}\int_{\mathbb{S}^2}\int_{-t-r}^{-t+r}A_I(\rho,\theta;\eps)\ d\rho dS_\theta+O(\eps^{-1}t^{-1/3+ C\eps}).}
Note that
\eq{\label{decaysphmeans:sphint}&\int_{\mathbb{S}^2}A_I(-t+r\omega\cdot\theta,\theta;\eps)\ dS_\omega =\int_{\mathbb{S}^2}A_I(-t+r\omega_3,\theta;\eps)\ dS_\omega\\
&=\int_0^{2\pi}\int_{0}^\pi A_I(-t+r\cos\alpha,\theta;\eps)\sin\alpha\ d\alpha d\beta=2\pi \int_{-1}^1A_I(-t+r\tau,\theta;\eps)\ d\tau\\
&=2\pi r^{-1} \int_{-t-r}^{-t+r}A_I(\rho,\theta;\eps)\ d\rho. }
In the first step, we notice that the integral is invariant under rotation. In the second step, we write the integral in the spherical coordinates. Moreover, recall from \eqref{prop:trzwtu:derivative} that
$2\partial_qU_I+A_I=O(\lra{q}^{-2+C\eps}e^{Cs})$,
so we have (recall that $r-t=-t^{2/3}$)
\fm{&\int_{\mathbb{S}^2}U_I(\eps\ln t-\delta,-t^{2/3},\omega;\eps)\ dS_\omega\\&=\int_{\mathbb{S}^2}\int_{-t-r}^{-t^{2/3}}\partial_qU_I(\eps\ln t-\delta,\rho,\theta;\eps)\ d\rho dS_\theta+O(t^{C\eps}\int_{-t-r}^{-t^{2/3}}\lra{\rho}^{-2+C\eps}\ d\rho+\eps^{-1}t^{-1/3+ C\eps})\\
&=\int_{\mathbb{S}^2}[U_I(\eps\ln t-\delta,-t^{2/3},\theta;\eps)-U_I(\eps\ln t-\delta,-2t+t^{2/3},\theta;\eps)]\ dS_\theta+O(\eps^{-1}t^{-1/3+ C\eps}).}
It follows that
\fm{\int_{\mathbb{S}^2} U_I(\eps\ln t-\delta,-2t+t^{2/3},\theta;\eps)\  dS_\theta=O(\eps^{-1}t^{-1/3+C\eps}),\qquad \forall t\geq 2e^{\delta/\eps}.}
Finally, by \eqref{trzwtu:esthat}, we have $\partial_qU_I=O(\lra{q}^{-1+C\eps}e^{Cs})$. By the mean value theorem, we have \fm{|U_I(\eps\ln t-\delta,-2t+t^{2/3},\theta;\eps)-U_I(\eps\ln t-\delta,-2t,\theta;\eps)|\lesssim t^{-1+C\eps}\cdot t^{2/3}\lesssim t^{-1/3+C\eps}.}
\end{proof}
\rmk{\rm Let us compare  \eqref{prop:decaysphmeans:c1} with the pointwise bounds for $U_I$ obtained in Proposition \ref{prop:trzwtu}. There it was proved that for all $(s,q,\omega)\in[0,\infty)\times\R\times\mathbb{S}^2$, we have $U_I=O(\eps^{-1}\lra{q}^{C\eps}e^{Cs})$. Thus, we have
\fm{|U_I(\eps\ln t-\delta,-2t,\omega;\eps)|\lesssim \eps^{-1}t^{C\eps},\qquad \forall t\geq 2e^{\delta/\eps},\omega\in\mathbb{S}^2.}
The estimate \eqref{prop:decaysphmeans:c1}, however, indicates that 
\fm{\inf_{\omega\in\mathbb{S}^2}|U_I(\eps\ln t-\delta,-2t,\omega;\eps)|\lesssim \eps^{-1}t^{-1/3+C\eps},\qquad\forall t\geq 2e^{\delta/\eps}.}
}
\rm

\subsection{Sufficient conditions for vanishing $\wh{A}$}\label{sec:prop:vanishingA}
We can now prove Theorem \ref{sdltthm}.
\prop{\label{prop:vanishingA}Suppose that the quasilinear wave equation \eqref{qwe} \emph{violates} the null condition. That is, we have $G(\omega)\not\equiv 0$ on $\mathbb{S}^2$, or equivalently, $g^{\alpha\beta}_0\neq0$ for some $\alpha,\beta=0,1,2,3$. Fix $(u_0,u_1)\in C_c^\infty(\R^3)$ and $0<\eps\ll1$. Let $u$ be the global $C^\infty$ solution to \eqref{qwe} with data $(\eps u_0,\eps u_1)$. Fix $\delta\in(0,1)$ and let $\wh{A}$ be the corresponding function defined by \eqref{prevac:defhata} in the asymptotic completeness result in Section \ref{secreview}.

Moreover, suppose that one of the following assumptions holds:
\begin{enumerate}[\rm a)]
    \item We have $\min_{(q,\omega)\in\mathbb{R}\times\mathbb{S}^2}\wh{A}(q,\omega;\eps)\geq 0$ or $\max_{(q,\omega)\in\mathbb{R}\times\mathbb{S}^2}\wh{A}(q,\omega;\eps)\leq 0$; 
    \item We have $\min_{ \omega\in \mathbb{S}^2}G(\omega)\geq 0$ and $\min\{\wh{A},0\}\in L^1_{q,\omega}(\R\times\mathbb{S}^2)$, or we have $\max_{ \omega\in \mathbb{S}^2}G(\omega)\leq 0$ and $\max\{\wh{A},0\}\in L^1_{q,\omega}(\R\times\mathbb{S}^2)$;
    \item Let $C_0>0$ be a constant uniform in all $\eps\ll1$ such that \fm{C_0\geq \sup_{(q,\omega)\in\R\times\mathbb{S}^2}|\frac{1}{2} G(\omega)\wh{A}(q,\omega;\eps)|.}
    Suppose that there exists a  constant $B_0>C_0$ that is also uniform in all $\eps\ll1$, such that  \fm{\lim_{q\to-\infty}\lra{q}^{1+B_0\eps}\sup_{\omega\in\mathbb{S}^2}|\wh{A}(q,\omega;\eps)|=0.}
\end{enumerate}
Then, as long as $\eps\ll1$ (depending on the constant $B_0$ in part c)),  we have $\wh{A}\equiv 0$ and thus $u\equiv 0$.}
\rm

\bigskip

To prove Proposition \ref{prop:vanishingA}, we will first prove that if one of the assumptions a), b), and c) holds, then we have $\wh{A}\equiv 0$. We will then show that $u\equiv 0$ if $\wh{A}\equiv 0$.

\subsubsection{The assumption a)}\footnote{\label{FN12}Note that the violation of the null condition is not used at all in this proof. Thus, the same conclusion holds no matter whether the null condition is satisfied or not.}\label{vanishingA:assumpa} Assume that $\min_{(q,\omega)\in\mathbb{R}\times\mathbb{S}^2}\wh{A}(q,\omega;\eps)\geq 0$ or $\max_{(q,\omega)\in\mathbb{R}\times\mathbb{S}^2}\wh{A}(q,\omega;\eps)\leq 0$. Now, we can rewrite \eqref{prop:decaysphmeans:c2} as
\fm{\int_{\mathbb{S}^2}\int_{\R}1_{q\in[-2t,R]}\cdot |\wh{A}(q,\omega;\eps)|\exp(\frac{1}{2}G(\omega)\wh{A}(q,\omega;\eps)\cdot (\eps\ln t-\delta))\ dqdS_\omega\lesssim \eps^{-1}t^{-1/3+C\eps},\  \forall t\geq 2e^{\delta/\eps}.}
Here the integrand is nonnegative everywhere. Besides, we have
\fm{\exp(\frac{1}{2}G(\omega)\wh{A}(q,\omega;\eps)\cdot (\eps\ln t-\delta))\geq \exp(-C(\eps\ln t-\delta))\gtrsim t^{-C\eps}.}
As a result, we have
\fm{\int_{\mathbb{S}^2}\int_{\R}1_{q\in[-2t,R]}\cdot |\wh{A}(q,\omega;\eps)|\ dqdS_\omega\lesssim \eps^{-1}t^{-1/3+C\eps}\cdot t^{C\eps},\qquad\forall t\geq 2e^{\delta/\eps}.}
We send $t\to\infty$ and apply the monotone convergence theorem. As long as $\eps\ll1$, we have $t^{-1/3+C\eps}\cdot t^{C\eps}\leq t^{-1/6}$. This yields
\fm{\int_{\mathbb{S}^2}\int_{\R} |\wh{A}(q,\omega;\eps)|\ dqdS_\omega=0}
and thus $\wh{A}=0$ almost everywhere. Since $\wh{A}$ is continuous, we conclude that $\wh{A}\equiv 0$.

\subsubsection{The assumption b)}\label{vanishingA:assumpb}
Assume that $\min_{ \omega\in \mathbb{S}^2}G(\omega)\geq 0$ and $\min\{\wh{A},0\}\in L^1_{q,\omega}(\R\times\mathbb{S}^2)$. Or, assume that $\max_{ \omega\in \mathbb{S}^2}G(\omega)\leq 0$ and $\max\{\wh{A},0\}\in L^1_{q,\omega}(\R\times\mathbb{S}^2)$. We now set
\fm{\mcl{O}_+&:=\{(q,\omega)\in\R\times\mathbb{S}^2:\ G(\omega)\wh{A}(q,\omega;\eps)>0\},\\
\mcl{O}_-&:=\{(q,\omega)\in\R\times\mathbb{S}^2:\ G(\omega)\wh{A}(q,\omega;\eps)<0\},\\
\mcl{O}_0&:=\{(q,\omega)\in\R\times\mathbb{S}^2:\ G(\omega)=0\}.}
Note that $\wh{A}=0$ in $(\R\times\mathbb{S}^2)\setminus(\mcl{O}_+\cup\mcl{O}_-\cup \mcl{O}_0)$ and that $\mcl{O}_+,\mcl{O}_-,\mcl{O}_0$ are pairwise disjoint. Thus, we can write the integral on the left hand side of \eqref{prop:decaysphmeans:c2} as the sum of
\fm{&\int_{\mcl{O}_+}\wh{A}(q,\omega;\eps)\cdot 1_{q\in[-2t,R]}\cdot\exp(\frac{1}{2}G(\omega)\wh{A}(q,\omega;\eps)\cdot(\eps\ln t-\delta))\ dqdS_\omega,\\
&\int_{\mcl{O}_-}\wh{A}(q,\omega;\eps)\cdot 1_{q\in[-2t,R]}\cdot\exp(\frac{1}{2}G(\omega)\wh{A}(q,\omega;\eps)\cdot(\eps\ln t-\delta))\ dqdS_\omega,\\
&\int_{\mcl{O}_0}\wh{A}(q,\omega;\eps)\cdot 1_{q\in[-2t,R]}\ dqdS_\omega.}
Since $G(\omega)\not\equiv 0$ on $\mathbb{S}^2$, the set $\{\omega\in\mathbb{S}^2:\ G(\omega)=0\}$ has surface measure $0$.  Thus, the last integral equals zero. 

If $\min_{\omega\in\mathbb{S}^2}G(\omega)\geq 0$, we have $\wh{A}|_{\mcl{O}_+}>0$ and $\wh{A}|_{\mcl{O}_-}<0$. If $\max_{\omega\in\mathbb{S}^2}G(\omega)\leq 0$, we have $\wh{A}|_{\mcl{O}_-}>0$ and $\wh{A}|_{\mcl{O}_+}<0$. Thus, by \eqref{prop:decaysphmeans:c2}, we have for all $t\geq 2t^{\delta/\eps}$
\fm{&|\int_{\mcl{O}_+}|\wh{A}|\cdot 1_{q\in[-2t,R]}\cdot\exp(\frac{1}{2}G(\omega)\wh{A} \cdot(\eps\ln t-\delta))\ dqdS_\omega\\
&\quad-\int_{\mcl{O}_-}|\wh{A}|\cdot 1_{q\in[-2t,R]}\cdot\exp(\frac{1}{2}G(\omega)\wh{A} \cdot(\eps\ln t-\delta))\ dqdS_\omega|\lesssim \eps^{-1}t^{-1/3+C\eps}.}
Under the assumption b), we have $\wh{A}\cdot 1_{\mcl{O}_-}=\wh{A}\cdot 1_{G(\omega)\wh{A}<0}\in L^{1}_{q,\omega}$. Moreover, we have $1\geq \exp(\frac{1}{2}G(\omega)\wh{A} \cdot(\eps\ln t-\delta))\to 0$ as $t\to\infty$ in $\mcl{O}_-$. Thus, by the Lebesgue dominated convergence theorem, we conclude that 
\fm{\int_{\mcl{O}_-}|\wh{A}|\cdot 1_{q\in[-2t,R]}\cdot\exp(\frac{1}{2}G(\omega)\wh{A} \cdot(\eps\ln t-\delta))\ dqdS_\omega\to 0,\qquad \text{as }t\to\infty.}As a result, we have
\fm{\int_{\mcl{O}_+}|\wh{A}|\cdot 1_{q\in[-2t,R]}\cdot\exp(\frac{1}{2}G(\omega)\wh{A} \cdot(\eps\ln t-\delta))\ dqdS_\omega\to 0,\qquad \text{as }t\to\infty.}
However, since $G(\omega)\wh{A}>0$ in $\mcl{O}_+$, by the monotone convergence theorem, we have
\fm{&\int_{\mcl{O}_+}|\wh{A}|\ dqdS_\omega=\lim_{t\to\infty}\int_{\mcl{O}_+}|\wh{A}|\cdot 1_{q\in[-2t,R]}\ dqdS_\omega\\
&\leq \lim_{t\to\infty}\int_{\mcl{O}_+}|\wh{A}|\cdot 1_{q\in[-2t,R]}\cdot\exp(\frac{1}{2}G(\omega)\wh{A} \cdot(\eps\ln t-\delta))\ dqdS_\omega=0.}
It follows that $\wh{A}=0$ almost everywhere in $\mcl{O}_+$. But the definition of $\mcl{O}_+$ implies that $\wh{A}\neq 0$ in $\mcl{O}_+$. Thus, $\mcl{O}_+$ has zero measure in $\R\times\mathbb{S}^2$. By \eqref{prop:decaysphmeans:c2} again, we have
\fm{&\int_{\mathbb{S}^2}\int_{\R}|\wh{A}|\cdot 1_{q\in[-2t,R]}\cdot\exp(\frac{1}{2}G(\omega)\wh{A} \cdot(\eps\ln t-\delta))\ dqdS_\omega\\
&=\int_{\mcl{O}_-}|\wh{A}|\cdot 1_{q\in[-2t,R]}\cdot\exp(\frac{1}{2}G(\omega)\wh{A} \cdot(\eps\ln t-\delta))\ dqdS_\omega\lesssim \eps^{-1}t^{-1/3+C\eps}.}
Following the proof in  Section \ref{vanishingA:assumpa}, we conclude that $\wh{A}\equiv 0$.

\subsubsection{The assumption c)} Let $C_0>0$ be a constant uniform in all $\eps\ll1$ such that \fm{C_0\geq \sup_{(q,\omega)\in\R\times\mathbb{S}^2}|\frac{1}{2} G(\omega)\wh{A}(q,\omega;\eps)|.}
Such a $C_0$ exists because we have $|\wh{A}|\lesssim\lra{q}^{-1+C\eps}$ with an implicit constant uniform in $\eps\ll1$. Suppose that there exists a  constant $B_0>C_0$ that is also uniform in all $\eps\ll1$, such that  \fm{\lim_{q\to-\infty}\lra{q}^{1+B_0\eps}\sup_{\omega\in\mathbb{S}^2}|\wh{A}(q,\omega;\eps)|=0.}
Recall that $\wh{A}$ is a continuous function and that $\wh{A}|_{q\geq R}\equiv 0$. 
From this limit, we have
\eq{\label{vanishingA:partc:Adecay}|\wh{A}(q,\omega;\eps)|&\lesssim_{B_0,\eps}\lra{q}^{-1-B_0\eps},\qquad \forall (q,\omega)\in\R\times\mathbb{S}^2.}
The implicit constant in \eqref{vanishingA:partc:Adecay} is allowed to depend on $\eps$. We remark that this implicit constant is not important in our proof for part c).

We now derive a pointwise decay for $u$ from \eqref{vanishingA:partc:Adecay}. Recall that \eqref{pfmthm:final} with $|I|=0$ yields
\fm{u(t,x;\eps)&=\frac{\eps }{2\pi}\int_{\mathbb{S}^2} \wh{A}(-t+x\cdot\theta,\theta;\eps)\ dS_\theta+O(\lra{t-|x|}^{-2}t^{C\eps})}
whenever $t\geq (2e^{\delta/\eps})^{1/3}$ and $|x|<t$.

\lem{\label{lem:secltthm:uptb}Fix $\nu_0\in(1/2,1)$. Then, for all $(t,x)\in\R^{1+3}$ with $t\geq (2e^{\delta/\eps})^{1/3}$ and $|x|\leq t-t^{\nu_0}/2$, we have
\fm{|u(t,x)|&\lesssim_{\eps,B_0,\nu_0}t^{-1}\lra{|x|-t}^{-B_0\eps}}
as long as $\eps\ll_{B_0,\nu_0}1$.}
\begin{proof}
We first estimate the spherical integral of $\wh{A}$ above.
By \eqref{vanishingA:partc:Adecay}, we have
\fm{\int_{\mathbb{S}^2}|\wh{A}(-t+x\cdot\theta,\theta;\eps)|\ dS_\theta&\lesssim_{B_0,\eps} \int_{\mathbb{S}^2}\lra{-t+x\cdot\theta}^{-1-B_0\eps}\ dS_\theta\lesssim |x|^{-1}\int_{-t-|x|}^{-t+|x|}\lra{\rho}^{-1-B_0\eps}\ d\rho.}
To obtain the last estimate, we express the spherical integral in the spherical coordinates. A similar computation has been done in \eqref{decaysphmeans:sphint}. If $|x|\leq t/2$, we use the mean value theorem to show that the right side is $O(\lra{|x|-t}^{-1-B_0\eps})=O(t^{-1-B_0\eps})$. If $t/2\leq |x|\leq t-t^{\nu_0}/2$, the right side is $O(t^{-1}\cdot(B_0\eps)^{-1}\lra{|x|-t}^{-B_0\eps})$. In summary, we have
\fm{\int_{\mathbb{S}^2}|\wh{A}(-t+x\cdot\theta,\theta;\eps)|\ dS_\theta&\lesssim_{B_0,\eps} t^{-1}\lra{|x|-t}^{-B_0\eps}.}
This estimate does not involve the parameter $\nu_0$.

It remains to estimate the remainder term $O(\lra{r-t}^{-2}t^{C\eps})$. For all $t\geq (2e^{\delta/\eps})^{1/3}$ and $|x|\leq t-t^{\nu_0}/2$, we have $\lra{r-t}\geq t^{\nu_0}/2$ and thus 
\fm{\lra{r-t}^{-2}t^{C\eps}\leq\lra{r-t}^{-B_0\eps}t^{C\eps}\cdot (t^{\nu_0}/2)^{-2+B_0\eps}\lesssim \lra{r-t}^{-B_0\eps}t^{-2\nu_0+(B_0+C)\eps}\leq \lra{r-t}^{-B_0\eps}t^{-1}.}
The last estimate holds because $2\nu_0>1$. We thus have $-2\nu_0+(B_0+C)\eps\leq -1$ as long as $\eps\ll_{B_0,\nu_0}1$. This finishes our proof.
\end{proof}\rm

In addition to \eqref{pfmthm:final}, we have a different way to approximate $u$. That is, we have \eqref{prevac:uapp} in Section \ref{secreview:secapp}. We also recall that  $\wt{u}=\wh{u}=\eps r^{-1}\wh{U}$. The bound in Lemma \ref{lem:secltthm:uptb} then implies a bound for $\wh{U}$.
\lem{\label{lem:est:wtu3}Fix $\nu_0\in(1/2,1)$. As long as $\eps\ll_{B_0,\nu_0}1$, we have
\fm{|\wh{U}(s,q,\omega;\eps)|\lesssim_{\eps,B_0,\nu_0}\lra{q}^{-B_0\eps},\qquad \forall s\geq \eps\ln 2,\ \omega\in\mathbb{S}^2,\ -qe^{-\frac{\nu_0(s+\delta)}{\eps}}\in[1/2,2]. }
The last condition for $q$ is equivalent to $q\in[-2t^{\nu_0},-t^{\nu_0}/2]$ with $t=e^{(s+\delta)/\eps}$. }
\begin{proof}By \eqref{prevac:uapp}, we have
\fm{&|u(t,x;\eps)-\eps |x|^{-1}\wh{U}(\eps\ln t-\delta,\wh{q}(t,x;\eps),x/|x|;\eps)|\lesssim_{\nu_0} \eps t^{-2+C\eps}\lra{t-|x|}\\
&\lesssim_{\nu_0} \eps t^{-2+(1+B_0\eps)\nu+C\eps}\lra{t-|x|}^{-B_0\eps}\lesssim_{\nu} \eps t^{-1}\lra{t-|x|}^{-B_0\eps}}
for all $(t,x)\in\Omega\cap\{|r-t|\leq 2t^{\nu_0}\}$. Since $\nu_0<1$, we have $-2+(1+B_0\eps)\nu_0+C\eps<-1$ for $\eps\ll_{\nu_0,B_0}1$.
By Lemma \ref{lem:secltthm:uptb}, for all $(t,x)\in\Omega\cap\{-2t^{\nu_0}\leq r-t\leq -t^{\nu_0}/2\}$, as long as $\eps\ll_{B_0,\nu_0}1$ we have
\fm{&|\wh{U}(\eps\ln t-\delta,\wh{q}(t,x;\eps),x/|x|;\eps)|\lesssim_{\nu_0}  |\eps^{-1}ru|+\eps^{-1}r\cdot  \eps t^{-1}\lra{t-|x|}^{-B_0\eps}\lesssim_{\eps,B_0,\nu_0} \lra{t-|x|}^{-B_0\eps}.}
Moreover, in the proof of Lemma \ref{prevac:lem:ahatqrt}, we have shown that for all $(t,x)\in\Omega\cap\{-2t^{\nu_0}\leq |x|-t\leq -t^{\nu_0}/2\}$,
\fm{&|\wh{U}(\eps\ln t-\delta,\wh{q}(t,x;\eps),x/|x|;\eps)-\wh{U}(\eps\ln t-\delta,|x|-t,x/|x|;\eps)|\\
&\lesssim (1+\ln\lra{|x|-t})\lra{|x|-t}^{-1}t^{C\eps}\lesssim\lra{t-|x|}^{-B_0\eps}\cdot \eps^{-1}\lra{|x|-t}^{-1+(B_0+C)\eps}t^{C\eps}\\
&\lesssim_{\eps,B_0,\nu_0}\lra{t-|x|}^{-B_0\eps}\cdot (t^{\nu_0}/2)^{-1+(B_0+C)\eps}t^{C\eps}\lesssim_{\eps,B_0,\nu_0}\lra{t-|x|}^{-B_0\eps} .} 
In the last step, we use $\nu_0>1/2$ and choose $\eps\ll_{\nu_0,B_0}1$. In summary, we have
\fm{&|\wh{U}(\eps\ln t-\delta,|x|-t,x/|x|;\eps)|\lesssim_{\eps,B_0,\nu_0} \lra{t-|x|}^{-B_0\eps}}
in $\Omega\cap\{-2t^{\nu_0}\leq |x|-t\leq -t^{\nu_0}/2\}$.

Finally, we notice that 
\fm{\{t\geq 2e^{\delta/\eps},-2t^{\nu_0}\leq r-t\leq -t^{\nu_0}/2\}\subset\Omega\cap\{-2t^{\nu_0}\leq r-t\leq -t^{\nu_0}/2\}.}
In fact, if $t\geq 2e^{\delta/\eps}$ and $-2t^{\nu_0}\leq r-t\leq -t^{\nu_0}/2$, then 
\fm{r-\frac{t+e^{\delta/\eps}}{2}-2R&=r-t+\frac{t-e^{\delta/\eps}}{2}-2R\geq \frac{t}{4}-(2t^{\nu_0}+2R)>0.}
Thus we have $(t,x)\in\Omega$ by the definition \eqref{secreview:defnOmega} of $\Omega$. This finishes the proof.
\end{proof}\rm

We can now show that $G(\omega)\wh{A}(q,\omega;\eps)=0$ everywhere. The value of $\nu_0$ in Lemmas \ref{lem:secltthm:uptb} and \ref{lem:est:wtu3} will be chosen in the proof of the next lemma.

\lem{\label{rmk:ga0}As long as $\eps\ll_{B_0}1$, we have $G(\omega)\wh{A}(q,\omega;\eps)=0$ for all $(q,\omega)\in\R\times\mathbb{S}^2$.}
\begin{proof}
Because of the finite speed of propagation, we only need to prove this lemma for $q\leq R$. There is nothing to prove when $G(\omega)=0$, so let us fix $\omega\in\mathbb{S}^2$ so that $G(\omega)\neq 0$.

Since $B_0>C_0$, we choose $\nu_0\in(1/2,1)$, depending only on $B_0$ and $C_0$ but not on $\eps$, such that $\nu_0B_0>C_0$. We now apply Lemmas \ref{lem:secltthm:uptb} and \ref{lem:est:wtu3} with this choice of $\nu_0$. This gives us a bound for $\wh{U}$; see Lemma \ref{lem:est:wtu3}.

We first show that $G\wh{A}\leq 0$ everywhere.  Set
\fm{I_+:=\{q\leq R:\ G(\omega)\wh{A}(q,\omega;\eps)>0\},\qquad I_-:=\{q\leq R:\ G(\omega)\wh{A}(q,\omega;\eps)< 0\}.}
Note that $\wh{A}\equiv 0$ on $(-\infty,R]\setminus(I_+\cup I_-)$.
By Lemma \ref{lem:est:wtu3} (with $s=\eps\ln t-\delta$ and $q=-t^{\nu_0}$), we have 
\eq{\label{lem:ga0est1}\abs{\int_{-t^{\nu_0}}^{R}\wh{A}(q,\omega;\eps)\exp(\frac{1}{2}G(\omega)\wh{A}(q,\omega;\eps)s)\ dq}=|\wt{U}(s,-t^{\nu_0},\omega;\eps)|\lesssim_{\eps,B_0,\nu_0} t^{-\nu_0 B_0\eps},\quad \forall t\geq 2e^{\delta/\eps}.}
The left side here can be written as $||\mcl{I}_+|-|\mcl{I}_-||$, where
\fm{\mcl{I}_\pm=\mcl{I}_\pm(t,\omega)=\int_{I_\pm\cap[-t^{\nu_0},R]}\wh{A}(q,\omega;\eps)\exp(\frac{1}{2}G(\omega)\wh{A}(q,\omega;\eps)s)\ dq.}
To see this, we notice that the sign of $\wh{A}$ remains unchanged in $I_+$ (or $I_-$). In other words, we have
\fm{\lim_{t\to\infty}||\mcl{I}_+|-|\mcl{I}_-||=0.}
In addition, since $G\wh{A}<0$ in $I_-$, we have
$\lim_{t\to\infty}\wh{A}\exp(\frac{1}{2}G\wh{A}s)=0$ whenever $q\in I_-$. Moreover, by \eqref{vanishingA:partc:Adecay}, we have
\fm{\int_\R |\wh{A}(q,\omega;\eps)|\ dq\lesssim_{\eps,B_0} \int_{-\infty}^{R}\lra{q}^{-1-B_0\eps}\ dq\lesssim_{\eps,B_0}1.}
As a result, we have
\fm{|1_{q\in I_-\cap[-t^{\nu_0},R]}\cdot \wh{A}(q,\omega;\eps)\exp(\frac{1}{2}G(\omega)\wh{A}(q,\omega;\eps)s)|\leq |\wh{A}(q,\omega;\eps)|\in L^1_q(\R).}By the Lebesgue dominated convergence theorem, we conclude that
$\lim_{t\to\infty}\mcl{I}_-=0$. This also forces $\lim_{t\to\infty}\mcl{I}_+=0$. However, since $ G\wh{A}s>0$ and since the sign of $\wh{A}(\cdot,\omega;\eps)$ remains unchanged on $I_+$, we have
\fm{\int_{I_+\cap[-t^{\nu_0},R]}|\wh{A}(q,\omega;\eps)|\ dq\leq |\mcl{I}_+|.}
By sending $t\to\infty$ on both sides and applying the monotone convergence theorem, we conclude that
\fm{\int_{I_+}|\wh{A}(q,\omega;\eps)|\ dq\leq \lim_{t\to\infty}|\mcl{I}_+|=0\Longrightarrow \wh{A}|_{I_+}\equiv 0.} This forces $I_+=\varnothing$. Thus, $G\wh{A}\leq 0$ everywhere.

Now, since $\mcl{I}_+=0$, by \eqref{lem:ga0est1} we have $|\mcl{I}_-|\lesssim_{\eps,B_0}  t^{-\nu_0 B_0\eps}$. Since \fm{C_0\geq \sup_{(q,\omega)\in\R\times\mathbb{S}^2}|\frac{1}{2}G(\omega)\wh{A}(q,\omega;\eps)|,} we have
\fm{\int_{I_-\cap[-t^{\nu_0},R]}|\wh{A}(q,\omega;\eps)|\ dq&\leq |\mcl{I}_-|\cdot\sup_{(q,\omega)\in\R\times\mathbb{S}^2}\exp(-\frac{1}{2}G(\omega)\wh{A}(q,\omega;\eps)s)\\
&\lesssim_{\eps,B_0} t^{-\nu_0 B_0\eps}\cdot \exp(C_0(\eps\ln t-\delta))\lesssim_{\eps,B_0} t^{(-\nu_0B_0+C_0)\eps}.}
The choice of $\nu_0$ implies that $\nu_0B_0>C_0$, so we have $\int_{I_-}|\wh{A}|\ dq=0$ by sending $t\to\infty$. Thus, we have $\wh{A}|_{I_-}\equiv 0$.
\end{proof}
\rm

Finally, we show that if $G(\omega)\wh{A}(q,\omega;\eps)=0$ everywhere and if $G(\omega)\not\equiv 0$ on $\mathbb{S}^2$, then $\wh{A}\equiv 0$. In fact, in Section \ref{vanishingA:assumpb}, we have mentioned that the set $\{\omega\in\mathbb{S}^2:\ G(\omega)=0\}$ has surface measure zero. In other words, we have $\wh{A}=0$ almost everywhere on $\mathbb{S}^2$. We thus conclude that $\wh{A}\equiv 0$ by its continuity.

\subsubsection{The final step}
To end the proof of Proposition \ref{prop:vanishingA}, we show that if $\wh{A}\equiv 0$, then $u\equiv 0$.

Suppose that $\wh{A}\equiv 0$. We first claim that $\norm{\partial u(t)}_{L^2(\R^3)}\lesssim \eps t^{-1/2+C\eps}$ for all $t\geq 2e^{\delta/\eps}$. Here we want to apply Lemma \ref{prevac:estmore:lem} instead of \eqref{prevac:uapp}, because Lemma \ref{prevac:estmore:lem} holds in the whole region $\Omega$. By \eqref{prevac:defhata}, we have
\fm{A(q,\omega)&=\wh{A}(F(q,\omega),\omega)=0,\qquad\forall (q,\omega)\in\R\times\mathbb{S}^2.}
Recall that $A_1A_2=-2A$ and that $A_1<-1$. Thus, we also have $A_2\equiv 0$ and thus $\wt{U}\equiv 0$. Here we recall the definitions of $A,A_1,A_2,\wt{U}$ in \eqref{prevac:limita} and \eqref{prevac:defnwtmuU}. By Lemma \ref{prevac:estmore:lem}, we have\fm{\sum_{|I|\leq 1}|Z^Iu|\lesssim \eps t^{-2+C\eps}\lra{r-t},\qquad \forall (t,x)\in\Omega.}
Thus, we have $|\partial u|\lesssim \eps t^{-2+C\eps}$ for all $(t,x)\in\Omega$ by \eqref{c1l2.1f1}.
Moreover, whenever $t\geq 2e^{\delta/\eps}$ but $(t,x)\notin\Omega$, we have $|x|\leq \frac{1}{2}(t+e^{\delta/\eps})+2R\leq 3t/4+2R\leq 7t/8$. Since $Z^Iu=O(\eps t^{-1+C\eps})$ by \eqref{secreview:uptb}, we have
\fm{|\partial u|&\lesssim \lra{r-t}^{-1}|Zu|\lesssim t^{-1}\cdot \eps t^{-1+C\eps}\lesssim \eps t^{-2+C\eps},\qquad t\geq 2e^{\delta/\eps},\ (t,x)\not\in\Omega.}
In summary, we have $|\partial u|\lesssim \eps t^{-2+C\eps}$ for all $t\geq 2e^{\delta/\eps}$ and all $x\in\R^3$. Since $u|_{r-t\geq R}\equiv 0$, we conclude that
\fm{\norm{\partial u(t)}_{L^2(\R^3)}\lesssim \eps t^{-2+C\eps}\cdot|B(0,t+R)|^{1/2}\lesssim \eps t^{-1/2+C\eps},\qquad \forall t\geq 2e^{\delta/\eps}.}
This finishes the proof of our claim.

Now, we apply the standard energy estimate from, e.g.,\ \cite[Proposition I.2.1]{MR2455195}. Since $u$ is a global solution to \eqref{qwe}, 
for each $0\leq t_1<t_2$, we have
\eq{\label{vanishingA:finalstep:energyest}\norm{\partial u(t_1)}_{L^2}\lesssim \norm{\partial u(t_2)}_{L^2}\exp(\int_{t_1}^{t_2}\norm{\partial (g^{**}(u))(\tau)}_{L^\infty}\ d\tau).}
Note that there is no inhomogeneous term in \eqref{qwe}. Moreover, the energy estimate in \cite[Proposition I.2.1]{MR2455195} is forward in time while \eqref{vanishingA:finalstep:energyest} is backward in time. To handle this difference, we simply apply the forward energy estimate to $w(t,x):=u(t_2-t,-x)$ which still solves \eqref{qwe}.

By the chain rule and the bound $|\partial u|\lesssim \eps \lra{t}^{-1}$ for all $t\geq 0$, we have 
\fm{\exp(\int_{t_1}^{t_2}\norm{\partial (g^{**}(u))(\tau)}_{L^\infty}\ d\tau)\lesssim\exp(\int_{t_1}^{t_2}C\eps (1+\tau)^{-1}\ d\tau)\lesssim (\frac{t_2+1}{t_1+1})^{C\eps}.}
And since $\norm{\partial u(t)}_{L^2}\lesssim \eps t^{-1/2+C\eps}$ for all $t\geq 2e^{\delta/\eps}$, we have
\fm{\norm{\partial u(t)}_{L^2}\lesssim \norm{\partial u(T)}_{L^2}\cdot (\frac{1+T}{1+t})^{C\eps}\lesssim \eps T^{-1/2+C\eps},\qquad \forall 0\leq t\leq T,\ T>2e^{\delta/\eps}.}
The left side is independent of $T$. Thus, by choosing $\eps\ll1$ and sending $T\to\infty$, we conclude that $\partial u(t)\equiv 0$ for each $t\geq 0$. Since $u(t,\cdot)$ has a compact support, we have  $u\equiv 0$ for all $(t,x)\in[0,\infty)\times\R^3$.

\subsection{Corollaries of Proposition \ref{prop:vanishingA}}\label{sec:pfltthm}

We finally prove Theorem \ref{ltthm_weak}. We will prove a slightly stronger version of this theorem.

\prop{\label{ltthm}Suppose that the quasilinear wave equation \eqref{qwe} \emph{violates} the null condition. Fix $(u_0,u_1)\in C_c^\infty(\R^3)$ and $0<\eps\ll1$. Let $u$ be the $C^\infty$ global solution to \eqref{qwe} with initial data $(\eps u_0,\eps u_1)$. Also fix $\delta\in(0,1)$. Then, there exist constants $B_0,B_1>0$ depending on $u_0,u_1,\delta$ (and not on $\eps$), such that the following conclusion holds. Suppose that one of the following assumptions holds:
\begin{enumerate}[\rm i)]
    \item Fix  $0<\nu_0<1$. For all $t\geq e^{\delta/\eps}$ and $\omega\in\mathbb{S}^2$, we have \eq{\label{ltthm:asu1}|(u_t-u_r)(t,(t-t^{\nu_0})\omega)|\lesssim \eps t^{-1}\lra{t^{\nu_0}}^{-1- B_0\eps}\sim \eps t^{-1-\nu_0-\nu_0B_0\eps};}
    \item  Fix  $0<\nu_0<1$. For all $t\geq e^{\delta/\eps}$ and all $x\in\R^3$ with  $|x|\in[t-2t^{\nu_0},t-t^{\nu_0}/2]$, we have   \eq{\label{ltthm:asu2} |u(t,x) |\lesssim\eps t^{-1-(\nu_0 B_0+B_1)\eps};}
    \item Fix  $0<\nu_0<1$. For all $t\geq e^{\delta/\eps}$ and $\omega\in\mathbb{S}^2$, we have \eq{\label{ltthm:asu30}\sum_{1\leq i<j\leq 3}|(\partial_t-\partial_r)\Omega_{ij}u(t,(t-t^{\nu_0})\omega)|\lesssim \eps t^{-1}\lra{t^{\nu_0}}^{-1- B_0\eps}\sim \eps t^{-1-\nu_0-\nu_0B_0\eps}.} Moreover, we have \eq{\label{ltthm:asu3}|u(t,0)|\lesssim\eps t^{-1-B_0\eps} \qquad \forall t\geq e^{\delta/\eps}.}
\end{enumerate}

Then, as long as $\eps\ll1$ which may depend on the constants in each of the assumptions, we have  $u
\equiv0$.}\\

\rm

In its proof, we will apply part c) of Proposition \ref{prop:vanishingA}. Let $C_0>0$ be the constant in part c) of Proposition \ref{prop:vanishingA}. Let  the $B_0$ in the statement of Proposition \ref{ltthm} be an arbitrary constant larger than $C_0$ (e.g. $B_0=C_0+1$). Following Section \ref{secreview}, we obtain a corresponding function $\wh{A}$.  It suffices to prove that 
\fm{ |\wh{A}(q,\omega;\eps)|\lesssim_{\eps,B_0} \lra{q}^{-1-B_0\eps},\qquad \forall (q,\omega)\in\R\times\mathbb{S}^2.}

\subsubsection{The assumption i)}

\lem{\label{ltthm:asu1:lem}Fix  $0<\nu_0<1$. Suppose that  \fm{|(u_t-u_r)(t,(t-t^{\nu_0})\omega)|\lesssim \eps t^{-1-\nu_0-\nu_0B_0\eps}\lra{r(t,\omega)-t}^{-1- B_0\eps},\qquad \forall t\geq e^{\delta/\eps},\ \omega\in\mathbb{S}^2.}
Then, we have
\fm{|\wh{A}(q,\omega;\eps)|&\lesssim_{B_0,\nu_0}\lra{q}^{-1-B_0\eps},\qquad \forall (q,\omega)\in\R\times\mathbb{S}^2.}}
\begin{proof}
It is easy to check that $(t,(t-t^{\nu_0})\omega)\in\Omega\cap\{|r-t|\lesssim t^{\nu_0}\}$ whenever $t\geq e^{\delta/\eps}$. Thus, by \eqref{prevac:uapp}, we have
\fm{|\partial (u-\wh{u})(t,(t-t^{\nu_0})\omega;\eps)|&\lesssim \lra{t^{\nu_0}}^{-1}|Z (u-\wh{u})(t,(t-t^{\nu_0})\omega;\eps)|\\
&\lesssim_{\nu_0}\eps t^{-2+C\eps}\lesssim \eps t^{-1-\nu_0-\nu_0B_0\eps}. }
In the last estimate, we use $-2+C\eps<-1-\nu_0-\nu_0B_0\eps$ as long as $\eps\ll_{\nu_0,B_0}1$. We thus have
\fm{|(\partial_t-\partial_r)  \wh{u} (t,(t-t^{\nu_0})\omega;\eps)|&\lesssim  \eps t^{-1-\nu_0-\nu_0B_0\eps}. }
By \eqref{prop:trzwtu:c:est0} in Proposition \ref{prop:trzwtu}, we also have
\fm{&(\partial_t-\partial_r)  \wh{u} (t,(t-t^{\nu_0})\omega;\eps)\\
&=-2\eps (t-t^{\nu_0})^{-1}\wh{A}(-t^{\nu_0},\omega;\eps)+O(\eps(1+\ln\lra{t^{\nu_0}})\lra{t^{\nu_0}}^{-2}t^{-1+C\eps}+\eps t^{-2+C\eps})\\
&=-2\eps (t-t^{\nu_0})^{-1}\wh{A}(-t^{\nu_0},\omega;\eps)+O(\eps t^{-1-\nu_0-\nu_0B_0\eps}).}
Previously, we have proved that $\eps t^{-2+C\eps}\lesssim \eps t^{-1}\lra{r(t,\omega)-t}^{-1-B_0\eps}$. Besides, we have $\eps(1+\ln\lra{t^{\nu_0}})\lra{t^{\nu_0}}^{-2}t^{-1+C\eps}\lesssim_{\nu_0} \eps t^{\nu_0/2-2\nu_0-1+C\eps}\leq \eps t^{-1-\nu_0-\nu_0B_0\eps}$. In summary, we have
\fm{|\wh{A}(-t^{\nu_0},\omega;\eps)|\lesssim \eps^{-1}(t-t^{\nu_0})\cdot \eps t^{-1-\nu_0-\nu_0B_0\eps}\lesssim \lra{t^{\nu_0}}^{-1-B_0\eps}. }
This estimate holds for all $t\geq e^{\delta/\eps}$ and $\omega\in\mathbb{S}^2$. Thus, we have
\fm{|\wh{A}(q,\omega;\eps)|\lesssim \lra{q}^{-1-B_0\eps},\qquad \forall q\leq -e^{\nu_0\delta/\eps},\omega\in\mathbb{S}^2.}
Recall that $\wh{A}|_{q\leq R}\equiv 0$. Besides, whenever $q\in[-e^{\nu_0\delta/\eps},R]$, we have \fm{|\wh{A}|\lesssim \lra{q}^{-1+C\eps}\lesssim \lra{q}^{-1-B_0\eps}\cdot \lra{e^{\nu_0\delta/\eps}}^{(B_0+C)\eps}\lesssim_{B_0}\lra{q}^{-1-B_0\eps}.}
This finishes the proof.
\end{proof}\rm

\subsubsection{The assumption ii)}
Assuming \eqref{ltthm:asu2}, we will prove that \eqref{ltthm:asu1} holds. Let $B_1>1$ be a constant to be chosen. By taking $q=-t^{\nu_0}$ and $q=-t^{\nu_0}-t^{\nu_0-B_1\eps}$ in \eqref{ltthm:asu2}, we obtain
\fm{|u(t,(t-t^{\nu_0})\omega)-u(t,(t-t^{\nu_0}-t^{\nu_0-B_1\eps})\omega)|&\lesssim \eps t^{-1-(\nu_0B_0+B_1)\eps}=\eps t^{-1-\nu_0-\nu_0B_0\eps}\cdot t^{\nu_0-B_1\eps}.}
By the mean value theorem, there exists $b_0\in(0,1)$ such that
\fm{|u_r(t,(t-t^{\nu_0}-b_0t^{\nu_0-B_1\eps})\omega)|&\lesssim \eps t^{-1-\nu_0-\nu_0B_0\eps}.}
Moreover, we notice that
\fm{|u_{rr}(t,x)|&\lesssim |\partial^2u|+r^{-1}|\partial u|\lesssim \eps t^{-1+C\eps}\lra{|x|-t}^{-2},\qquad\forall t\geq 1,\ |x|/t\in[1/2,2].}
It follows that
\fm{&|u_r(t,(t-t^{\nu_0})\omega)|\lesssim |u_r(t,(t-t^{\nu_0}-b_0t^{\nu_0-B_1\eps})\omega)|+\int_{t^{\nu_0}}^{t^{\nu_0}+b_0t^{\nu_0-B_1\eps}}|u_{rr}(t,(t-\rho)\omega)|\ d\rho\\
&\lesssim \eps t^{-1-\nu_0-\nu_0B_0\eps}+\int_{t^{\nu_0}}^{t^{\nu_0}+b_0t^{\nu_0-B_1\eps}}\eps t^{-1+C\eps}\lra{\rho}^{-2}\ d\rho\\
&\lesssim \eps t^{-1-\nu_0-\nu_0B_0\eps}+\eps t^{-1+C\eps}\cdot t^{\nu_0-B_1\eps}\cdot t^{-2\nu_0}\lesssim \eps t^{-1-\nu_0-\nu_0B_0\eps}+\eps t^{-1-\nu_0-(B_1-C)\eps}.}
To estimate the integral in the second row, we notice that $\lra{\rho}^{-2}\sim t^{-2\nu_0}$. Now, by choosing $B_1$ such that $B_1\geq B_0+C$, we conclude that $u_r(t,(t-t^{\nu_0})\omega)=O(\eps t^{-1-\nu_0-\nu_0B_0\eps})$. And since $u_t+u_r=O(\eps t^{-2+C\eps})$, we obtain \eqref{ltthm:asu1}.

\subsubsection{The assumption iii)}
\lem{\label{ltthm:asu3:lem}Fix  $0<\nu_0<1$. Suppose that  \fm{\sum_{1\leq i<j\leq 3}|(\partial_t-\partial_r)\Omega_{ij}u(t,(t-t^{\nu_0})\omega)|\lesssim \eps t^{-1}\lra{t^{\nu_0}}^{-1- B_0\eps}\sim \eps t^{-1-\nu_0-\nu_0B_0\eps},\qquad \forall t\geq e^{\delta/\eps},\omega\in\mathbb{S}^2.} 
Then, we have
\fm{\sum_{1\leq i<j\leq 3}|(\omega_i\partial_{\omega_j}-\omega_j\partial_{\omega_i})\wh{A}(q,\omega;\eps)|&\lesssim_{B_0,\nu_0}\lra{q}^{-1-B_0\eps},\qquad \forall (q,\omega)\in\R\times\mathbb{S}^2.}}
\begin{proof}
By following the proof in Lemma \ref{ltthm:asu1:lem}, we first obtain
\fm{\sum_{1\leq i<j\leq 3}|(\partial_t-\partial_r) \Omega_{ij} \wh{u} (t,(t-t^{\nu_0})\omega;\eps)|&\lesssim  \eps t^{-1-\nu_0-\nu_0B_0\eps}. }
If we chose a multiindex $I$ such that $Z^I=\Omega_{ij}$, then by \eqref{prop:trzwtu:inductformula:a}, we have $A_I=-2(\omega_i\partial_{\omega_j}-\omega_j\partial_{\omega_i})\wh{A}$.
By \eqref{prop:trzwtu:c:est0} in Proposition \ref{prop:trzwtu}, we thus have
\fm{&(\partial_t-\partial_r) \Omega_{ij} \wh{u} (t,(t-t^{\nu_0})\omega;\eps)\\
&=-2\eps (t-t^{\nu_0})^{-1}(\omega_i\partial_{\omega_j}-\omega_j\partial_{\omega_i})\wh{A}(-t^{\nu_0},\omega;\eps)+O(\eps(1+\ln\lra{t^{\nu_0}})\lra{t^{\nu_0}}^{-2}t^{-1+C\eps}+\eps t^{-2+C\eps})\\
&=-2\eps (t-t^{\nu_0})^{-1}(\omega_i\partial_{\omega_j}-\omega_j\partial_{\omega_i})\wh{A}(-t^{\nu_0},\omega;\eps)+O(\eps t^{-1-\nu_0-\nu_0B_0\eps}).}
The remaining proof is the same as that of Lemma \ref{ltthm:asu1:lem}.
\end{proof}\rm

By \eqref{ltthm:asu3} and \eqref{pfmthm:final} with $|I|=0$, we have
\fm{|\int_{\mathbb{S}^2}\wh{A}(-t,\theta;\eps)\ dS_\theta|&\lesssim \eps^{-1}t^{-2+C\eps}+t^{-1-B_0\eps}\lesssim t^{-1-B_0\eps}\qquad \forall t\geq e^{\delta/\eps}.}
Thus, for all $t\geq e^{\delta/\eps}$ and $\omega\in\mathbb{S}^2$, we have
\fm{|\wh{A}(-t,\omega;\eps)|&\leq|\wh{A}(-t,\omega;\eps)-\frac{1}{4\pi}\int_{\mathbb{S}^2}\wh{A}(-t,\theta;\eps)\ dS_\theta|+\frac{1}{4\pi}|\int_{\mathbb{S}^2}\wh{A}(-t,\theta;\eps)\ dS_\theta|\\
&\lesssim \int_{\mathbb{S}^2}|\wh{A}(-t,\omega;\eps)-\wh{A}(-t,\theta;\eps)|\ dS_\theta+t^{-1-B_0\eps}\lesssim_{B_0,\nu_0} t^{-1-B_0\eps}.}
In the last estimate, we apply Lemma \ref{ltthm:asu3:lem}. It follows that
\fm{|\wh{A}(q,\omega;\eps)|&\lesssim \lra{q}^{-1-B_0\eps},\qquad q\leq -e^{\delta/\eps},\omega\in\mathbb{S}^2.}
Using $\wh{A}=O(\lra{q}^{-1+C\eps})$ and $\wh{A}|_{q\geq R}\equiv 0$, we can once again prove $|\wh{A}(q,\omega;\eps)|\lesssim_{B_0} \lra{q}^{-1-B_0\eps}$ for all $(q,\omega)$.

\bibliography{paperltt}{}
\bibliographystyle{plain}

\end{document}